\documentclass[a4paper,twoside,tbtags]{article}
\usepackage{amsmath}\usepackage{amsfonts}\usepackage{amssymb}

\usepackage{amsthm}

\theoremstyle{definition}
\newtheorem{definition}{Definition}[section]

\newtheorem{theorem}{Theorem}[section]
\newtheorem{lemma}[theorem]{Lemma}

\newtheorem{proposition}[theorem]{Proposition}
\newtheorem{remark}{Remark}[section]

\usepackage{amsmath}\usepackage{amsfonts}\usepackage{amssymb} 
\usepackage{stackrel} \usepackage{latexsym}\usepackage{epsfig}\usepackage{oldgerm}

\newcommand\DN{\newcommand}

\DN\ws{weak solution}

\DN\lref[1]{Lemma~\ref{#1}}\DN\tref[1]{Theorem~\ref{#1}}\DN\pref[1]{Proposition~\ref{#1}}
\DN\sref[1]{Section~\ref{#1}}\DN\ssref[1]{Subsection~\ref{#1}}
\DN\dref[1]{Definition~\ref{#1}}\DN\rref[1]{Remark~\ref{#1}} \DN\corref[1]{Corollary~\ref{#1}}
\DN\eref[1]{Example~\ref{#1}}
\numberwithin{equation}{section}
\newcounter{Const} 
\DN\Ct{\refstepcounter{Const}c_{\theConst}}\DN\cref[1]{c_{\ref{#1}}}	

\DN\R{\mathbb{R}}\DN\N{\mathbb{N}}\DN\Q{\mathbb{Q}}\DN\C{\mathbb{C}}\DN\Z{\mathbb{Z}}

\DN\limi[1]{\lim_{#1\to\infty}} 	\DN\limz[1]{\lim_{#1\to0}}
\DN\limsupi[1]{\limsup_{#1\to\infty}} 	\DN\limsupz[1]{\limsup_{#1\to0}}
\DN\liminfi[1]{\liminf_{#1\to\infty}} 	\DN\liminfz[1]{\liminf_{#1\to0}}
\DN\sumii[1]{\sum_{#1=1}^{\infty}}\DN\sumi[1]{\sum_{#1=0}^{\infty}}

\DN\map[3]{#1\!:\!#2\!\to\!#3} \DN\ot{ \otimes } \DN\ts{ \times }
\DN\PD[2]{\frac{\partial#1}{\partial#2}}
\DN\half{\frac{1}{2}}
\DN\Rd{\R ^d}\DN\RdN{(\Rd )^{\N }}\DN\Rdm{\R ^{dm}}

\DN\bs{\bigskip}\DN\ms{\smallskip}

\DN\bigcupii[1]{\bigcup_{#1=1}^{\infty}}

\DN\As[1]{\textbf{(#1)}}

\DN\PF{\medskip \noindent {\em Proof. }} \DN\PFEND{\qed \medskip}

\DN\RA{\Rightarrow}

\DN\XRl{\mathbf{X}^{R,\La }}
\DN\Hrm{\mathcal{H}_r^m}
\DN\xx{\mathsf{x}} 
\DN\xxSR{\xx (\SR )}

\DN\labx{\lab (\xx )}
\DN\labi{\lab ^i}
\DN\labii{\lab ^{i+1}}
\DN\labix{\lab ^i(\xx )}
\DN\labiix{\lab ^{i+1}(\xx )}
\DN\Pls{{P}_{\lab (\xx )}^{\infty }}
\DN\Els{{E}_{\lab (\xx )}^{\infty }}

\DN\PP{P}
\DN\PPP{\mathbf{P}}
\DN\PPPP{\mathbb{P}}

\DN\PPx{\mathsf{P}_{\xx }}

\DN\Pbar{\bar{{P}}}
\DN\Pbars{\bar{{P}}_{\mathbf{s}}}
\DN\Ps{\P _{\mathbf{s}}}
\DN\PBr{P _{\mathrm{Br}}^{\infty}}
\DN\PbarsB{\Pbar _{\sB }} 
\DN\Psbf{\mathbf{P}_{\mathbf{s}}}

\DN\Pssfl{\mathsf{P}_{\xx }^{\La }}
\DN\PRssfl{\mathsf{P}_{R,\xx }^{\La }}
\DN\PRl{\mathsf{P}_{R}^{\La }}

\DN\Pml{\mathbf{P}_{\mu ^{\lab } } }
\DN\Pmul{\mathbf{P}_{\mul }} 

\DN\PPUone{\mathsf{P}^{\Os ,[1]}}
\DN\PPUm{\mathsf{P}^{\Os ,[m]}}
\DN\PPU{\mathsf{P}^{\Os }}
\DN\PPmU{\mathsf{P}_{\mu }^{\Os }}
\DN\PPxU{\mathsf{P}_{\xx }^{\Os }}
\DN\PPsU{\mathsf{P}_{\sss }^{\Os }}

\DN\Pmut{\mathsf{P}_{\mut }}

\DN\CS{C([0,\infty); S^{\N })}
\DN\CSm{C([0,\infty); S^{m})}
\DN\Xsm{\mathbf{X}^{\sss ,m}}
\DN\Xsn{\mathbf{X}^{\sss ,n}}
	\DN\XRi{X^{\Rsss ,i}}\DN\XRj{X^{\Rsss ,j}}
\DN\XRs{\mathbf{X}^{\Rsss }}
\DN\XRm{\mathbf{X}^{\Rsss ,m}}
\DN\XRmbar{\overline{\mathbf{X}}^{\Rsss ,m}}
\DN\Xsmbar{\overline{\mathbf{X}}^{\sss ,m}}
\DN\Xmbar{\overline{\mathbf{X}}^{m}}
\DN\XLa{\mathbf{X}}
\DN\XOs{\mathbf{X}}

\DN\XXXRit{(\XRit ,\XRidt )}
\DN\XXXRiu{(\XRiu ,\XRidu )}
\DN\XXXRiz{(\XRiz ,\XRidz )}
\DN\XXXsiu{(\Xsiu ,\Xsidu )}
\DN\XXXsit{(\Xsit ,\Xsidt )}
\DN\XXXIsit{(\XIsit ,\XIsidt )}
\DN\XXXsiz{(\Xsiz ,\Xsidz )}
\DN\XRit{\XRi _t}
\DN\XRjt{\XRj _t}
\DN\XRiu{\XRi _u}
\DN\XRju{\XRj _u}
\DN\XRiz{\XRi _0}
\DN\XXXit{(\Xit ,\Xidt )}
\DN\Xit{X_t^i} 

\DN\XRidt{\mathsf{X}_t^{\Rsss ,\diamond i}}
\DN\Xidt{\mathsf{X}_t^{\diamond i}}
\DN\Xsidt{\mathsf{X}_t^{\sss ,\diamond i}}
\DN\XIsidt{\mathsf{X}_t^{\infty\sss ,\diamond i}}
\DN\Xsidz{\mathsf{X}_0^{\sss ,\diamond i}}
\DN\XRidu{\mathsf{X}_u^{\Rsss ,\diamond i}}
\DN\XRidz{\mathsf{X}_0^{\Rsss ,\diamond i}}
\DN\Xidu{\mathsf{X}_u^{\sss ,\diamond i}}
\DN\Xsidu{\mathsf{X}_u^{\sss ,\diamond i}}

\DN\Xsit{\Xsi _t} 
\DN\XIsit{\XIsi _t} 
\DN\Xsjt{\Xsj _t} 
\DN\Xsiu{\Xsi _u} 
\DN\Xsju{\Xsj _u} 
\DN\Xsiz{\Xsi _0} 

\DN\Xsi{X^{\sss ,i}}
\DN\Xsid{X^{\sss ,i}}

\DN\XIsi{X^{\infty\sss ,i}}
\DN\Xsj{X^{\sss ,j}}

\DN\XR{\mathbf{X}^{R}}
\DN\BB{\mathsf{B}}

\DN\IrT{\mathit{I}_{r,T}}
\DN\EE{E}

\DN\nN{{R}}
\DN\FBL{F(\mathbf{X}^{\Rsss,m}, 
\mathbf{B}_{\rsp }^{\Rsss,m},
\mathbf{L}^{\Rsss,m})}

\DN\EDrlwr{(\Er ^{\La },\mathcal{D} _r^{\La } )}
\DN\Erlwr{\Er ^{\La }}
\DN\EDrupr{(\E ^{\Os },\mathcal{D} _r^{\Os } )}

\DN\0{\limi{r}\limi{s}\limi{\mathsf{p}}\limi{R}}
\DN\1{\Rsss ,\La } 
\DN\2{\Rsss } 
\DN\3{\muRs (d\mathsf{x}) } 
\DN\4{\mu (d\sss ) } 
\DN\5{1_{\Sr } (\XRit ) }
\DN\6{L^{1}(\sS \times \mathsf{S} ,\, \muone _r )} 
\DN\8{\sigma _a^{\Rsss,m}}
\DN\9{(\ER ^{\La},\mathcal{D} _R^{\La})} 

\DN\sigmaA{\sigma_a^{\sss,m}}

\DN\nc{\nabla _j \psi ( \mathbf{X}_u^{\Rsss , m } )}
\DN\phat{} 
\DN\QN{Q_{r,s,\mathsf{p}}^{\Rsss,m}}
\DN\RN{R_{r,s}^{\Rsss,m}}

\DN\muNl{\mu^{\Rsss}\circ \lab^{-1}}
\DN\mul{\mu \circ \lab ^{-1}}
\DN\muls{\mu ^{\sss }\circ \lab ^{-1}}
\DN\mus{\mu ^{\sss }}
\DN\muRs{\mu ^{\Rsss }}
\DN\muRsone{\mu^{\Rsss ,[1]}} 

\DN\XXX{\mathsf{X}}
\DN\diai{\diamond i}
\DN\diaj{\diamond j}
\DN\Xid{\XXX^{\diai }}

\DN\tail{\mathrm{tail}}

\DN\bN{\mathsf{b}^{\nN }} 
\DN\bNi{\mathsf{b}^{\nN ,i}}
\DN\bNrs{\mathsf{b}_{r,s}^{\nN }}
\DN\bNrsp{\mathsf{b}_{r,s,\p }^{\nN }}
\DN\bNrsq{\mathsf{b}_{r,s,\q }^{\nN }}
\DN\brsp{\mathsf{b}_{r,s,\p }}
\DN\brsq{\mathsf{b}_{r,s,\q }}
\DN\brs{\mathsf{b}_{r,s }}
\DN\br{\mathsf{b}_{r}}
\DN\bNrstail{\mathsf{b}_{r,s}^{\nN ,\tail }}
\DN\btail{\mathsf{b} ^{\tail }}
\DN\cNrsp{\mathsf{c}_{r,s,\p }^{\nN }}

\DN\bis{\mathsf{b}}
\DN\brsps{\brsp }
\DN\brsptail{\brsp ^{\tail }}
\DN\brss{\brs }
\DN\is{\infty \, \sss }
\DN\sm{\sss ,m}
\DN\si{\sss ,i}

\DN\pirc{\pi _r^c}
\DN\pisc{\pi _s^c}
\DN\pitc{\pi _t^c}
\DN\pir{\pi _r}
\DN\pis{\pi _s}

 \DN\NN{\mathcal{N}}

\DN\rgn{\rho _{\mathrm{Gin}}}
\DN\kg{\mathsf{K}_{\mathrm{Gin}}}
\DN\mug{\mu _{\mathrm{Gin}}}
\DN\mugx{\mu _{\mathrm{Gin},x}}
\DN\muSinb{\muone _{\mathrm{Sin}, \beta }}
\DN\muAi{\mu _{\mathrm{Ai}}}

\DN\mub{\mu _{\mathrm{Be},\alpha }}

\DN\muone{\mu ^{[1]}}
\DN\muN{\mu ^{\n }}
\DN\muNx{\mu ^{\n }_{x}}
\DN\muNzero{\mu ^{\n }_{0}}
\DN\muNone{\mu ^{\n ,{1}}}
\DN\muNonebar{\bar{\mu }^{\n ,{1}}_{\q }}
\DN\muNstar{\mu ^{\n *}}

\DN\rairybeta{\rho _{\mathrm{Ai},\beta }}
\DN\rairytwo{\rho _{\mathrm{Ai},2}}
\DN\rNone{\rho ^{\n ,1}}
\DN\rNtwo{\rho ^{\n ,2}}
\DN\rNnk{\rho ^{\n ,n+k}}

\DN\rhobar{\bar{\rho }}
\DN\rb{\rho _{\beta }^n}
\DN\rbone{\rho _{\beta }^1}
\DN\rbx{\rho _{1 }^n}
\DN\rby{\rho _{2 }^n}
\DN\rbz{\rho _{4 }^n}

\DN\rbNone{\rho _{\beta }^{\n ,1}}
\DN\rbNtwo{\rho _{\beta }^{\n ,2}}
\DN\rbNxone{\rho _{\beta ,x}^{\n ,1}}
\DN\rbNxtwo{\rho _{\beta ,x}^{\n ,2}}
\DN\rbNxn{\rho _{\beta ,x}^{\n ,n}}
\DN\rN{\rho ^{\n ,n}}
\DN\rbN{\rho _{\beta }^{\n ,n}}
\DN\rbNN{\rho _{\beta }^{\n ,N}}
\DN\rbNx{\rho _{1}^{\n ,n}}
\DN\rbNy{\rho _{2}^{\n ,n}}
\DN\rbNz{\rho _{4}^{\n ,n}}

\DN\rhohat{\hat{\rho}}

\DN\sS{S}

\DN\SSS{\mathbf{S}}
\DN\SSw{\mathsf{S}_{\mathrm{SDE}}}
\DN\SSSw{\mathbf{S}_{\mathrm{SDE}}}

\DN\SSr{\mathsf{S} _{r}}
\DN\SSrm{\mathsf{S} _r^m} 
\DN\SSrk{\mathsf{S} _r^k}
\DN\SSrNk{\mathsf{S} _{r}^{N-k}}
\DN\SSk{\mathsf{S} ^{[k]}}
\DN\SSksingle{\Ssi ^{k}}
\DN\Ssi{\mathsf{S} _{\mathrm{s.i.}}}
\DN\Ssip{\mathsf{S} _{\mathrm{s.i.}}^{+f}}
\DN\SSone{\mathsf{S} ^{\mathbf{1}}}
\DN\SoneSSS{\sS \ts \mathsf{S} }

\DN\Sk{\sS ^{k}}
\DN\Sm{\sS ^{m}}
\DN\Sp{\sS _{\q }}
\DN\Srk{\Sr ^{k}}
\DN\Srm{\Sr ^{m}}
\DN\Sr{\sS _{r}}
\DN\SR{\sS _{R}}
\DN\SN{\sS ^{\mathbb{N}}} 
\DN\Sz{\SSS_{0}} 
\DN\SSz{\mathsf{S} _{0}} 

\DN\HHa{\mathsf{H}_{\mathsf{a}}}
\DN\HHz{\mathsf{H}_{0}}

\DN\Akr{\mathsf{A}^{k}_{r}}
\DN\Akrr{\mathsf{A}^{k}_{r+1}}
\DN\AkrN{\mathsf{A}^{k}_{r,N}}

\DN\dlog{\mathsf{d}}
\DN\dmu{\dlog ^{\mu }}
\DN\dpsi{\dlog ^{\mupsi }}
\DN\dmuN{\dlog ^{\n }}
\DN\dmuone{\dlog ^{\mu ^{1}}}
\DN\dmuhat{\hat{\dlog } ^{\mu }}

\DN\dG{\dlog ^{\mug }}
\DN\dSine{\dlog ^{\mu _{\mathrm{sin},\beta }} }

\DN\sigmaXms{\sigma \xm }
\DN\aaa{\mathit{a}}
\DN\bb{\mathsf{b}}
\DN\bbrs{\bb }

 \DN\sss{\mathsf{s}}
 \DN\Rsss{R\sss }
 \DN\rspN{r,s,\mathsf{p} \in \N }
 \DN\rspRN{r,s,\mathsf{p},R \in \N }
 \DN\rsp{r,s,\mathsf{p}}
 \DN\mrspN{m,r,s,\mathsf{p} \in \N } 

\DN\CY{C^{\infty} (\Rd )\ot\dY } 

\DN\Damu{\mathcal{D}^{\amu }} 
\DN\di{\mathcal{D} _{\circ }}
\DN\dib{\mathcal{D} _{\circ \mathrm{b}}}
\DN\dimu{\di ^{\mu }}
\DN\dRcirc{\mathcal{D} _{R,\circ }}

\DN\dibone{C_{0}^{\infty} (\sS )\ot \dib }

\DN\Lmu{L^2(\mathsf{S} ,\mu )}
\DN\Lm{\Lmu } 
\DN\Lmuk{L^{2}(\muk )}

\DN\Llocp{L_{\mathrm{loc}}^{p}}
\DN\Llocq{L_{\mathrm{loc}}^{q}}
\DN\Lloctwo{L_{\mathrm{loc}}^{2}}
\DN\Llocone{L_{\mathrm{loc}}^{1}}

\DN\LlocG{\Lloctwo (\muone _{\mathrm{Gin}})}
\DN\LlocpG{\Llocp (\muone _{\mathrm{Gin}})}
\DN\LchiG{L^2 (\chin\muone _{\mathrm{Gin}})}
\DN\LlocSine{\Lloctwo (\muone _{\mathrm{sin},\beta })}
\DN\LchiSine{L ^2 (\chin \muone _{\mathrm{sin},\beta })}
\DN\LlocpSine{\Llocp (\muone _{\mathrm{sin},\beta })}

\DN\Emu{\mathcal{E}}
\DN\Ermu{\Er }
\DN\Ermum{\Er ^{m }}
\DN\dom{\mathcal{D} }
\DN\Capa{\mathrm{Cap}}

\DN\E{\mathcal{E} }
\DN\Ea{\mathcal{E}}
\DN\Eak{\mathcal{E}^{k}}
\DN\Eaone{\mathcal{E}^{1}}
\DN\Eazero{\mathcal{E}^{0}}

\DN\Er{\E _r} \DN\ER{\E _R}
\DN\Erhat{\hat{\E }_r} \DN\ERhat{\hat{\E }_R}

\DN\amu{\mathsf{a} ,\mu }
\DN\amuone{\mathsf{a} ,\muone }
\DN\DDD{\mathbb{D}}
\DN\DDDa{\mathbb{D}^{\mathsf{a} }}
\DN\DDDk{\mathbb{D}^{k}}
\DN\DDDzero{\mathbb{D}^{0}}
\DN\DDDr{\DDD _{r}}
\DN\DDDR{\DDD _{R}}

\DN\ulab{\mathfrak{u} }
\DN\lab{\mathfrak{l} } 

\DN\ulabm{\mathfrak{u} _m}
\DN\labm{\mathfrak{l} _m} 

\DN\ulabma{\mathfrak{u} _m^{\mathsf{a} }}
\DN\labma{\mathfrak{l} _m^{\mathsf{a} }} 

 \DN\kpath{\ulab _{\mathrm{path}}}
 \DN\upath{\ulab _{\mathrm{path}}}
 \DN\kkpath{\ulab _{\mathrm{k,path}}}
 \DN\lpath{\lab _{\mathrm{path}}} 
 \DN\klpath{\lab _{\mathrm{k,path}}}

\DN\muairy{\mu_{\mathrm{Ai}, 2}}
\DN\muairyone{\mu_{\mathrm{Ai}, 1}}
\DN\muairyfour{\mu_{\mathrm{Ai}, 4}}
\DN\muairybeta{\mu_{\mathrm{Ai}, \beta }}

\DN\muairyk{\mu_{\mathrm{Ai}, 2}^{[k]}}
\DN\muairyonek{\mu_{\mathrm{Ai}, 1}^{[k]}}
\DN\muairyfourk{\mu_{\mathrm{Ai}, 4}^{[k]}}
\DN\muairybetak{\muairybeta ^{[k]}}
\DN\muairybetax{\mu_{\mathrm{Ai}, \beta ,x}}

\DN\musin{\mu_{\mathrm{sin}, 2}}
\DN\musinone{\mu_{\mathrm{sin}, 1}}
\DN\musinfour{\mu_{\mathrm{sin}, 4}}
\DN\musinbeta{\mu_{\mathrm{sin}, \beta}}

 \DN\WS{C([0,\infty);\mathsf{S} )}

 \DN\WSN{C([0,\infty);\SN )}
 \DN\WSz{C([0,\infty);\Sz )}
 \DN\WH{C([0,\infty);\mathsf{H} )}
 \DN\WSsi{C([0,\infty);\Ssi )}
 \DN\WSNs{C ^{\mathbf{s}}([0,\infty);\SN )}
 \DN\WSNz{C ^{\mathbf{0}}([0,\infty);\SN )} 
\DN\mupsi{\mu _{\Psi }}

\DN\WW{\WTSN \times \WTSNz }

\DN\pq{\p ,\q }
\DN\pqk{\p ,\q ,\kk }
\DN\pqkl{\p ,\q ,\kk , l }

\DN\p{\mathsf{p}}\DN\q{\mathsf{q}}

\DN\ijn{_{i,j=1}^{n}}
\DN\Os{\mathsf{upr}}
\DN\La{\mathsf{lwr}}

\DN\dRN{\mathcal{D}_{\circ , R , N}^{\mu } }
\DN\dR{\mathcal{D}_{\circ , R }^{\mu } }
\DN\dmuq{\mathcal{D}_{\circ , R }^{\mu } }
\DN\7{the $ \Lmu $, $ \E _{R,1}^{\La}$, and $ \E _{R,1}^{\1 }$-norm for $ \mu $-a.s.\! $ \sss $}
\DN\seven{the $ \Lmu $, $ \E _{R',1}^{\La}$, and $ \E _{R',1}^{\1 }$-norm for $ \mu $-a.s.\! $ \sss $}

\begin{document}

\begin{center}\begin{Large} {\sf 
Uniqueness of Dirichlet forms related to \\
infinite systems of interacting Brownian motions
}\end{Large}\end{center} \bigskip 
\begin{center}
{\bf Yosuke Kawamoto$^{1}$, Hirofumi Osada$^{2}$, 
Hideki Tanemura$^{3}$}
\end{center}
 \pagestyle{myheadings}
 \markboth{{\tt Uniqueness of Dirichlet forms. \quad \today }
 \quad } 
{ \sf{\today: \quad Yosuke Kawamoto, Hirofumi Osada, Hideki Tanemura }}


\begin{abstract}
The Dirichlet forms related to various infinite systems of interacting Brownian motions are studied. 
For a given random point field $ \mu $, there exist two natural infinite-volume Dirichlet forms 
$ (\mathcal{E}^{\mathsf{upr}},\mathcal{D}^{\mathsf{upr}})$ and $(\mathcal{E}^{\mathsf{lwr}},\mathcal{D}^{\mathsf{lwr}})$ on 
$ L^2(\mathsf{S} ,\mu ) $ describing interacting Brownian motions each with unlabeled equilibrium state $ \mu $. 
The former is a decreasing limit of a scheme of such finite-volume Dirichlet forms, and the latter is an increasing limit of another scheme of such finite-volume Dirichlet forms. Furthermore, the latter is an extension of the former. 
We present a sufficient condition such that these two Dirichlet forms are the same. In the first main theorem (Theorem 3.1) the Markovian semi-group given by 
$(\mathcal{E}^{\mathsf{lwr}},\mathcal{D}^{\mathsf{lwr}})$ is associated with a natural infinite-dimensional stochastic differential equation (ISDE). In the second main theorem (Theorem 3.2), we prove that these Dirichlet forms coincide with each other by using the uniqueness of {\ws}s of ISDE. We apply Theorem 3.1 to stochastic dynamics arising from random matrix theory such as the sine, Bessel, and Ginibre interacting Brownian motions and interacting Brownian motions with Ruelle's class interaction potentials, and Theorem 3.2 to the sine$ _2$ interacting Brownian motion and interacting Brownian motions with Ruelle's class interaction potentials of $ C_0^3 $-class. 
\end{abstract}

{\small \tableofcontents }

\section{Introduction}\label{s:1} 

An infinite system of interacting Brownian motions in $ \Rd $ can be represented by an $ (\Rd )^{\N }$-valued stochastic process 
$\mathbf{X}=(X^i )_{i\in\mathbb{N}}$ \cite{lang.1,lang.2,o.dfa,o.rm}.
This process is realized using several probabilistic constructs such as stochastic differential equation, Dirichlet form theory, and martingale problems. 
Among them, the second author constructed processes in a general setting using the technique of Dirichlet forms \cite{o.dfa,o.rm}. 

Specifically, the Dirichlet form introduced, $(\E^{\Os},\mathcal{D} ^{\Os})$ is obtained by the smallest extension of the bilinear form $ (\Emu ,\dimu )$ on $ \Lm $ with domain $ \dimu $ defined by 
\begin{align} 
&\notag 
\Emu (f,g) = \int_{\mathsf{S} } \DDD [f,g](\sss) \, \mu (d\sss)
,\\ &\label{:11b}
\DDD [f,g] (\sss) = \frac{1}{2} 
\sum_{i=1}^{\infty} 
 \nabla_{s_i} \check{f} \cdot \nabla_{s_i} \check{g}
,\\&\notag 
\dimu = 
\{ f \in \di \cap \Lmu \, ;\, \Emu (f,f) < \infty \} 
,\end{align}
where $\di $ is the set of all local smooth functions on the (unlabeled) configuration space $\mathsf{S} $ introduced in \eqref{:20a},
$\check{f}$ is a symmetric function such that $\check{f}(s_1,s_2,\dots)=f(\sss)$, $\cdot$ is the inner product in $\Rd$, $ \mathsf{s}=\sum_i \delta_{s_i}$ 
denotes a configuration, and $ \mu $ is a random point field on $ \Rd $. 
If we take $ \mu $ to be the Poisson random point field, the intensity of which is the Lebesgue measure, then the diffusion given by the Dirichlet form $ (\E^{\Os},\mathcal{D} ^{\Os})$ is the unlabeled Brownian motion $ \mathsf{B}$ such that $ \mathsf{B}_t = \sum_{i=1}^{\infty} \delta_{B_t^i}$, where $ \{ B^i \}_{i=1}^{\infty} $ is a system of independent copies of the standard Brownian motion.

This Dirichlet form is a decreasing limit of Dirichlet forms associated with finite systems of interacting Brownian motions in bounded domains $ \SR = \{ x\in\Rd ; |x|\le R \} $ with a boundary condition. 
Because of the boundary condition, Brownian particles that touch the boundary disappear. 
Also, particles enter the domain from the boundary according to the reversible measure $\mu $.

In contrast, Lang constructed the infinite system of Brownian motions as a limit of stochastic dynamics in bounded domains $ \SR $ by considering finite systems with another boundary condition \cite{lang.1,lang.2}. 
In his finite systems, a particle hitting the boundary is reflected and hence the number of particles in the domain is invariant.
His process is associated with the Dirichlet form 
$ (\E^{\La},\mathcal{D} ^{\La})$ that is the increasing limit of the Dirichlet forms associated with finite systems with the reflecting boundary condition.

In this paper, we discuss the relation between 
these Dirichlet forms, 
$(\E^{\Os},\mathcal{D} ^{\Os})$ and $(\E^{\La},\mathcal{D} ^{\La})$. 
The main purpose of this paper is to give a sufficient condition for 
\begin{align}\label{:11e}&
(\E^{\La},\mathcal{D} ^{\La})=(\E^{\Os},\mathcal{D} ^{\Os}) 
.\end{align}
By construction the inequality 
\begin{align}\label{:11f}&
(\E^{\La},\mathcal{D} ^{\La})\le (\E^{\Os},\mathcal{D} ^{\Os}) 
\end{align}
always holds whereas \eqref{:11e} does not necessarily hold in general. 
Although the problem is quite natural and general, little is known about the equality \eqref{:11e}. 
To the best of our knowledge, the unique example for which the equality \eqref{:11e} holds is the system of 
hard-core Brownian balls 
proved by the third author \cite{T3}.

The study of infinite systems of interacting Brownian motions was initiated by Lang \cite{lang.1,lang.2} and continued by 
Fritz \cite{Fr}, the third author \cite{T2}, and others. 
In their respective work, the free potential $ \Phi $ is assumed to be zero and 
the interaction potentials $\Psi $ are of class $ C_0^3 (\Rd )$ 
or exponentially decay at infinity and satisfying the super-stability in the sense of Ruelle. 
The infinite-dimensional stochastic differential equation (ISDE) is given by 
\begin{align}& \label{:12a}
dX^i_t = dB^i_t - \frac{\beta}{2} \sum ^{\infty}_{j=1,\,j\ne i} 
\nabla \Psi (X ^i_t-X^j_t) dt 
. \end{align}
Here $ \beta > 0 $ is an inverse temperature. 
Lang constructed a solution for the $ \mu $-a.s.\! $ \xx $ 
unlabeled initial point, where $ \mu $ is a grand canonical Gibbs measure with interaction potential $ \Psi $. 

Indeed, Lang and others solved the ISDE as a limit of solutions of finite-volume stochastic differential equations (SDE), describing particles in $ \SR $ with reflecting boundary condition on $ \partial \SR $. That is, the SDE is given by 
\begin{align}\label{:12b} 
 d X_t^i=& dB_t^i 
- \frac{\beta }{2} \sum_{j\not=i}^{\xxSR } 
\nabla \Psi (X_t^i -X_t^j )dt 
- \frac{\beta }{2} \sum_{j> \xxSR }^{\infty}
\nabla \Psi (X_t^i - x_j )dt
\\ \notag & \quad \ \, \, + 
\half \mathbf{n}^R(X_t^i )dL_t^{R,i} 
\ \quad \ \ \text{ for } 1\le i \le \xxSR 
,\\ \notag 
X_t^i =& x_i\quad \quad \quad \quad \quad \quad \quad \quad \quad \quad 
\, \text{ for } i> \xxSR 
\end{align}
with the initial condition $ \mathbf{X}_0 =(x_i)_{i=1}^{\infty} $ 
 such that $ |x_i | < |x_{i+1} |$ for all $ i \in\N $, and 
$ \xxSR $ coincides with the number of particles in $ \SR $. 
The process $L^{R,i}= \{L_t^{R,i}\} $ denotes the local time-type drift arising from the reflecting boundary condition on 
$ \partial \SR $ (see \eqref{:29n} for $L^{R,i}$) and 
$ \mathbf{n}^R (x)$ is the inward normal vector at 
$ x \in \partial \SR $. 

In contrast, the labeled diffusion in $ \SR $ 
introduced in \cite{o.dfa} is given by the SDE
\begin{align}\label{:12e}
d X_t^i =&dB_t^i 
- \frac{\beta }{2} \sum_{j\not=i,\ X_t^j \in \SR }
\nabla \Psi (X_t^i - X_t^j )dt 
\end{align}
with the foregoing boundary condition. %
These SDEs have thus different boundary conditions. 
The solutions of \eqref{:12b} are non-ergodic, 
whereas the solutions of \eqref{:12e} are ergodic. 
Indeed, the system in \eqref{:12b} keeps the initial number of particles in $ \SR $. 
In the second dynamics, the number of particles in $ \SR $ varies. 
The state space of solutions in \eqref{:12e} therefore consists of a unique ergodic component (regarded as 
$ \{ \bigcup_{m=0}^{\infty} (\SR ^{\mathrm{int}})^m \}$-valued process, where $ (\SR ^{\mathrm{int}})^0 =\{\emptyset \} $ and 
$ \SR ^{\mathrm{int}}$ is the interior of $ \SR $). 

Let $ (\E_R^{\La},\mathcal{D} _R^{\La})$ be the Dirichlet form on $ \Lm $ associated with \eqref{:12b}, that is, the Dirichlet form $ (\E_R^{\La},\mathcal{D} _R^{\La})$ 
describes the motion of unlabeled dynamics related to \eqref{:12b}. 
Let $ (\E^{\Os},\mathcal{D} _R^{\Os})$ be the Dirichlet form 
on $ \Lm $ associated with \eqref{:12e}. 
Here we use the notation $ (\E^{\Os},\mathcal{D} _R^{\Os})$ rather than 
$ (\E_R^{\Os},\mathcal{D} _R^{\Os})$ because 
$ (\E^{\Os},\mathcal{D} _R^{\Os})$ is the closure of 
$ (\E , \dimu \cap \mathcal{B}_R(\mathsf{S} ))$ 
(see \lref{l:21} for notational details), 
whereas $ (\E_R^{\La},\mathcal{D} _R^{\La})$ is the closure with respect to the energy form $ \E_R $ different from $ \E $. 
As we shall see later, these two Dirichlet forms satisfy the relation 
\begin{align}&\notag 
(\E_R^{\La},\mathcal{D} _R^{\La}) \le (\E^{\Os},\mathcal{D} _R^{\Os})
.\end{align}
Furthermore, as $ R\to \infty $, 
$ \{ (\E_R^{\La},\mathcal{D} _R^{\La}) \} $ is an increasing scheme of Dirichlet forms, whereas $ \{ (\E^{\Os},\mathcal{D} _R^{\Os}) \} $ is decreasing. 
This fact implies the obvious relation \eqref{:11f}. 

The difference in these schemes lies in the boundary condition. 
Therefore, our task is to control the effect of the boundary condition to prove it becomes negligible as $ R \to \infty$. 

The main examples of our models have a logarithmic interaction potential, which is a very long rang potential that has quite strong long-range effect. 
We emphasize that the ISDEs arising from random matrix theory usually have logarithmic interaction potentials, and hence this class of interacting Brownian motions is significant. 

Loosely speaking, a solution $ X $ of an SDE is called a weak solution if it is a pair $ (X ,B )$ of a continuous process $ X $ and Brownian motion $ B $ that satisfies the SDE. We call $ X $ is a strong solution if, in addition, $ X $ is a function of Brownian motion $ B $ \cite{IW}. 

In the first main theorem (\tref{l:31}), we shall prove that any limit point of the {\ws}s of \eqref{:12b} is a {\ws}\ of the ISDE \eqref{:12a} satisfying well-behaved properties (see \tref{l:31}). 
The limit points of the {\ws}s of \eqref{:12e} were proved to satisfy the ISDE \eqref{:12a} with the same well-behaved properties 
\cite{o.dfa,o.isde}. 
Hence, assuming the uniqueness of {\ws}s of \eqref{:12a} under the foregoing well-behaved properties, we proved that 
these two limits of the {\ws}s are the same. 
This establishes the coincidence of the two Dirichlet forms 
$ (\E^{\La},\mathcal{D} ^{\La}) $ and $(\E^{\Os},\mathcal{D} ^{\Os}) $ 
in the second main theorem (\tref{l:32}). 
 
The motivation of our work lies in the recent rapid and outstanding progress of random matrix theory in proving that the random point fields arising from Gaussian random matrices (invariant random matrices) such as sine$ _{\beta }$, 
Bessel$ _{\beta }$, and Ginibre random point fields, are universal. 
Indeed, these random point fields are obtained as scaling limits of eigenvalue distributions of a quite general class of random matrices and also log gases with general free potentials. 
Once this static universality is established, it is natural to pursue its dynamical counterpart. 
In a forthcoming paper, the first and second authors will prove that the natural reversible stochastic dynamics associated with these random point fields are also universal objects. 
Examples of universal stochastic dynamics are the sine, Bessel, and Ginibre interacting Brownian motions (see \sref{s:8}). 
They are limits of the stochastic dynamics related to $ N $-particle systems with reversible random point fields that converge to those universal random point fields mentioned above. 
This result is established by assuming the limits of the lower and upper Dirichlet forms in \eqref{:11e} are equal in addition to a certain strong convergence of the random point fields. 
Hence the second main theorem (Theorem 3.2) plays a crucial role in the dynamical universality of random matrices in the sense given above. 

The organization of the paper is as follows: 
In \sref{s:2}, we prepare the two schemes of the Dirichlet forms describing interacting Brownian motions, and quote related results. 
In \sref{s:3}, we state the main theorems (\tref{l:31} and \tref{l:32}). 
In \sref{s:4}, we prove \tref{l:31}. 
In \sref{s:5}, we prove \tref{l:32} and comment on a generalization to the uniformly elliptic case. 
In \sref{s:7}, we construct cut-off coefficients $ \mathsf{b}_{r,s,\p }$ appearing in \As{A6}. 
In \sref{s:8}, we present examples. 
In \sref{s:9}, we give concluding remarks with some open questions. 
\section{Preliminaries}\label{s:2} 
\subsection{Two schemes of Dirichlet forms}\label{s:21} 

Let $ \sS $ be a closed set in $ \Rd $ with interior 
$ \sS ^{\mathrm{int}}$ which is a connected open set 
satisfying $ \overline{\sS ^{\mathrm{int}}} = \sS $ and 
the boundary $ \partial \sS $ having Lebesgue measure zero. 

A configuration $ \sss = \sum_i \delta _{s_i}$ on $ \sS $ is 
a Radon measure on $ \sS $ consisting of delta masses $\delta_{s_i}$. 
Let $ \mathsf{S} $ be the configuration space over $ \sS $. 
Then, by definition, $ \mathsf{S} $ is the set given by 
\begin{align} & \label{:20a}
\mathsf{S} = \{ \sss = \sum_{i} \delta _{s_i}\, ;\, 
 \text{ $\sss ( K ) < \infty $ for all compact set $ K $} \} 
.\end{align}
By convention, we regard the zero measure as an element of $ \mathsf{S} $. 
We endow $ \mathsf{S} $ with the vague topology, which makes $ \mathsf{S} $ to be a Polish space. 

A probability measure $\mu $ on $(\mathsf{S} , \mathcal{B}(\mathsf{S} ))$ is called a random point field on $S$. We assume that $ \mu $ is supported on the set consisting of infinitely-many particles: 
\begin{align}&\label{:20b}
 \mu (\{\mathsf{s}\in \mathsf{S} ; \mathsf{s}(\sS ) = \infty \} ) = 1
.\end{align}

Let $ \sS _r = \{ s \in \sS \, ;\, | s | \le r \} $ and 
$ \Srm = \Sr \ts \cdots \ts \Sr $ be the $ m$-product of $ \Sr $. 
Let $ \mathsf{S} _r^{m} = \{ \sss \in \mathsf{S} \, ;\, \sss (\Sr ) = m \} $ 
for $r,m\in\N$. 
We set maps $ \map{\pir , \pirc }{\mathsf{S} }{\mathsf{S} }$ such that 
\begin{align}&\notag 
\text{$ \pir (\sss) = \sss (\cdot \cap \Sr )$\quad and \quad 
$ \pirc (\sss) = \sss (\cdot \cap \Sr ^c)$}
.\end{align} 
For $\sss\in\mathsf{S}_r^m$, 
$\mathbf{x}_r^m(\sss)=(x_r^{i}(\sss))_{i=1}^m \in S_r^m$
is called a $S_r^m$-coordinate of $\sss$ 
if $\pi_r(\sss)=\sum_{i=1}^m \delta_{x_r^{i}(\sss)}$.

\smallskip 

For a function $f: \mathsf{S} \to \R$ we set 
$ f_{r}^m (\mathsf{s},\mathbf{x}) = f_{r,\sss}^m (\mathbf{x}) $ 
such that $f_{r}^m : \mathsf{S} \times S_r^m \to \mathbb{R}$ 
and that $ f_{r,\sss}^m $, with $r, m \in\N $, satisfies 

\noindent 
(1) $f_{r,\sss}^m$ is a permutation invariant function on $S_r^m$ for each $\sss\in\mathsf{S} $.
\\
(2) $f_{r,\sss(1)}^m =f_{r,\sss(2)}^m $ 
if $\pirc (\sss(1))=\pirc (\sss(2))$ for 
$\sss(1), \sss(2)\in \mathsf{S}_ r^m$.
\\
(3) $f_{r,\sss}^m(\mathbf{x}_r^m(\sss))=f(\sss)$ for $\sss\in \mathsf{S}_ r^m$, where $\mathbf{x}_r^m(\sss)$ is a $S_r^m$-coordinate of $\sss$.
\\
(4) $f_{r,\sss}^m(s)=0$ for $\sss\notin \mathsf{S}_ r^m$.

\ms 

The function $f_{r,\sss}^m$ is called the $S_r^m$-representation of $f$. 
Note that $f_{r,\sss}^m$ is unique and 
$f(\sss)=\sum_{m=0}^\infty f_{r,\sss}^m(\mathbf{x}_r^m(\sss))$ for each $ r \in \N $. 
When $f$ is $\sigma[\pir ]$-measurable, the $S_r^m$-representations are independent of $\sss $. In this case we often write $f_r^m$ instead of $f_{r,\sss}^m$. 
We set 
\begin{align} &\notag 
\mathcal{B}_r(\mathsf{S} )=\{ f: \mathsf{S} \to\R \,;\, 
\mbox{ $f$ is $\sigma [\pir ]$-measurable}\}
,\quad 
\mathcal{B}_\infty(\mathsf{S} )=\bigcup_{r=1}^\infty
\mathcal{B}_r(\mathsf{S} )
,\\&\notag 
\di=\{ f\in \mathcal{B}_\infty(\mathsf{S} ) \,;\,
\mbox{ $f_{r,\sss}^m$ is smooth on $S_r^m$ for all $m,r\in\N$, $\sss\in \mathsf{S} $}\}
.\end{align}
Note that $\di \cap \Lm $ is dense in $\Lm $ and $\di \subset C(\mathsf{S} )$, 
where $ C(\mathsf{S} )$ is the set of all continuous functions on $ \mathsf{S} $. 
We remark that, if $ f \in \di $ and $ f $ is $\sigma[\pir ]$-measurable, then 
$f_{r,\sss}^m (x_1,\ldots,x_m)$ is constant in $ x_m $ 
on the boundary $ \partial \Sr $ for each 
$ (x_1,\ldots,x_{m-1})\in \Sr ^{m-1} $, 
and its derivatives vanishes on $ \partial \Srm $. 

For $ \sss=\sum_i \delta_{s_i}$ we set 
$ \nabla_{s_i} = (\PD{}{s_{i1}},\ldots,\PD{}{s_{id}})$. 
For $f,g \in \di $ let 
\begin{align} \label{:21a}
&\mathbb{D}_r^{m} [f,g] (\sss) = 
\begin{cases}
 \frac{1}{2} \sum_{i\,;\, s_i\in S_r} 
\nabla_{s_i} f_{r,\sss}^m(\mathbf{x}_r^m(\sss)) \cdot \nabla_{s_i} 
 g_{r,\sss}^m(\mathbf{x}_r^m(\sss)) \quad \text{ for }
\sss\in \mathsf{S}_ r^m 
,\\0 
\quad \text{ for } \sss\notin \mathsf{S}_ r^m
.\end{cases}
\end{align}
Moreover, we set 
\begin{align}\label{:21b}&
\mathbb{D}_r=\sum_{m=1}^\infty\mathbb{D}_r^{m}
.\end{align}
Note that $\mathbb{D}_r^{m}[f,g]$ is independent of the choice of the $S_r^m$-coordinate $\mathbf{x}_r^m(\sss)$ and is well-defined.
We now define bilinear forms on $\di $: 
\begin{align} & \label{:21c} 
\Ermum (f,g) =\int_{\mathsf{S} }\mathbb{D}_r^{m}[f,g](\sss)\mu(d\sss)
 \quad \mbox{ and }\quad
\Ermu (f,g) = 
\int_{\mathsf{S} }\mathbb{D}_r [f,g](\sss)\mu(d\sss) 
.\end{align}
Then clearly $ \Ermu =\sum_{m=1}^\infty \Ermum $.

Let $ (\Emu ,\dimu )$ be a bilinear form 
on $ \Lm $ with domain $ \dimu $ defined by 
\begin{align} 
\label{:21e}& 
\Emu (f,f) 
= \lim_{r\to\infty}\Ermu (f,f) 
,\\ \notag &
\dimu = 
\{ f \in \di \cap \Lmu \, ;\, \Emu (f,f) < \infty \} 
.\end{align}
We note that $ \Ermu (f,f)$ is nondecreasing in $ r $, and hence the limit in \eqref{:21e} exists. We assume 
\begin{align} &\label{:21f}
(\Ermum ,\dimu )
\mbox{ is closable on $\Lm $ for each $m,r \in\mathbb{N}$}
.\end{align}
We present later a sufficient condition regarding \eqref{:21f}; see \As{A1} in \sref{s:22}. 

\begin{lemma}[{\cite[Lemma 2.2, Theorem 2]{o.dfa}}] \label{l:21}
Assume \eqref{:21f}. Then the following hold. 
\\
\noindent(1) 
$(\Emu ,\dimu \cap \mathcal{B}_r(\mathsf{S} ))$ and 
$(\Ermu ,\dimu )$
are closable on $\Lm $ for each $r$.
\\\thetag{2} 
$(\Emu , \dimu )$ is closable on $\Lm $. 
\end{lemma}

For symmetric bilinear forms 
$ (\E ,\mathcal{D} )$ and $ (\mathcal{E}',\mathcal{D}' )$ we write 
 $ (\E ,\mathcal{D} ) \le (\mathcal{E}',\mathcal{D}' )$ 
if 
$ \mathcal{D}\supset \mathcal{D}' $ and 
$ \mathcal{E}(f,f)\le \mathcal{E}'(f,f) $ 
for all $f \in \mathcal{D}' $. 
We say a sequence of symmetric bilinear forms 
$ \{ (\mathcal{E}_n,\mathcal{D}_n ) \}_{n\in\N} $
 is increasing if 
$ (\mathcal{E}_n,\mathcal{D}_n ) \le 
(\mathcal{E}_{n+1},\mathcal{D}_{n+1} ) $ for all $ n $. 
Replacing $ \le $ by $ \ge $, we call 
$ \{ (\mathcal{E}_n,\mathcal{D}_n ) \}_{n\in\N} $ decreasing. 
Let $\EDrupr $ and $\EDrlwr $ 
denote the closures of 
$(\Emu ,\dimu \cap \mathcal{B}_r(\mathsf{S} ))$ and 
$(\Ermu ,\dimu )$ on $ \Lm $,
respectively. Then we quote: 
\begin{lemma} [{\cite[Lemma 2.2]{o.dfa}}] \label{l:22} 
Assume \eqref{:21f}. Then 
\\\thetag{1}
$\{\EDrlwr \}_{r\in\mathbb{N}}$ is increasing.
\\\thetag{2}
$\{\EDrupr \}_{r\in\mathbb{N}}$ is decreasing.
\end{lemma}

Let $(\E^{\La}, \mathcal{D} ^{\La})$ be the increasing limit of 
$\{\EDrlwr \}_{r\in\mathbb{N}}$, that is, 
\begin{align}&\notag 
\E ^{\La}(f,f)=\limi{r} \Erlwr (f,f),
\; \mbox{ and } \; \mathcal{D} ^{\La}=
\{f\in \bigcap_{r\in\N }
\mathcal{D} _r^{\La } \,;\, \limi{r}\Erlwr (f,f)<\infty \} 
.\end{align}
Then $(\E^{\La}, \mathcal{D} ^{\La})$ is closed on 
$\Lm $ by construction. 
Recall that $(\Emu, \dimu )$ is closable on $\Lm $ 
by \lref{l:21}. 
We then denote by $(\E^{\Os}, \mathcal{D} ^{\Os})$ the closure of 
$(\Emu, \dimu )$ on $\Lm $. 
Let $G_{r,\alpha }^{\La} $, $G_{\alpha }^{\La} $, 
$G_{r,\alpha }^{\Os} $, and $G_{\alpha }^{\Os} $ 
be the resolvents of $\EDrlwr $, $(\E^{\La}, \mathcal{D} ^{\La})$,
$\EDrupr $, and $(\E^{\Os}, \mathcal{D} ^{\Os})$ on $\Lm $, respectively. 
\begin{lemma}[{\cite[Lemma 2.1, Theorem 2]{o.dfa}}] \label{l:23}
Assume \eqref{:21f}. Then 
\\\thetag{1} 
$ \{ G_{r,\alpha }^{\La}\}_{r\in\N} $ converges to $G_{\alpha }^{\La} $ 
strongly in $\Lm $ as $ r \to \infty $ for each $\alpha>0$.
\\
\thetag{2} 
$ \{ G_{r,\alpha }^{\Os}\}_{r\in\N} $ converges to $G_{\alpha }^{\Os} $ 
strongly in $\Lm $ as $ r \to \infty $ for each $\alpha>0$.
\end{lemma}

By construction we have for each $ r $
\begin{align}\label{:23a}&
\EDrlwr \le \EDrupr 
.\end{align}
Hence taking $ r \to \infty $ we see that 
\begin{align}\label{:23b}&
(\mathcal{E}^{\La}, \mathcal{D} ^{\La}) \le 
(\mathcal{E}^{\Os}, \mathcal{D}^{\Os})
.\end{align}
We call 
$ (\mathcal{E}^{\La}, \mathcal{D} ^{\La})$ the lower Dirichlet form and 
$ (\mathcal{E}^{\Os}, \mathcal{D}^{\Os})$ 
the upper Dirichlet form. We also call 
$ \{ \EDrlwr \}_{r\in\N } $ a lower scheme and 
$ \{ \EDrupr \}_{r\in\N }$ an upper scheme. 
The relations \eqref{:23a} and \eqref{:23b} justify 
the names of these schemes. 
\subsection{Quasi-Gibbs measures, unlabeled diffusions, and labeled dynamics}
\label{s:22}

Let $ \Lambda_r $ be the Poisson random point field 
 whose intensity is the Lebesgue measure on $ \Sr $ and set 
 $ \Lambda_r^m = \Lambda_r (\cdot \cap \mathsf{S} _r^{m} ) $. 
Let $ \map{\Phi }{\sS }{\R \cup \{ \infty \} }$ and 
$ \map{\Psi }{\sS ^2 }{\R \cup \{ \infty \} }$ be measurable functions such that 
$ \Psi (x,y) = \Psi (y,x)$. Following \cite{o.rm,o.rm2} we quote: 
\begin{definition}\label{d:21}
 A random point field $ \mu $ is called 
 a $ ( \Phi , \Psi )$-quasi Gibbs measure if its regular conditional probabilities 
\begin{align}&\notag 
 \mu _{r,\sss }^{m} = 
\mu (\, \pir (\xx ) \in \cdot \, | \, \pirc (\mathsf{x}) = 
 \pirc (\mathsf{\sss }),\, \mathsf{x}(\Sr ) = m )
\end{align}
satisfy, for all $r,m\in \mathbb{N}$ and $ \mu $-a.s.\! $ \sss $, 
\begin{align}\label{:24b}&
\cref{;10Q}^{-1} e^{-\Hrm(\mathsf{x}) } \Lambda_r^{m} (d\mathsf{x}) \le 
 \mu _{r,\sss }^{m} (d\mathsf{x}) \le 
\cref{;10Q} e^{-\Hrm(\mathsf{x}) } \Lambda_r^{m} (d\mathsf{x}) 
.\end{align}
 Here $ \Ct \label{;10Q} = \cref{;10Q} (r,m,\sss )$ 
 is a positive constant depending on $ r ,\, m ,\, \sss $. 
 For two measures $ \mu , \nu $ on a $ \sigma $-field $ \mathcal{F} $, we write 
$ \mu \le \nu $ if $ \mu (A) \le \nu (A) $ for all $ A \in \mathcal{F} $. 
Moreover, $ \Hrm $ is the Hamiltonian on $ \Sr $ defined by 
\begin{align}& \notag 
\Hrm (\mathsf{x}) = \sum_{x_i\in \Sr \atop 1\le i \le m } \Phi (x_i) 
+ \sum_{ x_j, x_k \in \Sr \atop 1\le j < k \le m } \Psi (x_j,x_k)
\quad \text{ for } \mathsf{x} = \sum_i \delta_{x_i} 
.\end{align}
\end{definition}

\begin{remark}\label{r:A2} \thetag{1} 
From \eqref{:24b}, we see that for all $r,m\in \mathbb{N}$ and 
$ \mu $-a.s.\! $ \sss $, 
$ \mu _{r,\sss }^{m} (d\mathsf{x}) $ have (unlabeled) Radon--Nikodym densities 
$ m_{r,\sss }^{m} (\mathsf{x}) $ with respect to $ \Lambda_r^m $. 
Clearly, the canonical Gibbs measures $ \mu $ with potentials $ (\Phi , \Psi )$ 
are quasi-Gibbs measures, and their densities 
$ m_{r,\sss }^{m} (\mathsf{x}) $ with respect to $ \Lambda_r^{m} $ 
are given by the Dobrushin--Lanford--Ruelle (DLR) equation. 
That is, for $ \mu $-a.s.\! $ \sss = \sum_j \delta_{s _j} $ 
\begin{align}&\notag 
m_{r,\sss }^{m} (\mathsf{x}) = \frac{1}{\mathcal{Z}_{r,\sss }^{m} }
\exp \{ - \Hrm (\mathsf{x}) - 
\sum_{ x_i \in \Sr ,\, s _j \in \Sr ^c \atop 1\le i \le m} 
\Psi (x_i, s _j)\} 
.\end{align}
Here $ \mathcal{Z}_{r,\sss }^{m} $ is the normalizing constant. 
For random point fields appearing in random matrix theory, interaction potentials are logarithmic functions, 
where the DLR equations do not make sense as stated because the term 
$ \sum_{ x_i \in \Sr ,\, s _j \in \Sr ^c } \Psi (x_i, s _j) $ diverges. 
The notion of a quasi-Gibbs measure still makes sense for logarithmic potentials.\\\thetag{2} 
We refer to \cite{o.rm,o.rm2} for sufficient conditions of quasi-Gibbs measures. 
These conditions give us the quasi-Gibbs property of random point fields appearing in random matrix theory, such as 
sine$_{\beta} $ and Bessel$_{2, \alpha } $ ($ 1\le \alpha $), and 
Ginibre random point fields \cite{h-o.bes,o.rm,o.rm2}. 
\end{remark}

We make the following assumption. 

\smallskip 
\noindent 
\As{A1} $ \mu $ is a $ ( \Phi , \Psi )$-quasi Gibbs measure. 
Furthermore, there exists 
upper semi-continuous functions $ (\hat{\Phi }, \hat{\Psi })$ 
 and positive constants $ \Ct \label{;A21}$ and $ \Ct \label{;A22} $ satisfying 
\begin{align}&\notag 
\cref{;A21}^{-1} \hat{\Phi } (x)\le \Phi (x) \le \cref{;A21} \hat{\Phi }(x) ,\quad 
\cref{;A22}^{-1} \hat{\Psi } (x,y)\le \Psi (x,y) \le \cref{;A22} \hat{\Psi }(x,y) 
,\end{align}
where $ \Psi $ and $ \hat{\Psi }$ satisfy 
$ \Psi (x,y)=\Psi (y,x)$ and $ \hat{\Psi }(x,y) = \hat{\Psi }(y,x)$. 

\ms 

If these interaction potentials are translation invariant, we often write $ \Psi (x,y) = \Psi (x-y)$ and 
$ \hat{\Psi } (x,y) = \hat{\Psi }(x-y)$. 
The importance of \As{A1} lies in the fact that it gives a sufficient condition of the basic assumption \eqref{:21f}. We quote: 
\begin{lemma}[{\cite[45-46pp]{o.rm}}] \label{l:24}
Assume \As{A1}. Then $(\Ermum ,\dimu )$ is closable 
 on $\Lm $ for each $m,r \in\mathbb{N}$. 
In particular, 
$(\Ermu ,\dimu )$ is closable on $\Lm $. 
\end{lemma} 

We now recall two basic notions on random point fields: 
correlation functions and density functions. 

A symmetric and locally integrable function $\rho^n : S^n\to [0,\infty)$ is called the $n$-point correlation function of a random point field $\mu $ on $S$ with respect to the Lebesgue measure if $\rho^n$ satisfies 
\begin{align}&\notag 
\int_{A_1^{k_1}\times\cdots \ts A_m^{k_m}}\rho^n(x_1,\dots,x_n)dx_1\cdots dx_n
=\int_{\mathsf{S} }\prod_{i=1}^m \frac{\sss (A_i)!}{(\sss (A_i)-k_i)!}d\mu
\end{align}
for any sequence of disjoint bounded measurable sets $A_1,\dots,A_m\in \mathcal{B}(S)$ and a sequence of positive integers $k_1,\dots, k_m$ satisfying $k_1+\cdots + k_m=n$. When $\sss (A_i)-k_i<0$, according to our interpretation, $\sss (A_i)!/(\sss (A_i)-k_i)!=0$ by convention. 
We assume that $ \mu $ has the $ n $-point correlation function $ \rho ^n $ for each $ n \in\N $.

A symmetric function 
$\sigma _r^k : \Srk \to [0,\infty)$ is called the $k$-point density function of a random point field $\mu $ on $\Sr $ with respect to the Lebesgue measure if for all non-negative, bounded 
$ \sigma [\pir ]$-measurable function $ f $ with $ \Srk $-representation $ f_r^k $ 
\begin{align}&\notag 
\frac{1}{k! }\int_{\Srk } f_r^k \sigma _r^k d\mathbf{x}^k = 
\int_{\SSrk } f d \mu 
.\end{align}
Here $ \mathsf{S} _r^{k} = \{ \sss \in \mathsf{S} \, ;\, \sss (\Sr ) = k \} $ as before. 
We make the following assumption. 

\ms

\noindent \As{A2} \ \ $ \mu $ has correlation functions with respect to the Lebesgue measure 
of any order and $ \mu $ satisfies for each $ m, r \in \N $ 
\begin{align}
\label{:51h}&
\sum_{k=m}^{\infty} \frac{k!}{(k-m)!} \mu (\mathsf{S} _r^k ) < \infty 
.\end{align}
Clearly, \eqref{:51h} is equivalent to $ \int_{\Sr }\rho ^m (x) dx < \infty $ for all $ r \in \mathbb{N}$ 
whenever the $ m $-point correlation function $ \rho ^m $ of $ \mu $ 
with respect to the Lebesgue measure exists. 
Under the assumptions of \As{A1} and \eqref{:51h} we deduce that $ \mu $ 
has correlation functions and Campbell measures of any order. 
Hence we can dispense with the existence of correlation functions in \As{A2}.

\ms 

A family of probability measures $ \{ \PPx \}_{\xx \in \mathsf{S} } $ on 
$ \WS $ is called a diffusion if the canonical process 
$ \mathsf{X}=\{ \mathsf{X}_t \} $ under $ \PPx $ is a continuous process 
having a strong Markov property starting at $ \xx $. 
Furthermore, $ \{ \PPx \}_{\xx \in \mathsf{S} } $ is called conservative 
if it has an invariant probability measure. 

Assume \As{A1}. Then we deduce from \lref{l:21} and \lref{l:24} 
that the non-negative form $ (\Emu ,\dimu )$ is closable on 
$ \Lmu $. Therefore, let 
$ (\E^{\Os} ,\mathcal{D} ^{\Os} ) $ be its closure on $ \Lmu $. 
The next result is a refinement of \cite[119p.\! Corollary 1]{o.dfa} and can be proved in a similar fashion. 
%
\begin{lemma} \label{l:25}
Assume \As{A1} and \As{A2}. Then there exists a $ \mu $-reversible diffusion 
$\{ \PPxU \}_{\xx \in \mathsf{S} } $ associated with the Dirichlet form $ (\E^{\Os} ,\mathcal{D} ^{\Os} ) $ 
on $ \Lmu $. 
\end{lemma}

\PF
The assumption \As{A2} corresponds to \thetag{A.2} in 
\cite{o.dfa} (we write \thetag{A.2} below). 
It was assumed in \thetag{A.2}, in addition to \As{A2}, 
 the boundedness of density functions of all order on each $ \Sr $. 
In \cite{o.dfa}, \thetag{A.2} was used only in the proof of 
\thetag{2.2} in Lemma 2.4 (see \cite[125p]{o.dfa}). 
Moreover, \thetag{2.2} in \cite{o.dfa} is used only to prove 
Lemma 2.4 \thetag{3} in \cite{o.dfa}, which is the claim such that $ \dimu $ is dense in $ \Lm $. 
Hence our task is to prove this under \As{A2}. 
For this purpose we recall a mollifier on $ \mathsf{S} $ 
introduced in \cite{o.dfa}. 

Let $ \map{j}{\Rd }{\R }$ 
be a non-negative, smooth function such that 
$ \int_{\Rd } j(x)dx = 1$, $ j (x)=0 $ for $ |x|\ge 1/2$. 
Let $ j_{\varepsilon } (x)= \varepsilon ^d j (x/\varepsilon )$ and 
$ j_{\varepsilon } ^i ((x_1,\ldots,x_i)) = 
\prod_{j=1}^i j_{\varepsilon } (x_j)$. 

For a $ \sigma [\pi_r]$-measurable, bounded function $ f $ 
we set 
\begin{align}\label{:92a}&
\mathcal{J}_{r,\varepsilon }f (\mathsf{s}) = 
\begin{cases}
 j_{\varepsilon } ^i * f_r^i (\mathbf{x}_r^i (\mathsf{s}))
&\text{ for } \mathsf{s} \in \SSr^i \quad (i\ge 1),
\\
f (\mathsf{s}) 
&\text{ for } \mathsf{s} \in \SSr ^0 
.\end{cases}
\end{align}
Here $ f_r^i$ is $ \sS _r^i$-representation of $ f $ and 
$ \mathbf{x}_r^i (\mathsf{s})$ is the $ \sS _r^i$-coordinate 
introduced in \sref{s:21}. 
Moreover, $ * $ denotes the convolution in $ (\Rd )^i$, that is, 
$$ j_{\varepsilon } ^i * f_r^i (\mathbf{x}) = 
\int_{(\Rd )^i} j_{\varepsilon } ^i (\mathbf{x}-\mathbf{y})
f _r^i (\mathbf{y}) d\mathbf{y},$$
where we set 
$ f_r^i (\mathbf{x}) = 0 $ for $ \mathbf{x}\not\in \sS _r^i$.

We now prove that $ \dimu $ is dense in $ \Lm $. 

Let $ 0 < \delta < r $ ($\delta \in \R $) and $ f \in C_b(\mathsf{S} )\cap \mathcal{B}_{r-\delta } $. 
Then $ f $ is a bounded continuous, $ \sigma [\pi_{r-\delta }]$-measurable function on $ \mathsf{S} $ by definition. 
From \cite[Lemma 2.4 \thetag{2.1}]{o.dfa} we see $ \mathcal{J}_{r,\varepsilon }f \in \dimu $ 
for $ 0 < \varepsilon < \delta $. Moreover, because $ f \in C_b(\mathsf{S} ) $, we see from \eqref{:92a} that 
\begin{align}\label{:92A}&
\limz{\varepsilon } \mathcal{J}_{r,\varepsilon }f (\mathsf{s}) = 
f (\mathsf{s}) \quad \text{ for each } \mathsf{s}
,\\\label{:92p}&
\sup_{\mathsf{s}\in \mathsf{S} }
|\mathcal{J}_{r,\varepsilon }f (\mathsf{s}) |
\le 
\sup_{\mathsf{s}\in \mathsf{S} }
|f (\mathsf{s}) | < \infty 
.\end{align}
From \eqref{:92A}--\eqref{:92p} we can apply the Lebesgue convergence theorem to obtain 
\begin{align}\label{:92b}&
\limz{\varepsilon }
\int_{\mathsf{S} } | \mathcal{J}_{r,\varepsilon }f (\mathsf{s}) - 
f (\mathsf{s}) |^2 \mu (d\mathsf{s}) =
\int_{\mathsf{S} } 
\limz{\varepsilon }
| \mathcal{J}_{r,\varepsilon }f (\mathsf{s}) - 
f (\mathsf{s}) |^2 \mu (d\mathsf{s}) 
= 0 
\end{align}
for each $ r \in \N $. Because 
\begin{align*}&
\bigcup_{r=1}^{\infty} \bigcup_{0 < \delta < r ,\, \delta \in \R }
C_b(\mathsf{S} )\cap \mathcal{B}_{r-\delta } 
\end{align*}
is dense in $ \Lm $, we see $ \dimu $ is dense in $ \Lm $ from \eqref{:92b}. 
\PFEND

We note that \As{A2} is used to guarantee the existence of the diffusion. 
The $ \mu $-reversibility of the diffusion follows from 
$ 1 \in \mathcal{D} ^{\Os} $ and symmetry of $ (\E^{\Os} ,\mathcal{D} ^{\Os} ) $. 

By construction, such a family of diffusion measures 
$ \PPU = \{\PPxU \}$ with quasi-continuity in $ \xx $ is unique 
for quasi-everywhere starting point $ \xx $. 
Equivalently, there exists a set $ \mathsf{S} _0 $ such that 
the complement of $ \mathsf{S} _0$ has capacity zero and 
the family of diffusion measures $ \PPU = \{\PPxU \}$ associated 
with the Dirichlet space above with quasi-continuity in $ \xx $ 
is unique for all $ \xx \in \mathsf{S} _0 $ 
and $ \PPxU (\mathsf{X}_t \in \mathsf{S} _0 \text{ for all }t ) = 1 $ 
for all $ \xx \in \mathsf{S} _0 $.

We next lift the unlabeled dynamics $\mathsf{X}$ to a labeled dynamics 
$\mathbf{X}=(X^{i})_{i\in\N}$. 
Under $ \PPU = \{\PPxU \}$, we can write $ \mathsf{X}_t = \sum_{i=1}^{\infty} \delta_{X_t^i}$, 
where each $ X^i=\{ X_t^i \} $ is a continuous process with time parameter of the form 
$ [0,b)$ or $ (a,b)$. 
We call each $ X^i$ a tagged particle of $ \mathsf{X}$. 
%

For a given unlabeled dynamics $ \mathsf{X}=\sum_i \delta_{X^i}$, 
we call the collection of tagged particles $\mathbf{X}=(X^{i})_{i\in\N }$ labeled dynamics. 
Note that for a given $ \mathsf{X}$, there exist plural labeled dynamics in general. 
We next give a condition such that $\mathbf{X}=(X^{i})_{i\in\N}$ is determined uniquely up to initial labeling. For this purpose, we impose the following conditions: 

\bs 
\noindent \As{A3} 
Under $ \PPU = \{\PPxU \}$, each tagged particle 
$\{X^{i}\}_{i\in\N }$ does not collide with another. 
Furthermore, $\{X^{i}\}_{i\in\N }$ does not hit the boundary $ \partial \sS $ of $ \sS $. 

\bs 
This condition is equivalent to 
both the capacity of multiple points and 
that of configurations with particles on the boundary $ \partial \sS $ being zero: 
\begin{align}&\notag 
\Capa (\{ \mathsf{s}\in \mathsf{S} ; \mathsf{s}(\{ x \} ) 
\ge 2 \text{ for some } x \in \sS \} ) = 0 
,\\&\notag 
\Capa (\{ \mathsf{s}\in \mathsf{S} ; \mathsf{s}(\partial \sS ) 
\ge 1 \} ) = 0 
.\end{align}
Here $ \Capa $ denotes the one-capacity with respect 
to the Dirichlet space $ (\Emu ,\dom )$ on $ \Lm $ 
(see \cite{fot.2} for the definition of capacity). 

By \As{A3}, we write $ \mathsf{X}_t = \sum_i \delta_{X_t^i} $ such that 
$ X^i \in C(I_i;\sS )$, where $ I_i = [0,b)$ or $ I_i = (a,b)$ for some $ a, b \in (0,\infty ]$. 
We remark that $ \{X^i \}$ is unique up to labeling.

\bs 

Let $ \mathrm{Erf} (t) = 
\int_t^{\infty} (1/\sqrt{2\pi }) e^{-|x|^2/2} dx $ be the error function. We further assume: 

\ms

\noindent \As{A4} 
There exists a $ T>0 $ such that, for each $ R>0 $, 
\begin{align}\label{:26y}&
\liminf_{r\to\infty} \, \mathrm{Erf} \Big(\frac{r}{\sqrt{(r+R)T}}\Big) \, 
\int_{|x|\le {r+R}} \rho ^1 (x)dx = 0 
.\end{align}
For each $ r,T \in \N $, the following holds. 
\begin{align}\label{:28x}
& \int_{\sS }\mathrm{Erf} \Big( \frac{|x|-r}{T}\Big) \rho ^1 (x) dx 
< \infty 
.\end{align}

\medskip

We say non-explosion of each tagged particle holds if the right end of $ I_i $ is infinity for each $ i $, and 
non-entering of each tagged particle holds if $ I_i =[0,b_i)$ for each $ i $.

\begin{lemma} \label{l:92} 
Assume \As{A1}--\As{A3} and \eqref{:26y}. Then the non-explosion and the non-entering 
of each tagged particle hold. 
Furthermore, $ \mathrm{Cap}(\{\mathsf{s} \in \mathsf{S} \, ; \mathsf{s}(\sS ) < \infty \} ) = 0 $. 
\end{lemma}%
\PF 
Let $ X^i $ be a tagged particle with interval $ I_i = [0,b_i)$. Then by \cite[Theorem 2.5]{o.tp} we have 
\begin{align}
\label{:92aT}&
\PPmU ( \, \sup \{|X _t^i |;\, t \le T \} < \infty \, \text{ for all } T, i \in \mathbb{N} \, ) = 1 
,\end{align}
where $ \PPmU = \int _{\mathsf{S} } \PPxU \mu (d\xx )$. 
Then \eqref{:92aT} implies the non-explosion of the tagged particles and 
$ I_i = [0,\infty )$. 

Next suppose that $ I_i = (a,b)$. 
Then, applying the Markov property of the diffusion $ \PPU = \{ \PPxU \} $ 
at time $ a' \in (a,b)$ and using the preceding argument together with stationarity of the diffusion, 
we deduce $ I_i = (a,\infty )$. 
Because of reversibility, we see that such open intervals do not exist. 
Hence, we obtain $ I_i = [0,\infty )$ for all $ i $. 
This implies the non-entering of the tagged particles. 

From the non-explosion and the non-entering of the tagged particles combined with \eqref{:20b}, 
we obtain $ \mathrm{Cap}(\{\mathsf{s} \in \mathsf{S} \, ; \mathsf{s}(\sS ) < \infty \} ) = 0 $. 
This completes the proof. 
\PFEND

We thus see from \As{A3} and \As{A4} that under $ \PPU =\{\PPxU \}$ 
the tagged particles of 
$\{X^{i}\}_{i\in\N }$ neither collide each other nor hit the boundary 
$ \partial \sS $ nor explode nor enter. 
The second condition \eqref{:28x} in \As{A4} will be used in the proof of \lref{l:28}. 

We call $\ulab$ the unlabeling map if $\ulab((x_i))=\sum_i \delta_{x_i}$.
We call $\lab $ a label if $\map{\lab}{\mathsf{S} }{\SN }$ is a measurable 
map defined for $\mu $-a.s. $\xx $, and $\ulab \circ \labx = \xx $. 

\begin{lemma}[\cite{o.tp,o.rm}]\label{l:26}
Assume \As{A1}--\As{A4}. Let $ \lab $ be a label. 
Then under $ \PPU = \{\PPxU \}$ 
there exists a unique, labeled dynamics 
$\mathbf{X}=(X^{i})_{i\in\N}\in \CS $ such that 
$\mathbf{X}_0=\lab (\mathsf{X}_0)$ and that 
$\mathsf{X}_t =\sum_{i\in\N}\delta_{X_t^{i}}$ for all $ t $. 
\end{lemma}

Once the initial label $ \lab $ is assigned, the particles are marked forever because they neither collide nor explode nor enter. 
We hence determine the labeled dynamics $ \mathbf{X} $ from the unlabeled dynamics $ \mathsf{X} $ and the label $ \lab $ 
uniquely. We have thus had a natural correspondence between 
$ \mathbf{X}$ and $ (\mathsf{X},\lab)$ under the conditions 
\As{A3} and \As{A4}. 
We remark here that $ \mathbf{X}_t \not= \lab (\mathsf{X}_t )$ for $ t > 0 $ in general.

We introduce the $ m $-labeled process. 
Let $\mathbf{X}^m = (X^i)_{i=1}^m $ and 
$ \mathsf{X}^{m*} = \sum_{m<j } \delta_{X^j} $. 
We call the pair $ (\mathbf{X}^m,\mathsf{X}^{m*} )$ the $ m $-labeled process. 
We shall present the Dirichlet form associated with the $ m $-labeled process 
$ (\mathbf{X}^m,\mathsf{X}^{m*} )$. 

Let $ \mu^{[m]} $ be the $ m $-Campbell measure of $ \mu $ defined as 
\begin{align}
\label{:90A}&
\mu^{[m]} (d\mathbf{x}d\mathsf{y}) = 
\rho^m (\mathbf{x}) \mu_{\mathbf{x}} ( d\mathsf{y})d\mathbf{x} 
,\end{align}
where $ \rho^m $ is the $ m $-point correlation function of $ \mu $ 
with respect to the Lebesgue measure $ d\mathbf{x} $ on $ \Sm $, and 
$ \mu_{\mathbf{x}} $ is the reduced Palm measure conditioned at $ \mathbf{x} \in \Sm $. 
Let 
\begin{align} \label{:90b} &
\E ^{\mu^{[m]}} (f,g) = \int_{\sS ^m \times \mathsf{S}} \Big\{
\frac{1}{2}\sum_{j=1}^{m} 
( \PD{f}{x_j}, \PD{g}{x_j})_{\Rd } + 
\DDD [f,g] \Big\} \mu^{[m]} (d\mathbf{x} d\mathsf{y})
,\end{align}
where $ \PD{}{x_j}$ is the nabla in $ \Rd $ and $ \DDD $ defined on $ \di $
is given by \eqref{:11b}, which is naturally 
regarded as the carr\'{e} du champ on $ \sS ^m \times \mathsf{S}$, and 
\begin{align} \label{:90c}&
\di ^{\mu ^{[m]}}=\{ f \in C_0^{\infty}(\sS ^m)\otimes \di \, ;\, 
\E ^{\mu^{[m]}} (f,f) < \infty , \, f \in L^2(\sS ^m \times \mathsf{S}, \mu^{[m]} )  \} 
.\end{align}

The closablity of the bilinear form $ (\E ^{\mu^{[m]}} , C_0^{\infty}(\sS ^m)\otimes \di ) $ 
on $ L^2(\sS ^m \times \mathsf{S}, \mu^{[m]} ) $ follows from \As{A1}--\As{A2}. 
We can prove this in a similar fashion as the case $ m =0 $ as \lref{l:51} (\cite{o.dfa}). 
We then denote by $ (\E ^{\mu ^{[m]}} ,\mathcal{D} ^{\mu ^{[m]}} )$ its closure. 
The quasi-regularity of $ (\E ^{\mu ^{[m]}} ,\mathcal{D} ^{\mu ^{[m]}} )$ 
is proved by \cite{o.tp} for $ \mu $ with bounded correlation functions (see \cite{mr} for quasi-regularity). 
The generalization to \As{A2} is derived by the same argument in the proof of Lemma \ref{l:25}. 
From \eqref{:90b}, we easily see  the Dirichlet form $ (\E ^{\mu ^{[m]}} ,\mathcal{D} ^{\mu ^{[m]}} )$ is local. 
When $ m=0 $, we interpret $ \mu^{[0]} = \mu $ and $ (\E ^{\mu ^{[0]}} ,\mathcal{D} ^{\mu ^{[0]}} )$ 
as the closure of $ (\E, \di^\mu) $,
which coincides with $ (\E^{\Os} ,\mathcal{D} ^{\Os} ) $.

Let $ \PPUm = \{\PPUm _{(\mathbf{s}^m,\mathsf{s}^{m*})} \} $ 
denote the diffusion measures associated with the 
$ m $-labeled Dirichlet form $ (\E ^{\mu^{[m]}} , \mathcal{D} ^{\mu ^{[m]}} ) $ 
on $ L^2(\sS ^m \times \mathsf{S}, \mu^{[m]} ) $. 
(see \cite{o.tp}). 
We quote: 
\begin{lemma}[{\cite[Theorems 2.4, 2.5]{{o.tp}}}] \label{l:90}
Assume that \As{A1}--\As{A4} hold. Let $ \lab $ be a label and let 
$ (\mathbf{X}^m, \mathsf{X}^{m*}) $ be the associated $ m $-labeled dynamics under $ \PPU $. Then 
\begin{align}\label{:90a}&
\PPUm _{(\mathbf{s}^{m},\mathsf{s}^{m*})} = 
\PPsU \circ (\mathbf{X}^m, \mathsf{X}^{m*})^{-1} 
.\end{align}
Here 
$ \mathsf{s}=\ulab ((\mathbf{s}^m,\mathsf{s}^{m*})) = 
\sum_{i=1}^m\delta_{s_i}+\mathsf{s}^{m*}$. 
\end{lemma}

Note that $ \PPsU $ in the right hand side is independent of $ m \in \N $. 
Hence \eqref{:90a} gives a sequence of coupled $ \Sm \ts \mathsf{S} $-valued continuous processes 
with distributions 
$ \PPUm _{(\mathbf{s}^m,\mathsf{s}^{m*})} $. 
In this sense, there exists a natural coupling among the $ m $-labeled Dirichlet forms 
$ (\E ^{\mu ^{[m]}} ,\mathcal{D} ^{\mu ^{[m]}}) $ on $ L^2(\sS ^m\ts \mathsf{S} ,\mu^{[m]} ) $. 
This coupling is a key point of the construction of {\ws}s of ISDE in \cite{o.isde}.

We can regard $ \mathbf{X}^m$ as a Dirichlet process of the diffusion 
$ (\mathbf{X}^m, \mathsf{X}^{m*})$ associated with the $ m $-labeled Dirichlet space. 
That is, one can regard 
$ A_t^{[\mathbf{x_m}]}:= \mathbf{X}_t^m - \mathbf{X}_0^m $ 
as a $ dm$-dimensional additive functional given by 
the composition of $ (\mathbf{X}^m, \mathsf{X}^{m*})$ with 
the coordinate function $ \mathbf{x}^m =(x_1,\ldots,x_m) \in (\Rd )^m$. 
Although $ \mathbf{X}^m$ can be regarded as an additive functional of the unlabeled process 
$ \mathsf{X}=\sum_i \delta _{X^i}$, $ \mathbf{X}^m$ 
is no longer a Dirichlet process in this case. 
Indeed, as a function of $ \mathsf{X} $, we cannot identify 
$ \mathbf{X}_t^m$ without tracing the trajectory of 
$ \mathsf{X}_s = \sum_i \delta_{X_s^i}$ for $ s \in [0,t]$. 

Once $ \mathbf{X}^m$ can be regarded as a Dirichlet process, 
we can apply the It$ \hat{\mathrm{o}}$ formula (Fukushima decomposition) and 
Lyons--Zheng decomposition to $ \mathbf{X}^m$, 
which are important in proving the results in the subsequent subsections.

\subsection{ISDE-representation: Logarithmic derivative}\label{s:23}
We next present the ISDE describing the labeled 
dynamics given by \lref{l:26}. The key notion for this is 
the logarithmic derivative of $ \mu $ to be introduced below. 

We set $ L_{\mathrm{loc}}^1(\sS \times \mathsf{S} ,\muone )= 
\bigcap_{r=1}^{\infty} L^1(\sS \times \mathsf{S} ,\muone _r )$, where 
$ \muone _r (\,\cdot\,)= \muone (\,\cdot\, \cap \Sr \times \mathsf{S} )$. 
Let $ \dib = \{ g \in \di ; \text{ $ g $ is bounded} \} $. 
We set 
\begin{align*}&
\dibone 
=\{ \sum_{i=1}^m f_i(x)g_i(\mathsf{y}) \, ;\, f_i\in 
C_{0}^{\infty}(\sS ),\ g_i \in \dib ,\, m \in \N 
\} 
.\end{align*}
We now recall the notion of the logarithmic derivative of $ \mu $ 
\cite{o.isde}. 
\begin{definition} \label{d:22}
An $ \Rd $-valued function 
$ \dmu \in L_{\mathrm{loc}}^1(\sS \times \mathsf{S} ,\muone )^d $ is called 
 {\em the logarithmic derivative} of $\mu $ if, for all 
$\varphi \in \dibone $, 
\begin{align}\label{:27a}&
\int _{\sS \times \mathsf{S} } \dmu (x,\mathsf{y})\varphi (x,\mathsf{y}) 
\muone (dx d\mathsf{y}) = - \int _{\sS \times \mathsf{S} } 
 \nabla_x \varphi (x,\mathsf{y}) \muone (dx d\mathsf{y}) 
.\end{align}
\end{definition}%

We make the following assumption:

\smallskip 

\noindent \As{A5}
$\mu $ has a logarithmic derivative $\dlog^\mu $. 
\smallskip 

\noindent 
The next lemma reveals the importance of logarithmic derivative. 
\begin{lemma}[\cite{o.isde}] \label{l:27}
Assume \As{A1}--\As{A5}. 
Let $ \mathbf{X}$ and $ \lab $ be as in \lref{l:26}. 
Assume 
\begin{align}\label{:27A}&
\bb = \frac{1}{2}\dmu 
.\end{align}
Then there exists $\mathsf{S}_ 0\subset \mathsf{S} $ such that $\mu(\mathsf{S}_ 0)=1$ 
and that the labeled dynamics $\mathbf{X} =(X^i)_{i\in\N }$ under 
$ \PPxU $ solves the ISDE for each $\xx \in \mathsf{S}_ 0 $ 
\begin{align}\label{:27b}
& dX_t^{i}=dB_t^i + \bb\XXXit dt, 
\quad i\in \N 
,\\ \label{:27c}
& \mathbf{X}_0=\labx 
,\end{align}
where $\mathbf{B}=(B^i)_{i\in\N}$ is the $(\R^d)^{\N}$-valued standard Brownian motion, and 
$ \mathsf{X}_t^{\diamond i }= \sum_{j\not=i}\delta_{X_t^{j}}$ for 
$ \mathbf{X}=(X^i)_{i\in\N}$. 
\end{lemma}

Let $\mathbf{X}=(X^i)_{i\in\N}$ be a {\ws}\ of ISDE \eqref{:27b} 
and denote by $ \mathsf{X}_t = \sum_{i=1}^{\infty} \delta_{X_t^i}$ 
the associated unlabeled process. Let $\mu_t$ be the distribution of 
$ \mathsf{X}_t $. 
We make the following assumptions on $ \mathbf{X}$. 

\smallskip

\noindent \As{$\mu $-AC} 
The $\mu $-absolutely continuous condition is satisfied. 
That is, if $ \mathsf{X}_0 = \mu $ in law, then 
\begin{align}\label{:28a}&
\mu_t \prec \mu \quad \mbox{ for all $t > 0$}
,\end{align}
where $\mu_t \prec \mu$ means $\mu_t$ is absolutely continuous with respect to $\mu $.

\smallskip 

\noindent \As{NBJ} 
The no-big-jump condition is satisfied. 
That is, if $ \mathsf{X}_0 = \mu $ in law, then 
\begin{align} \label{:28b} &
P( \IrT (\mathbf{X})<\infty )=1 \quad \text{ for each $ r, T\in \N$} 
,\end{align}
where $\IrT $ is the maximal label with which the particle intersects $S_r$ defined by
\begin{align}\label{:28c}&
\IrT (\mathbf{X}) = \sup \{i\in \N \,;\, |X_t^i|\le r \mbox{ for some $0\le t\le T$}\}
.\end{align}

\begin{lemma}\label{l:28}
Assume \As{A1}--\As{A5}. 
Then under $ \PPU = \{ \PPxU \}_{\xx \in\mathsf{S} }$ the labeled dynamics 
$\mathbf{X}$ satisfies the conditions \As{$\mu $-AC} and \As{NBJ}. 
\end{lemma}
\PF
Because the unlabeled dynamics $ \mathsf{X}$ is $ \mu $-reversible, 
$ \mu_t = \mu $ for all $ t $. Hence \As{$\mu $-AC} is obvious. 
The second claim follows from the Lyons--Zheng decomposition and \eqref{:28x} 
as we explain below.

We review a Lyons--Zheng decomposition recalling it 
from Theorem 5.7.3 in \cite[296p]{fot.2} in the form we need. 
While we believe the name Lyons--Zheng decomposition is now common, 
it was not used in \cite{fot.2}. 

We note that by \lref{l:92} no tagged particles $ \mathbf{X}=(X^i)$ explode. 
Let $ i $ be fixed and let $ m \in \N $ such that $ i \le m $. 
We use the Dirichlet form of $ m $-labeled process 
$ (\mathbf{X}^m,\mathsf{X}^{m*})$ given by \eqref{:90b}. 
Then the coordinate function $ x_i $, regarded as a function defined on 
$ \sS ^m \times \mathsf{S}$, is locally in the domain of the Dirichlet form. 
Then the additive functional $ A_t^{[x_i]} := X_t^i -X_0^i$ is a Dirichlet process. Hence one can apply the Fukushima decomposition to $ A_t^{[x_i]}$ 
to have 
\begin{align}& \notag 
A_t^{[x_i]} = M_t^{[x_i]} + N_t^{[x_i]}
,\end{align}
where $ M^{[x_i]}$ is the continuous local martingale additive functional and 
$ N^{[x_i]}$ is the continuous additive functional of locally zero energy. 
See Theorem 5.5.1 in \cite[273p]{fot.2} for detail. 
In our situation it is easy to see that $ M_t^{[x_i]} = B_t^i $, where $ B^i $ is 
a $ d $-dimensional Brownian motion. 

The Lyons--Zheng decomposition is a further decomposition of $ A_t^{[x_i]} $ 
consisting of two martingales. For $ A^{[x_i]}$, we have for each 
$ 0 \le t \le T $, 
\begin{align}\label{:94b}&
A_t^{[x_i]} = \frac{1}{2}M_t^{[x_i]} + \frac{1}{2} \hat{M}_t^{[x_i]} 
\quad \text{ for $ \PPUm _{\mu ^{[m]}} $-a.s.}
,\end{align}
where $ \hat{M}^{[x_i]}$ is the time reversal of $ M^{[x_i]}$ on $ [0,T]$ 
such that 
\begin{align*}&
\hat{M}_t^{[x_i]} = M_{T-t}^{[x_i]} (r_T)- M_{T}^{[x_i]} (r_T)
.\end{align*}
Here we set $ \map{r_T}{C([0,T]; \sS ^m \times \mathsf{S} )}
{C([0,T]; \sS ^m \times \mathsf{S} )}$ such that $ r_T(w )(t) = w (T-t)$. 
Such a decomposition is valid for general Dirichlet forms and 
their Dirichlet processes with functions locally in the domain of the Dirichlet form. This decomposition is called the Lyons--Zheng decomposition. 

In our case, because of the coupling in \lref{l:90} and subsequent argument, 
we have for each $ i \in \N $ 
\begin{align}
\label{:49a}&
\Xit - X _0 ^i 
=\half B_t^{i} + \half (B_{T-t}^{i} (r_{T}) - B_{T}^{i} (r_{T}) )
\quad \text{ for $ 0 \le t \le T $, $ \PPmU $-a.s. }
\end{align}
%
Hence from the coupling in \lref{l:90}, we see 
\begin{align}\label{:49b}
 \sum_{i=1}^{\infty} \PPmU ( \{\inf_{t\in[0,T]}|X_t^i| \le r \}) 
= &
 \int_{S\times \mathsf{S} } \mu ^{[1]}
(dx d\mathsf{s}) \PPUone _{(x,\mathsf{s})} ( \{\inf_{t\in[0,T]}|X_t^1| \le r \})
\\ \notag = &\PPUone _{\mu^{[1]}}( \{\inf_{t\in[0,T]}|X_t^1| \le r \})
.\end{align}
Then we have from \eqref{:49a}, the martingale inequality, and \eqref{:28x} 
\begin{align}\label{:44o}
 & \PPUone _{\mu ^{[1]}} ( \{\inf_{t\in[0,T]}|X_t^1| \le r \}) \le \, 
\PPUone _{\mu ^{[1]}} (\{\sup_{t\in[0,T]}|X_t^1- X_0^1 | \ge | X_0^1 | - r \}) 
\\ \notag \le \, &
\PPUone _{\mu ^{[1]}} (\{ \sup_{t\in[0,T]} |B_t^1| \ge | X_0^1 | - r \}) + 
\PPUone _{\mu ^{[1]}} (\{ \sup_{t\in[0,T]} |\ \hat{B}_t^1| \ge | X_0^1 | - r \}) 
\\ \notag = \, & 2 
\PPUone _{\mu ^{[1]}} (\{ \sup_{t\in[0,T]} |B_t^1| \ge | X_0^1 | - r \}) 
\\ \notag \le &2
 \int_{\sS \ts \mathsf{S} } 
\mathrm{Erf} ( \frac{| x |-r}{{\cref{;44}}T} ) \mu^{[1]} (dxd\mathsf{s}) 
 \\ \notag =&2 
\int_{\sS }\mathrm{Erf} ( \frac{|x|-r}{{\cref{;44}}T}) \rho ^1 (x) dx 
< \, \infty 
,\end{align}
where $ \Ct \label{;44}$ is a positive constant. Here we used \eqref{:49a} 
for the second line, the martingale inequality for the fourth line, and \eqref{:28x} for the last line. 
Hence putting \eqref{:49b} and \eqref{:44o} together and 
using Borel--Cantelli's lemma we have for each $ r \in \N $ 
\begin{align}&\label{:49h}
\PPmU ( \limsup_{i\to\infty }\{\inf_{t\in[0,T]}|X_t^i| \le r \}) =0 
.\end{align}
Note that by \eqref{:28c} 
\begin{align*}&
\PPmU ( \IrT (\mathbf{X})<\infty ) = 
1 - \PPmU ( \limsup_{i\to\infty }\{\inf_{t\in[0,T]}|X_t^i| \le r \}) 
.\end{align*}
Hence we obtain \As{NBJ} from this and \eqref{:49h}. 
\PFEND

%

\subsection{Finite systems in $ \SR $ of interacting Brownian motions with reflecting boundary condition} \label{s:24}

We give the SDE representation of the 
unlabeled process $ \mathsf{X}$ associated with 
the Dirichlet form $\9 $ on $ \Lm $. 
We denote by $ \PRl =\{\PRssfl \}$ the family of the diffusion measures given by $\9 $ on $ \Lm $. 
The Dirichlet form $\9 $ is dominated by $ (\E^{\Os},\mathcal{D} ^{\Os})$, that is, 
\begin{align}& \label{:29A}
 \9 \le (\E^{\Os},\mathcal{D} ^{\Os})
.\end{align}
Then from \eqref{:29A} we see that the capacity of $\9 $ is 
dominated by that of $ (\E^{\Os},\mathcal{D} ^{\Os})$. Hence non-collision of 
tagged particles under $ \PRl $ follows from that of the limit diffusion $\mathsf{X}$ given by $ (\E^{\Os},\mathcal{D} ^{\Os})$, which is assumed by \As{A3}. 
With the same reason, tagged particles under $ \PRl $ do not hit the set $ (\partial \sS )\cap \SR $. 
Non-explosion and non-entering of tagged particles under $ \PRl $ are obvious because they are reflecting diffusion on $ \SR $ and 
frozen outside $ \SR $. 

We now denote by $\mathbf{X}=(X^{i})_{i=1}^{\infty}$ 
the labeled process 
associated with $ \mathsf{X}$ and the label $ \lab $. Then 
$ \mathsf{X}=\{ \mathsf{X}_t \} $ is given by 
$ \mathsf{X}_t = \sum_{i=1}^{\infty} \delta_{X_t^{i}}$ 
from $ \mathbf{X}=(X^{i})_{i=1}^{\infty}$. 
By definition $ \mathbf{X}_0 = \lab (\mathsf{X}_0)$. 
The process $\mathbf{X}$ under $ \PRl $ 
describes the system of interacting Brownian motions in which
\begin{enumerate}
\item each particle in $S_R$ moves in ${S_R}$ and when it hits the boundary $\partial S_R$, it reflects and enter the domain $S_R$ immediately, 
\item the particles out of $S_R$ stay the initial positions forever.
\end{enumerate}
We denote by $\muRs $ the regular conditional distribution defined by
\begin{align}\label{:29a}&
\muRs (\,\cdot\,)=\mu (\, \cdot \, |\, \sigma [\pi_R^c] )(\sss)
\quad \text{ for $\mu $-a.s. $\sss$}
.\end{align}
Let $\mathsf{S} ^{\Rsss} = \{\mathsf{y}\in\mathsf{S} \,;\,\pi_R^c(\mathsf{y})=\pi_R^c(\sss)\}$. 
Then $\muRs $ is a probability measure supported on $ \mathsf{S} ^{\Rsss} $. 
Since we condition $\muRs $ outside $ \SR $, we regard $\muRs $ as a random point field on $ \SR $. 
Let $\muRsone $ be the 1-Campbell measure of $\muRs $. Then we have 
\begin{align}\label{:29b}&
\muRsone (dxd\mathsf{y})
=\rho^{\Rsss ,1}(x)\mu_x^{\Rsss}(d\mathsf{y})dx 
\quad \text{ for } x \in S_R 
,\end{align}
where $ \rho^{\Rsss ,1}$ is the one-point correlation function of 
$ \muRs $ and $ \mu_x^{\Rsss} $ is the reduced Palm measure of $ \muRs $ conditioned at $ x $. 
By the Green formula, we see for all $\varphi \in \dibone $, 
\begin{align}\label{:29c}
&- \int _{\sS _R \times \mathsf{S} } 
 \nabla_x \varphi (x,\mathsf{y}) 
 \muRsone (dx d\mathsf{y}) =
 \int _{\sS _R \times \mathsf{S} } 
 \dmu (x,\mathsf{y})\varphi (x,\mathsf{y}) 
 \muRsone (dx d\mathsf{y}) 
 \\ \notag 
 &\qquad\qquad \qquad + \int _{\partial\sS _R \times \mathsf{S} }
 \varphi (x,\mathsf{y})\mathbf{n}^R(x) 
\mathcal{S}_R(dx)\mu_x^{\Rsss}(d\mathsf{y})
,\end{align}
where $\mathcal{S}_R $ is the Lebesgue surface measure on the boundary 
$\partial \sS _R$ and $\mathbf{n}^R(x)$ is the inward normal unit vector at $x\in \partial \sS _R$. 
Hence for $\mu $-a.s. $\sss$ and for $\mu_x^{\Rsss}$-a.s. $\mathsf{y}$, 
the logarithmic derivative $ \mathsf{d}^{\Rsss} $ 
of $\muRs $ coincides with the sum of 
\begin{align}\notag 
& 
 \dmu(x, \pi_R(\mathsf{y})+ \pi_R^c(\sss))
\quad \text{ for }x\in\SR 
\end{align}
and a singular part associated with the the boundary $\partial \SR $. 
We then obtain informally 
\begin{align}\label{:29e}& 
 \mathsf{d}^{\Rsss} (x,\mathsf{y}) 
= 1_{\SR }(x)\dmu (x,\pi_R(\mathsf{y})+ \pi_R^c(\sss)) 
+ \mathbf{n}^R(x) 1_{\partial \SR }(x) \delta_{x}
.\end{align}
Here we naturally extend the domain of 
$ \mathsf{d}^{\Rsss} (x,\mathsf{y})$ to $ \sS \ts \mathsf{S} $ by taking 
$ \mathsf{d}^{\Rsss} (x,\mathsf{y}) = 0 $ for 
$ x \not\in \SR $. This is reasonable because particles outside 
$ \SR $ are fixed. 

By definition $ \xxSR $ coincides with the number of particles in 
$ \SR $ for a given configuration $ \xx $. 
From the Green formula \eqref{:29c} we see that 
$\mathbf{X} = (X^{i})_{i=1}^{\infty}$ is 
the system of infinite number of particles such that 
only particles in $ \SR $ move and satisfies the following SDE: 
For $\mu $-a.s. $\sss =\sum_{i}\delta_{s_i}$ and 
for $ \muRs $-a.s.\! $ \mathsf{x}=\sum_{i}\delta_{x_i}$ 
\begin{align}\label{:29m} 
d X_t^i =& dB_t^i + 
\half \dmu (X_t^i ,\mathsf{X}_t^{\diamond i })dt
 + 
\half \mathbf{n}^R(X_t^i )dL_t^{R,i} , 
\; 1\le i \le \xxSR , 
\\\label{:29n} 
dL_t^{R,i}=& \mathbf{1}_{\partial S_R}(X_t^i )dL_t^{R,i}
,\; 1\le i \le \xxSR 
,\\ \label{:29o}
X_t^{i}=& X_0^{i}, \quad i> \xxSR ,
\\ \label{:29p}
\mathbf{X}_0=& \labx 
,\end{align}
where $ \lab (\mathsf{s}) = (s_i)_{i=1}^{\infty}$, 
$ \lab (\mathsf{x}) = (x_i)_{i=1}^{\infty}$, 
$ \mathsf{X}_t^{\diamond i }= \sum_{j\not=i} \delta_{X_t^j }$, and 
$L^{R,i}= \{L_t^{R,i}\} $ are non-negative increasing processes; 
see for instance \cite{Chen93}. 
The particles outside $ \SR $ are frozen by \eqref{:29o}. 
Hence $ L_t^{R,i}=0$ for $ i > \xxSR $. 
We remark that $ s_i = x_i$ for all $ i > \xxSR $ 
for $ \mus $-a.s.\! $ \mathsf{x}$. 

Let $ \9 $ be the Dirichlet form introduced in \lref{l:22}. 
Then we can easily deduce from \As{A1} and \As{A2} that 
$ \9 $ is a local, quasi-regular Dirichlet form on 
$ \Lm $ and that there exists the associated diffusion $ \mathsf{X}$. The capacity for $ \9 $ is dominated by 
that for $ (\E ^{\Os},\mathcal{D} ^{\Os})$. Hence from \As{A3} we deduce that 
$ \mathsf{X}$ has also non-collision property. Clearly, each tagged particle of $ \mathsf{X}$ does not explode because of the definition of $ \ER ^{\La }$. 
We have thus obtained the labeled process $ \mathbf{X}$ 
from the unlabeled process $ \mathsf{X}$ and the label $ \lab $. 
To clarify the dependence on $ \9 $ we write $ \mathbf{X} = \XRl $. 
Using the Fukushima decomposition and taking \eqref{:29e} into account, we see that $ \mathbf{X} = \XRl $ is a {\ws}\ of SDE \eqref{:29m}--\eqref{:29p}.

Let $ \muRs $ be as \eqref{:29a}. Then 
\begin{align}\label{:29q}&
\mu ( \,\cdot\, ) = \int_{\mathsf{S} } \muRs ( \,\cdot\, ) 
\mu (d\sss )
\end{align}
and $ \muRs (\mathsf{S} ^{\Rsss} )=1$, where 
$\mathsf{S} ^{\Rsss} = \{\mathsf{y}\in\mathsf{S} \,;\,\pi_R^c(\mathsf{y})=\pi_R^c(\sss)\}$ 
as before. We set 
\begin{align}\label{:29r}&
\ER ^{\2 } (f,g) := \int_{\mathsf{S} } \DDDR [f,g] d\muRs 
= \int_{\mathsf{S} } \DDDR [f_{\Rsss },g_{\Rsss }] d\muRs 
,\\ \notag &
\di ^{\2 } = 
\{ f \in \di \cap L^2 (\mathsf{S} ,\muRs ) 
\, ;\, \ER ^{\2 } (f,f) < \infty \} 
,\end{align}
where for a function $ h $ on $ \mathsf{S} $, we define the function $ h_{\Rsss }$ on $ \mathsf{S} $ 
such that 
\begin{align}
\label{:29RR}&
 h_{\Rsss }(\cdot ) = 
h(\pi_R ( \cdot ) + \pi_R^c (\sss )) 
.\end{align}
The second equality in \eqref{:29r} is clear because 
for $ \muRs $-a.s.\,$ \mathsf{x}$ 
\begin{align}&\notag 
\DDDR [f,g](\mathsf{x})
=
\DDDR [f,g]_{\Rsss }(\mathsf{x})
=
\DDDR [f_{\Rsss },g_{\Rsss }](\mathsf{x}) 
.\end{align}
From \As{A1} we easily deduce that 
$ (\ER ^{\2 },\di ^{\2 })$ is closable on $ L^2 (\mathsf{S} ,\muRs )$. 
We then denote by $(\ER ^{\1 },\mathcal{D} _R^{\1 })$ 
the closure of $ (\ER ^{\2 },\di ^{\2 })$ on $ L^2 (\mathsf{S} ,\muRs )$. 
Furthermore, we see that $(\ER ^{\1 },\mathcal{D} _R^{\1 })$ is 
a local, quasi-regular Dirichlet form on $ L^2(\mathsf{S} ^{\Rsss} ,\muRs )$. 
Hence there exists a diffusion $ \mathsf{X}$ associated with 
$(\ER ^{\1 },\mathcal{D} _R^{\1 })$ on $ L^2(\mathsf{S} ^{\Rsss} ,\muRs )$. 
With the same reason for $ \9 $, we obtain the labeled dynamics 
$ \mathbf{X}$ from the unlabeled dynamics $ \mathsf{X}$. 
To clarify the dependence on $(\ER ^{\1 },\mathcal{D} _R^{\1 })$ 
we write $ \mathbf{X}^{\Rsss ,\La }$. 
Using the Fukushima decomposition, we deduce that 
the associated labeled diffusion is a {\ws}\ of the SDE \eqref{:29m}--\eqref{:29p} 
for $ \muRs $-a.s.\,$ \mathsf{x}$ for $ \mu $-a.s.\,$ \mathsf{s}$.


\bigskip

Let $\mathcal{T}(\mathsf{S} ) $ be the tail $\sigma$-field of the configuration space $\mathsf{S} $:
\begin{align}
\label{:30a}&
\mathcal{T}(\mathsf{S} )= \bigcap_{R=1}^\infty \sigma [\pi_R^c]
.\end{align}
Let $\mus $ be the regular conditional probability conditioned by the tail $\sigma$-field 
$ \mathcal{T}(\mathsf{S} )$ defined as
\begin{align}\label{:30b}&
\mus (\cdot)=\mu(\cdot | \mathcal{T}(\mathsf{S} ))(\sss)
.\end{align}
As $\mathsf{S} $ is a Polish space, such a regular conditional probability exists, and satisfies 
\begin{align}\label{:30c}&
\mu(\cdot)=\int_{\mathsf{S} }\mus (\cdot)\mu(d\sss)
.\end{align}
From the martingale convergence theorem we see that, for $\mu $-a.s.\! $\sss$, 
\begin{align}\label{:30z}&
\lim_{R\to\infty}\muRs (A) =\mus (A)
\quad \text{ for any } A \in \mathcal{B}(\mathsf{S} ) 
.\end{align}
This implies the weak convergence of $ \{ \muRs \} $ to $ \mus $ for $\mu $-a.s.\! $\sss$.

For given increasing sequences of natural numbers $ \mathbf{a}=\{ a(k) \}_{k=1}^{\infty} $ we set 
\begin{align}& \notag 
\mathsf{K}(\mathbf{a}) = 
\{ \mathsf{s}\in \mathsf{S} \, ;\, \sss (\sS _k )\le a(k) \ \text{ for all } k \} 
.\end{align}
Note that $ \mathsf{K}(\mathbf{a}) $ is compact and that, 
furthermore, a subset $ \mathsf{A}$ in $ \mathsf{S} $ is relatively compact if and only if 
there exist such a sequence $ \mathbf{a}$ satisfying 
$ \mathsf{A} \subset \mathsf{K}(\mathbf{a})$. 
Let $ \mathbf{a}[r]=\{ a[r] (k) \} $ 
be a family of increasing sequences of natural numbers such that 
\begin{align}& \notag
a[r] (k)< a[r+1] (k) \quad \text{ for all }k
.\end{align} 
Then $ \mathsf{K}(\mathbf{a}[r]) \subset 
\mathsf{K}(\mathbf{a}[r+1]) $. 
Because of the compactness criteria in $ \mathsf{S} $ as above and compact regularity of probability measures on Polish spaces, 
we can take $ \mathbf{a}[r]$ such that 
\begin{align}\label{:71h}&
\mu (\mathsf{K}(\mathbf{a}[r])^c)\downarrow 0 
\quad \text{ as } r \to \infty 
.\end{align}
Hence from \eqref{:30c} and the monotone convergence theorem 
\begin{align*}&
\int_{\mathsf{S} } \limi{r} \mus (\mathsf{K}(\mathbf{a}[r])^c) \mu (d\sss ) = 
\limi{r} \int_{\mathsf{S} } \mus (\mathsf{K}(\mathbf{a}[r])^c) \mu (d\sss ) = 
\limi{r}\mu (\mathsf{K}(\mathbf{a}[r])^c) = 0 
.\end{align*}
This implies $ \mus (\mathsf{K}(\mathbf{a}[r])^c) \downarrow 0 $ 
as $ r \to \infty $ for $ \mu $-a.s.\! $ \mathsf{s}$. 
For each $ R \in \N $ we similarly have 
$ \muRs (\mathsf{K}(\mathbf{a}[r])^c) \downarrow 0 $ 
as $ r \to \infty $ for $ \mu $-a.s.\! $ \mathsf{s}$. 
Then from these and \eqref{:30z} we obtain 
\begin{align}\label{:71i}&
\limi{r}\sup_{R \in\N } 
\muRs (\mathsf{K}(\mathbf{a}[r])^c) = 0 
\quad \text{ for $ \mu $-a.s.\! $ \mathsf{s}$}
.\end{align}
We now introduce a family of cut-off functions 
$ \{ \varpi _{\mathbf{a}[r]} \}_{r\in\N } $. 
\begin{lemma} \label{l:2X}
We set $ \mathbf{a}_+[r] = \{1 + a[r](k+1)\}_{k=1}^{\infty}$ for 
$ \mathbf{a}[r]=\{ a [r](k) \}_{k=1}^{\infty} $.
For each $ r $ there exists $ \varpi _{\mathbf{a}[r]} \in C_0(\mathsf{S} )$ 
satisfying the following: 
\begin{align}
\label{:2X1} & 
0 \le \varpi _{\mathbf{a}[r]} \le 1 
,\\
\label{:71k}&
\varpi _{\mathbf{a}[r]} (\sss ) = 
\begin{cases}
1 & \text{ for } \sss \in \mathsf{K}(\mathbf{a}[r]) \\
0 & \text{ for } \sss \in \mathsf{K}(\mathbf{a}_+[r]^c )
,\end{cases}
\\
\label{:91a}&
\DDD [\varpi _{\mathbf{a}[r]},\varpi _{\mathbf{a}[r]} ]
(\mathsf{s}) \le 2 \quad \text{ for all }\mathsf{s}\in\mathsf{S} 
.\end{align}
\end{lemma}

\PF
Let $ \theta \in C^{\infty}(\R )$ such that 
$ 0 \le \theta (t) \le 1 $ for all $ t \in \R $ and 
$ \theta (t) = 1 $ for $ t \le 0 $ and $ \theta (t) = 0 $ for $ t \ge 1 $. 
Furthermore, we assume $ |\theta ' (t)| \le 2 $ for all $ t $. 

Let $ \mathsf{s}= \sum_i \delta_{s_i}$. 
We take a label $ \lab = (\labi )$ such that $ |\labi (\mathsf{s})| \le |\labii (\mathsf{s})|$ for all $ i $. 
We set 
\begin{align}\label{:9b}&
\mathbf{d}_{\mathbf{a}[r]} (\mathsf{s})= 
\big\{ \sum_{k=1}^{\infty} \sum_{i\in J_{k,\mathsf{s}} (\mathbf{a}[r]) }
(k-|\labi (\mathsf{s})|)^2
 \big\}^{1/2} 
,\end{align}
where 
$ J_{k,\mathsf{s}} (\mathbf{a}[r]) = \{ i \, ;\, i > a[r](k) ,\, 
\labi (\mathsf{s}) \in \sS _k \} $. 
Let 
\begin{align}\label{:91d}& 
\varpi _{\mathbf{a}[r]} (\mathsf{s})
= \theta \circ \mathbf{d}_{\mathbf{a}[r]} (\mathsf{s})
.\end{align}
Then $ \varpi _{\mathbf{a}[r]} $ satisfies $ \varpi _{\mathbf{a}[r]} \in C_0(\mathsf{S} )$, 
\eqref{:2X1}, and \eqref{:71k}. 

A straightforward calculation shows 
\begin{align}\label{:91e} 
\DDD [\varpi _{\mathbf{a}[r]} ,\varpi _{\mathbf{a}[r]} ] 
 (\mathsf{s}) = & \half 
\Big\{ \frac{\theta' (\mathbf{d}_{\mathbf{a}[r]} (\mathsf{s}))}
{\mathbf{d}_{\mathbf{a}[r]} (\mathsf{s})} \Big\}^2 
 \sum_{k=1}^{\infty} \sum_{i\in J_{k,\mathsf{s}} (\mathbf{a}[r]) }
(k-|\labi (\mathsf{s})|)^2
\\ \notag = & \half 
\big(\theta' (\mathbf{d}_{\mathbf{a}[r]} (\mathsf{s}))\big)^2 \le 2 
.\end{align}
We thus see that $ \varpi _{\mathbf{a}[r]} $ satisfies \eqref{:91a}. 
\PFEND

Let 
$$
\mathcal{D}_{\bullet\bullet} =\{ f \varpi _{\mathbf{a}[r]} ;\, f \in \dimu ,\, r \in \N \} 
.$$ 
Then $ \mathcal{D}_{\bullet\bullet}$ is a subset of 
$ \Lmu $, $ \mathcal{D} _R^{\La} $, and $ \mathcal{D} _R^{\1 } $ for $ \mu $-a.s.\! $ \sss $ because 
$ \varpi _{\mathbf{a}[r]} \in C_0(\mathsf{S} )$ and satisfies \eqref{:2X1}--\eqref{:91a} by \lref{l:2X}. 

We set the $ \E _{R,1} ^{\La } $-norm of $ f $ by the square root of 
$ \ER ^{\La } (f,f)+ (f,f)_{L^2(\mathsf{S} ,\mu)}$, and 
the $ \E _{R,1}^{\1 } $-norm of $ f $ similarly. 
The next result shows the importance of $ \mathcal{D}_{\bullet\bullet} $. 
\begin{lemma} \label{l:2J}
$ \mathcal{D}_{\bullet\bullet}$ is dense in $ \Lmu $, 
$ \mathcal{D} _R^{\La} $, and $ \mathcal{D} _R^{\1 } $ with respect to the 
$ \Lmu $, $ \E _{R,1}^{\La}$, and $ \E _{R,1}^{\1 }$-norm for $ \mu $-a.s.\! $ \sss $, respectively. 
\end{lemma}

\PF 
We deduce from \eqref{:71h} that $ \mathcal{D}_{\bullet\bullet}$ is dense in $ \dimu $ 
with respect to the $ \Lmu $-norm. 
We have already proved that $ \dimu $ is dense in $ \Lm $ in the proof of \lref{l:25}. 
Hence we conclude that $ \mathcal{D}_{\bullet\bullet}$ is dense in $ \Lmu $. 

Recall that by definition 
\begin{align}& \label{:2Jz}
\E _R^{\La } (f,f) = \int_{\mathsf{S} } \DDDR [f,f] d\mu \quad \text{ for $ f \in \dimu $}
.\end{align}
We set $ \DDDR [f,f] = \DDDR [f]$. 
Then we have for $ f \in \dimu $ 
\begin{align}\notag 
\DDDR [ f-f \varpi _{\mathbf{a}[r]}] &
=\DDDR [f (1- \varpi _{\mathbf{a}[r]}) ] 
\\ \notag &
\le 
2 \big\{\DDDR [f] |1- \varpi _{\mathbf{a}[r]} | ^2 + 
|f|^2\DDDR [1- \varpi _{\mathbf{a}[r]} ] \big\} 
.\end{align}
From this combined with \eqref{:71h}, \eqref{:2X1}--\eqref{:91a}, and \eqref{:2Jz}, 
we see that $ \mathcal{D}_{\bullet\bullet}$ is dense in $ \mathcal{D} _R^{\La} $ 
with respect to the $ \E _{R,1}^{\La }$-norm. 
Replacing \eqref{:71h} by \eqref{:71i}, we deduce that 
$ \mathcal{D}_{\bullet\bullet}$ is dense in $ \mathcal{D} _R^{\1 } $ 
with respect to the $ \E _{R,1}^{\1 }$-norm for $ \mu $-a.s.\! $ \sss $ similarly. 
\PFEND
\begin{proposition}\label{l:2Y}
There exists a countable subset $ \mathcal{D}_{\bullet} $ of 
$ \mathcal{D}_{\bullet\bullet}$ such that $ \mathcal{D}_{\bullet} $ is dense in 
$ \Lmu $, $ \mathcal{D} _R^{\La} $, and $ \mathcal{D} _R^{\1 } $ for $ \mu $-a.s.\! $ \sss $ with respect to 
\7, respectively. 
\end{proposition}

\PF
A set $ A $ with pseudometric $ d $ is called separable if there exists 
a countable subset $ B \subset A $ being dense in $ A $ with respect to $ d $. 
In general, if $ A $ is separable with respect to $ d $, 
then any subset $ A_0 $ of $ A $ is also separable with respect to $ d $. 

As we see in \lref{l:2J}, $ \mathcal{D}_{\bullet\bullet}$ is a dense subset of 
$ \Lmu $, $ \mathcal{D} _R^{\La} $, and $ \mathcal{D} _R^{\1 } $ for $ \mu $-a.s.\! $ \sss $ 
with respect to the respective norms. 
Then it is enough for the claim to prove 
$ \mathcal{D}_{\bullet\bullet}$ is separable with respect to the 
$ \Lmu $, $ \E _{R,1}^{\La}$, and $ \E _{R,1}^{\1 }$-norm for $ \mu $-a.s.\! $ \sss $, respectively. 

 Let $ \dmuq = \dimu \cap \mathcal{B}_R $ and 
$ \mathcal{D}_{\bullet\bullet R} = \{ f \varpi _{\mathbf{a}[r]} ; f \in \dmuq , r \in \N \} $. 
We shall prove that the set $ \mathcal{D}_{\bullet\bullet R}$ is separable 
with respect to \7, respectively. 

For $ f \in \dR $ let $ f_R^m $ be the $ S_R^m $-representation of $ f $. 
Note that $ f_R^m $ is independent of $ \sss $ because $ f $ is $ \sigma [\pi_R ]$-measurable. 
We thus regard $ f_R^m $ as a function defined on $ S_R^m $. We see 
$ f_R^m \in C^{\infty}(S_R^m )$ from $ f \in \dR $. 
We can regard $ \dR $ as a subset of $\sum_{m=0}^{\infty} C^{\infty}(S_R^m )$ by the correspondence 
$ f \mapsto (f^m)_{m=0}^{\infty}$. 
Here by convention $ C^{\infty}(S_R^0 )$ is the set of constants corresponding to the functions defined on 
the subset $ \mathsf{S}_R^0 $. 
Let $ d_{R,N}$ ($ N\in \N \cup \{ \infty \} $) 
be the pseudometric on $\sum_{m=0}^{\infty} C^{\infty}(S_R^m )$ given by 
\begin{align*}&
d_{R,N} (f,g) = \sum_{m=0}^N \Big\{ \sup_{x \in \sS _R^m} |f_R^m(x)-g_R^m(x)| 
+ \sup_{x \in \sS _R^m} |\partial f_R^m(x)- \partial g_R^m(x)| \Big\}
.\end{align*}
Then $ \dR $ is separable with respect to the pseudometric $ d_{R,N}$ for each $ N \in \N $. 
Hence we take a countable dense subset $ \mathbb{F}_{R,N}$ of $ \dR $ 
with respect to $ d_{R,N}$. Let 
\begin{align}
\label{:2Yb}&
\widetilde{\mathbb{F}}_{R,N}= 
\{ f \varpi _{\mathbf{a}[r]} ; f \in \mathbb{F}_{R,N} , r \in \N \}
.\end{align}
Then 
$ \widetilde{\mathbb{F}}_{R,N} $ is a countable dense subset of 
$ \mathcal{D}_{\bullet\bullet R}$ with respect to $ d_{R,N}$. 
Let 
\begin{align}
\label{:2Yd}&
 \widetilde{\mathbb{F}}_R = \bigcup_{N\in\N }\widetilde{\mathbb{F}}_{R,N}
.\end{align}
Then $ \widetilde{\mathbb{F}}_R $ is a countable dense subset of 
$ \mathcal{D}_{\bullet\bullet R}$ with respect to $ d_{R,N}$ for each $ N \in \N $. 

Let $ f \varpi _{\mathbf{a}[r]}$ and 
$ f' \varpi _{\mathbf{a}[r']} $ be arbitrary elements of 
$ \mathcal{D}_{\bullet\bullet R} $. 
For each $ r , r' $ there exists $ N \in \N $ such that 
$ \varpi _{\mathbf{a}[r]} (\mathsf{s}) = \varpi _{\mathbf{a}[r']} (\mathsf{s}) = 0 $ on 
$ \sum_{m=N+1}^{\infty} \mathsf{S}_R^m $. 
Then clearly 
$$ 
 f \varpi _{\mathbf{a}[r]} (\mathsf{s}) = f' \varpi _{\mathbf{a}[r']} (\mathsf{s}) = 0 \quad \text{ on }
 \sum_{m=N+1}^{\infty} \mathsf{S}_R^m 
.$$
Hence 
\begin{align}
\label{:2Ya}&
d_{R,\infty}( f \varpi _{\mathbf{a}[r]} , f' \varpi _{\mathbf{a}[r']} ) = 
d_{R,N} ( f \varpi _{\mathbf{a}[r]} , f' \varpi _{\mathbf{a}[r']} )
.\end{align}
From \eqref{:2Yb}--\eqref{:2Ya} we see that $ \widetilde{\mathbb{F}}_R $ is dense in 
$ \mathcal{D}_{\bullet\bullet R} $ with respect to the metric $ d_{R,\infty} $. 

We note that, if $ R' \le R $, then the metric $ d_{R,\infty} $ on $ \mathcal{D}_{\bullet\bullet R} $ 
is stronger than the metrics given by 
\seven. 
Hence we deduce that $ \widetilde{\mathbb{F}}_R $ is dense in 
$ \mathcal{D}_{\bullet\bullet R} $ with respect to these metrics. 
We thus see that for any $ R' \le R $ the set $ \widetilde{\mathbb{F}}_R $ is a countable dense subset of 
$ \mathcal{D}_{\bullet\bullet R}$ with respect to the metrics given by 
\seven. 

Recall that the sets $ \mathcal{D}_{\bullet\bullet R} $ $ (R \in \N )$ are increasing in $ R $. 
Hence we see that 
\begin{align}
\label{:2Yc}&
 \mathcal{D}_{\bullet\bullet} = \bigcup_{R=R'}^{\infty} \mathcal{D}_{\bullet\bullet R} 
.\end{align}
We now take $ \mathcal{D}_{\bullet}:= \cup_{R=R'}\widetilde{\mathbb{F}}_{R} $. 
Then from the argument as above we see that, for any $ R \ge R'$, 
the set $ \mathcal{D}_{\bullet} $ contains a countable dense subset 
$ \widetilde{\mathbb{F}}_{R} $ in 
$ \mathcal{D}_{\bullet\bullet R}$ with respect to 
\seven. 
From this and \eqref{:2Yc} we see that 
$ \mathcal{D}_{\bullet} $ is a countable dense subset of 
$ \mathcal{D}_{\bullet\bullet}$ with respect to \seven. 
Retaking $ R'$ as $ R $ implies the claim. 
\PFEND

Let $ \XRl $ be the {\ws}\ of the SDE \eqref{:29m}--\eqref{:29p} given by the Dirichlet form 
$(\ER ^{\La },\mathcal{D} _R^{\La })$ on $ L^2 (\mathsf{S} ,\mu )$. 
Recall the disintegration of $ \mu $ given by \eqref{:29q}: 
\begin{align}\label{:29P}&
\mu = \int_{\mathsf{S} } \muRs \mu (d\sss ) 
.\end{align}
Taking \eqref{:29P} into account, we introduce a new {\ws}\ of \eqref{:29m}--\eqref{:29p}. 
For $ \mu $-a.s.\,$ \mathsf{s}$, we denote by $ \mathbf{X}^{\Rsss ,\La }$ the {\ws}\ of the SDE 
\eqref{:29m}--\eqref{:29p} given by the Dirichlet form 
$(\ER ^{\1 },\mathcal{D} _R^{\1 })$ on $ L^2 (\mathsf{S},\muRs )$. 

The next result provides a dynamical counterpart of the decomposition \eqref{:29P}. 
\begin{proposition} \label{l:29} 
For $ \mu $-a.s.\,$ \mathsf{s}$, $ \XRl = \mathbf{X}^{\Rsss ,\La }$ in distribution 
for $ \muRs $-a.s.\,$ \mathsf{x}$. 
\end{proposition}
\PF
Let $ T_{R,t}^{\La } $ and $ T_{R,t}^{\Rsss ,\La }$ be 
the semi-groups given by the Dirichlet forms 
$ \9 $ and $ (\ER ^{\1 },\mathcal{D} _R ^{\1 })$ on $ \Lm $ and $ L^2 (\mathsf{S} ,\muRs )$, respectively. 
Then to prove \pref{l:29} it is sufficient to prove the coincidence of these two semi-groups.

Let $ \mathcal{D}_{\bullet} $ be the set given by \pref{l:2Y}. 
Then $ \mathcal{D}_{\bullet}$ is a countable dense subset of 
$ \Lmu $, $ \mathcal{D} _R^{\La} $, and 
$ \mathcal{D} _R^{\1 } $ for $ \mu $-a.s.\! $ \sss $ 
with respect to $ \Lmu $, $ \E _{R,1}^{\La}$, and 
$ \E _{R,1}^{\1 } $-norm for $ \mu $-a.s.\! $ \sss $, respectively. 

From \eqref{:21c}, \eqref{:29q} and \eqref{:29r} we have for $ f,g \in \mathcal{D}_{\bullet} $ 
\begin{align}\label{:29s}&
\ER ^{\La } (f,g) = \int_{\mathsf{S} } \ER ^{\1 } 
(f_{\Rsss },g_{\Rsss }) \mu (d\sss)
.\end{align}
Here for a function $ h $ on $ \mathsf{S} $ let $ h_{\Rsss }$ be as \eqref{:29RR}. 
Then we see for $ f,g \in \mathcal{D}_{\bullet} $ 
\begin{align}\label{:29t}&
\int_{\mathsf{S} } f(\mathsf{s}) g (\mathsf{s})\4 -
\int_{\mathsf{S} } T_{R,t}^{\La } f(\mathsf{s}) g (\mathsf{s})\4 
=\int_0^t
\ER ^{\La }(T_{R,u}^{\La } f,g)du 
,\\ \label{:29tt}&
\int_{\mathsf{S} } f_{\Rsss }(\mathsf{x}) g_{\Rsss } (\mathsf{x})\3 -
\int_{\mathsf{S} }T_{R,t}^{\Rsss ,\La }
( f_{\Rsss }) (\mathsf{x}) g_{\Rsss } (\mathsf{x})
\3 
\\ \notag & \quad \quad \quad \quad \quad \quad \quad \quad \quad \quad \quad \quad \quad \quad \quad 
=\int_0^t
\ER ^{\1 } (T_{R,u}^{\Rsss ,\La } ( f_{\Rsss }) ,g_{\Rsss })du 
\end{align}
for $ \mu $-a.s.\! $ \sss $. 
We recall that $ \mathcal{D}_{\bullet} $ is countable. 
Hence, for $ \mu $-a.s.\! $ \sss $, \eqref{:29tt} holds for all 
$ f,g \in \mathcal{D}_{\bullet} $. 
Because the particles outside $ \SR $ are frozen, we easily see 
that for $ f,g \in \mathcal{D}_{\bullet} $ 
and for $ \mu $-a.s.\,$ \mathsf{s}$ 
\begin{align}\label{:29u}&
T_{R,t}^{\La }( f_{\Rsss })\, (\mathsf{x} ) =
(T_{R,t}^{\La } f )_{\Rsss }\, (\mathsf{x} ) 
&& \text{ for $\muRs $-a.s.\! $ \mathsf{x}$}
,\\ \label{:29v}&
T_{R,t}^{\Rsss ,\La } (f_{\Rsss })\, (\mathsf{x} ) =
(T_{R,t}^{\Rsss ,\La } f)_{\Rsss }( \mathsf{x}) 
&& \text{ for $\muRs $-a.s.\! $ \mathsf{x}$}
.\end{align}
In the following we write 
$ T_{R,t}^{\La } f_{\Rsss } = T_{R,t}^{\La }( f_{\Rsss }) $ and 
 $ T_{R,t}^{\Rsss ,\La } f_{\Rsss } = T_{R,t}^{\Rsss ,\La } (f_{\Rsss })$. 

From \eqref{:29q} and the definition $ h_{\Rsss }(\cdot ) = 
h(\pi_R ( \cdot ) + \pi_R^c (\sss )) $ we have 
\begin{align}\label{:29x}&
\int_{\mathsf{S} } f (\mathsf{s}) g (\mathsf{s}) \4 =
\int_{\mathsf{S} }\Big\{ \int_{\mathsf{S} } 
f_{\Rsss }(\mathsf{x}) g_{\Rsss } (\mathsf{x})\3 \Big\}\4 
.\end{align}
Replacing $ f $ with $ T_{R,t}^{\La } f $ in \eqref{:29x} and 
using \eqref{:29u} we have 
\begin{align}\label{:29xx}
\int_{\mathsf{S} } T_{R,t}^{\La } f(\mathsf{s}) g (\mathsf{s})\4 = &
 \int_{\mathsf{S} } \int_{\mathsf{S} } 
(T_{R,t}^{\La }f )_{\Rsss }(\mathsf{x}) g_{\Rsss } (\mathsf{x})
\3 \4 
\\ \notag 
= &
 \int_{\mathsf{S} } \int_{\mathsf{S} } 
 T_{R,t}^{\La }f _{\Rsss }(\mathsf{x}) g_{\Rsss } (\mathsf{x})
\3 \4 
.\end{align}
Note that 
$ \ER ^{\La }(T_{R,u}^{\La } f,T_{R,u}^{\La } f) \le 
\ER ^{\La }( f,f) $. 
With the same reason as \eqref{:29s} and from \eqref{:29u} 
we have for all $ f,g \in \mathcal{D}_{\bullet} $ 
\begin{align}\label{:29y}
\ER ^{\La }(T_{R,u}^{\La } f,g) &= 
\int_{\mathsf{S} } \ER ^{\1 } 
((T_{R,u}^{\La } f )_{\Rsss }, g _{\Rsss }) \4 
\\ \notag &
= 
\int_{\mathsf{S} } \ER ^{\1 } 
(T_{R,u}^{\La } f _{\Rsss },g _{\Rsss }) \4 
.\end{align}
Putting \eqref{:29x}--\eqref{:29y} into \eqref{:29t} and using 
the Fubini theorem we obtain 
\begin{align}&\notag 
\int_{\mathsf{S} }\Big\{ \int_{\mathsf{S} } f_{\Rsss }(\mathsf{x}) g_{\Rsss } (\mathsf{x})\3 \Big\}\4 -
\int_{\mathsf{S} }\Big\{ \int_{\mathsf{S} }T_{R,t}^{\La }
f_{\Rsss }(\mathsf{x}) g_{\Rsss } (\mathsf{x})
\3 \Big\}\4 
\\ \notag &
= \int_{\mathsf{S} }\Big\{ \int_0^t
\ER ^{\1 } (T_{R,u}^{\La } 
f_{\Rsss },g_{\Rsss }) du \Big\}\4 
.\end{align}
From this combined with \eqref{:29q}, \eqref{:29r}, and the definition of $ \DDDR $ we see 
for $ \mu $-a.s.\! $ \sss $ 
\begin{align}&\notag 
\int_{\mathsf{S} } f_{\Rsss }(\mathsf{x}) g_{\Rsss } (\mathsf{x})\3 -
\int_{\mathsf{S} }T_{R,t}^{\La }
f_{\Rsss }(\mathsf{x}) g_{\Rsss } (\mathsf{x})
\3 
=\int_0^t
\ER ^{\1 } (T_{R,u}^{\La } 
f_{\Rsss },g_{\Rsss })du 
.\end{align}
Comparing this with \eqref{:29tt}, we see that $ f \mapsto T_{R,t}^{\La }f_{\Rsss }$ is the semi-group 
associated with the Dirichlet form $(\ER ^{\1 },\mathcal{D} _R^{\1 })$ on $ L^2 (\mathsf{S} ,\muRs )$. 
We note here $ f = f_{\Rsss }$ for $ \muRs $-a.s. 
Therefore, for $ \mu $-a.s.\! $ \mathsf{s}$, 
\begin{align}\label{:29w}&
T_{R,t}^{\La } f (\mathsf{x}) = 
T_{R,t}^{\La } f_{\Rsss } (\mathsf{x}) = 
T_{R,t}^{\Rsss ,\La } f(\mathsf{x})
\quad \text{ for $\muRs $-a.s.\! $ \mathsf{x}$}
.\end{align}
Here we used 
$ \muRs (\mathsf{S} ^{\Rsss} )=1$, where 
$\mathsf{S} ^{\Rsss} = \{\mathsf{y}\in\mathsf{S} \,;\,\pi_R^c(\mathsf{y})=\pi_R^c(\sss)\}$ 
as before. 
%

The {\ws}s $ \XRl $ and 
$ \mathbf{X}^{\Rsss ,\La }$ are associated with 
the semi-groups $ T_{R,t}^{\La } $ and 
$ T_{R,t}^{\Rsss ,\La }$, respectively. 
Hence from \eqref{:29w} we deduce that 
these are equivalent in distribution. 
We thus see that the {\ws}s of SDE \eqref{:29m}--\eqref{:29p} given by these Dirichlet forms are the same. 
\PFEND

\section{Statements of the main results} \label{s:3}

In this section we present our main results. 
Let $ \bb \in L_{\mathrm{loc}}^1(\sS \times \mathsf{S} ,\muone )$ 
be a coefficient of ISDE \eqref{:31b} below. 
We introduce cut-off coefficients $\bb _{r,s,\p}$ of $ \bb $. 
Let $ C_b (S\times\mathsf{S} )$ be 
the set of all bounded continuous functions on $ S\times\mathsf{S} $. 
Then the main requirements for them are the following: 

\ms 

\noindent \As{A6} 
$\bb _{r,s,\p} \in C_b (S\times\mathsf{S} ) $ for each $ \rsp \in\N $ 
with $ r < s $, and 
\begin{align} \label{:31t} & 
\limi{r}\limi{s}\limi{\p }\sup_{R\ge r+s+1}
\parallel \bb _{r,s,\p} - \bbrs 
\parallel_{L_{\mathrm{loc}}^1(\sS \times \mathsf{S} ,\, \muRsone )}
 = 0 \quad \text{ for $ \mu $-a.s.\! $ \sss $}
,\end{align}
where $\muRsone $ is the one-Campbell measure of $\muRs $, and 
$ L_{\mathrm{loc}}^1(\sS \times \mathsf{S} ,\, \muRsone )$ is equipped with the semi-norm 
such that 
$ L^1(\sS _k \times \mathsf{S} ,\, \muRsone )$ for each $ k \in \N $.

We shall present a sufficient condition of \As{A6} in \sref{s:7}. 

Let $ P_{\labx }^{R}$ be the distribution of the {\ws}\ of 
SDE \eqref{:29m}--\eqref{:29p} given by the Dirichlet form 
$ \9 $ on $ \Lm $. 
The first main theorem of this paper is the following.
\begin{theorem}\label{l:31} 
Assume that \As{A1}--\As{A6} hold. 
Then the sequence $ \{P_{\labx }^{R} \}_{R\in\N }$ converges weakly in $\CS $ to $ \Pls $ for $ \mu $-a.s.\! $ \xx $, that is, 
for any $ F \in C_b(\CS )$ 
\begin{align}\label{:31a}&
\limi{R} \int_{\CS } F dP_{\labx }^{R} = 
\int_{\CS } F d\Pls 
.\end{align}
For $ \mu $-a.s.\! $ \xx $, the process 
$\XLa =(X^{i})_{i\in\N} $ under $ \Pls $ is a {\ws}\ of the ISDE 
\begin{align}\label{:31b}
& dX_t^{i}=dB_t^i + \bb\XXXit dt \quad \text{ for } i\in \N 
,\\ \label{:31c}
& \mathbf{X}_0=\labx 
\end{align} 
satisfying conditions \As{$\mu $-AC} and \As{NBJ}. 
Furthermore, $\XLa =(X^{i})_{i\in\N} $ under $ \Pls $ is 
associated with the resolvent $\{G^{\La}_{\alpha } \}$ 
of the Dirichlet form $(\E^{\La}, \mathcal{D} ^{\La})$ on $ \Lm $ 
in the following sense: 
\begin{align}\label{:31d}&
G^{\La}_{\alpha } (f) (\mathsf{x}) = 
\Els [\int_0^{\infty} e^{-\alpha t} f (\ulab (\mathbf{X}_t))dt]
,\end{align}
where $ f \in \Lm $, 
$ \ulab (\mathbf{X}_t) = \sum_{i=1}^{\infty} \delta_{X_t^{i}}$, and 
$ \Els $ is the expectation with respect to $ \Pls $. 
\end{theorem}

Because $(\E^{\La}, \mathcal{D} ^{\La})$ on $ \Lm $ is a Dirichlet form, there exists the associated Markovian semi-group on $ \Lm $ whose resolvent 
is $ G^{\La}_{\alpha } $ in \tref{l:31}. 
We have however not yet constructed the associated diffusion. Only a stationary Markov process is thus constructed at this stage. 
In general, we have to prove the locality and quasi-regularity of the Dirichlet form 
$(\E^{\La}, \mathcal{D} ^{\La})$ on $ \Lm $ for the existence of the associated diffusion. 

The next theorem establishes the existence of the associated diffusion by proving the identity between the Dirichlet forms 
 $(\E^{\La}, \mathcal{D} ^{\La})$ and $(\E^{\Os}, \mathcal{D} ^{\Os})$. 

We introduce another Dirichlet form $(\E^+, \mathcal{D} ^+)$. 
Recall that $ \bb = \frac{1}{2}\dmu $ by \eqref{:27A}, where 
$ \dmu $ is the logarithmic derivative of $ \mu $ defined by 
\eqref{:27a}. 
Put 
\begin{align} \notag 
\mathcal{D} ^+=& \Big\{ f\in \Lm ; \, 
\mbox{there exists $ f' \in L^2(S\times\mathsf{S} ,\muone)^d$ such that }
\\ \notag 
&\quad \quad - \int _{\sS \times \mathsf{S} } 
f(\delta_x+\mathsf{y})
\{ \nabla_x \varphi (x,\mathsf{y}) + \dmu (x,\mathsf{y})\varphi(x,\mathsf{y}) \} 
\muone (dx d\mathsf{y}) 
\\ \notag 
& \quad \quad \quad \quad =
 \int _{\sS \times \mathsf{S} } 
 f' (x,\mathsf{y})\varphi (x,\mathsf{y})
 \muone (dx d\mathsf{y}) 
\mbox{ for $\varphi\in 
C_{0}^{\infty}(\sS ^{\mathrm{int}} )\otimes \dib $}
\Big\}
.\end{align}
Let us denote the distributional derivative $f' $ by $D_x f$ and set 
\begin{align}&\notag 
\E^+ (f,g)= 
\int_{S\times \mathsf{S} } \frac{1}{2} 
D_x f \cdot D_x g \;\muone (dxd\mathsf{y}),
\quad f,g \in \mathcal{D} ^+
.\end{align}
We now state our second main theorem. 
\begin{theorem}\label{l:32} 
Assume \As{A1}--\As{A6}. 
Assume that a family of {\ws}s of ISDE \eqref{:27b}--\eqref{:27c} 
defined for $ \mu $-a.s.\! $ \mathsf{x}$ 
satisfying \As{$\mu $-AC} and \As{NBJ} 
are unique in law for $ \mu $-a.s.\! $ \mathsf{x}$. Then 
\begin{align}\label{:32c}&
(\E^{\Os},\mathcal{D} ^{\Os})=(\E^{\La},\mathcal{D} ^{\La})= (\E^+, \mathcal{D} ^+) 
.\end{align}
\end{theorem}

\smallskip 

$ \bullet $ If $ \Psi \in C_0^3 (\Rd )$ or hard-core potential of the form $ \Psi (x) = 1_{U}(x)$, where $ U = \{ |x| \le r \} $, 
and $ \Psi $ is a Ruelle-class interaction potential, 
then Lang \cite{lang.1,lang.2}, Fritz \cite{Fr}, 
Tanemura \cite{T2}, and others proved the pathwise uniqueness of stationary solutions for the associated grand canonical Gibbs measures (see \sref{s:85}). 
The potential $ \Psi (x) = 1_{U}(x)$ does not satisfies \As{A6}. 
It is plausible that \tref{l:31} can be generalized to be applicable to this type of potentials.

$ \bullet $ For sine$ _{\mathrm{\beta}}$ random point field 
$ \mu_{\mathrm{sine},\beta}$ with $ \beta \ge 1$, 
Tsai \cite{tsai.14} proved the pathwise uniqueness of {\ws}s 
and the existence of strong solutions 
for $ \mu_{\mathrm{sine},\beta}$-a.s.\!\! $ \mathsf{s}$. 
In \cite{o.isde}, the existence of {\ws}s whose unlabeled dynamics are $ \mu_{\mathrm{sine},\beta}$-reversible was proved 
for $ \beta = 1,2,4$. 
Combining these, 
we see that the assumptions in \tref{l:32} are fulfilled for 
$ \mu_{\mathrm{sine},\beta}$ with $ \beta = 1,2,4$ 
(see \sref{s:81}).


\section{Proof of \tref{l:31} } \label{s:4}

In this section we prove \tref{l:31}. 
Let $\mus $ be as in \eqref{:30b}. 

\begin{lemma} \label{l:41}
Assume that $\mu $ satisfies \As{A1}--\As{A6}. 
Then $\mus $ satisfies \As{A1}--\As{A6} 
for $\mu $-a.s.\! $\sss$.
\end{lemma}
\PF
This lemma follows from disintegration of $ \mu $ on $\mus $, and also the disintegration of their correlation functions and density functions, and Fubini's theorem. 
\PFEND

Let $\XRs =(\XRi )_{i=1}^{\infty}$ 
be the labeled diffusion process starting at $\labx $ 
whose unlabeled process is associated with the Dirichlet form 
$(\ER ^{\1 },\mathcal{D} _R^{\1 })$ introduced in \sref{s:24}. 
Although we wrote $ \XRs $ as $ \mathbf{X}^{\Rsss ,\La }$ in \sref{s:24}, 
we omit $ \mathrm{lwr}$ to simplify the notation. 
This omission yields no confusion because of \pref{l:29}. 

To clarify the dependence on $ R $ and $ \sss $, we write 
$\XRs =(\XRi )_{i=1}^{\infty}$ instead of 
$\mathbf{X} =(X^i )_{i=1}^{\infty}$ below. 
Recall that $ \xxSR $ equals the number of particles in $ \SR $. 
Suppose that $\xxSR \ge m$. We set the $m$-labeled process 
$ \mathbf{X}^{\Rsss,[m]}$ such that
\begin{align}\label{:41e}&
\mathbf{X}_t^{\Rsss,[m]}= 
(X_t^{\Rsss ,1},X_t^{\Rsss ,2},\dots,X_t^{\Rsss ,m}, 
\sum_{j=m+1}^\infty \delta_{\XRjt })
,\end{align}
where we freeze particles outside $ \SR $. 
Hence, $ \XRit = \XRiz $ for all $ t $ if $ i > \xxSR $. 
From \eqref{:41e} we have consistency such that, 
if we denote by $ X_t^{\Rsss,[m],i}$ the $ i $-th component of 
$\mathbf{X}_t^{\Rsss,[m]}$ from the beginning 
for $ 1 \le i \le m $ to clarify the dependence on $ m $, then 
\begin{align}&\notag 
X_t^{\Rsss,[m],i}=X_t^{\Rsss,[m+1],i} = \XRit 
\quad (i=1,2,\dots,m)
.\end{align}
It is known \cite{o.tp} that 
$\mathbf{X}^{\Rsss,[m]}$ is the diffusion process 
associated with the Dirichlet form 
\begin{align}&\notag 
\mathcal{E}^{\Rsss,[m]}(f,g)
=\int_{S_R^m\times\mathsf{S} }\{
\frac{1}{2}\sum_{i=1}^m \nabla_i f\cdot \nabla_i g
+\mathbb{D}_R[f,g] \}(\mathbf{x},\sss)
\mu^{\Rsss,[m]}(d\mathbf{x}d\sss)
\end{align}
on 
$ L^2(S^m\times\mathsf{S} ,\mu^{\Rsss,[m]}) $, 
where the domain $\mathcal{D} ^{\Rsss,[m]}$ 
is taken as the closure of 
\begin{align}&\notag 
\Big\{ f \in C_0^\infty(S^m)\otimes \di 
 ;\, 
\mathcal{E}^{\Rsss,[m]}(f,f) < \infty \Big\}
\cap L^2(S^m\times\mathsf{S} ,\mu^{\Rsss,[m]}) 
.\end{align}

We set $ f_i(\mathbf{x},\sss)=x_i\otimes 1 $. 
We can thus write for $ 1\le i \le m $
\begin{align}&\notag 
\XRit -\XRiz =
f_i(\mathbf{X}_t^{\Rsss,[m]})-f_i(\mathbf{X}_0^{\Rsss,[m]})=: 
A_t^{[f_i],[m]}
.\end{align}
Because the coordinate function $x_i=x_i\otimes 1$ belongs to 
$\mathcal{D} ^{\Rsss,[m]}$, 
$ A^{[f_i],[m]}$ is an additive functional of the $m$-labeled diffusion $ \mathbf{X}^{\Rsss,[m]}$ (see \cite{fot.2} for additive functional). 
We remark here that the $ m $-point correlation function of 
$ \muRs $ vanishes outside $ \SR $.

Applying the Fukushima decomposition to $ f_i $, 
the additive functional $A_t^{[f_i],[m]} $ can be decomposed as a sum of 
a unique continuous local martingale additive functional 
$M^{\Rsss,i}$ and an additive functional of zero energy 
$N^{\Rsss,i}$: 
\begin{align}&\notag 
\XRit -\XRiz =M_t^{\Rsss,i}+N_t^{\Rsss,i}
.\end{align}
We refer to \cite[Theorem 5.2.2]{fot.2} for the Fukushima decomposition. 

We recall another decomposition of $A_t^{[f_i],[m]} $ 
called the Lyons--Zheng decomposition \cite[Theorem 5.7.3]{fot.2}. 
Let $ \map{r_T}{C([0,T];\sS )}{C([0,T];\sS )}$ be such that 
$ r_T (X )_t = X_{T-t}$. 
Suppose that the distribution of $ \mathbf{X}_0^{\Rsss,[m]} $ 
is $ \mu^{\Rsss,[m]} $, or more generally, 
absolutely continuous with respect to $ \mu^{\Rsss,[m]}$. 
Then from the Lyons--Zheng decomposition we obtain 
\begin{align} \label{:41j}
\XRit -\XRiz 
&=\half M_t^{\Rsss,i} + \half (M_{T-t}^{\Rsss,i} (r_{T}) - 
M_{T}^{\Rsss,i} (r_{T}) )
\quad \text{ a.s}
.\end{align}

From \eqref{:29m} and $ f_i(\mathbf{x},\sss)=x_i\otimes 1 $ 
we see that $ M^{\Rsss,i} = B^i$ for $ 1 \le i \le \xxSR $, and hence
\eqref{:41j} becomes a simple form. 
That is, for $ 1 \le i \le \xxSR $ 
\begin{align}\label{:41k}&
\XRit -\XRiz 
=\half B_t^{i} + \half (B_{T-t}^{i} (r_{T}) - B_{T}^{i} (r_{T}) )
.\end{align}
For $ \xxSR < i < \infty $ we have 
$ \XRit = \XRiz = \labix $ by definition. 
Thus \eqref{:41k} is enough for our purpose. 
The decomposition \eqref{:41k} will be the main tool in this section. 

We set the maximal module variable 
$\XRmbar $ of the first $m$-particles by
\begin{align}\label{:42z}&
\XRmbar 
=\max_{1\le i\le m}\sup_{t\in [0,T]}|\XRit | 
.\end{align}
Throughout this section we fix $ T \in \N $. 
From \eqref{:41k} and \eqref{:42z} we obtain
\begin{lemma} \label{l:42} 
Assume that the distribution of $ \XRs _0 $ is $ \muRs \circ\lab^{-1}$. 
Then there exists a positive constant $ \Ct \label{;53}$ such that for $0\le t,u\le T$ 
\begin{equation}\label{:42b}
\sup_{R\in\N}\sum_{i=1}^m \EE [|\XRit -\XRiu |^4 ] 
\le \cref{;53} m |t-u|^2
.\end{equation}
Furthermore, for each $ m\in\N$ 
\begin{equation}\label{:42a}
\lim_{a\to\infty}\liminf_{R\to\infty} \PP 
(\XRmbar \le a)=1
,\end{equation}
and for each $r\in\N$ 
\begin{equation}\label{:42c}
\lim_{ l \to\infty} \inf_{R\in\N}
\PP (\IrT (\XRs )\le l )=1
,\end{equation}
where $ \IrT $ is defined by \eqref{:28c}. 
\end{lemma}
\PF 
From \eqref{:41k}, we obtain 
 \begin{align} \label{:42e}
 2 |\XRit - \XRiz | 
 & \le |B^i _t | + | B^i _{T-t} (r_{T}) - B^i _{T} (r_{T}) | 
 \quad \text{ a.s}
 .\end{align}
From \eqref{:42e} we easily obtain \eqref{:42b}. 

Recall that 
$ \labx =(\labix )_{i\in\N } \in \SN $ 
is a label. 
From \eqref{:30z} we obtain for $ A \in \mathcal{B}(\SN ) $ 
\begin{align}\label{:41c}&
\lim_{R\to\infty}\muRs \circ \lab ^{-1}(A)
=\mus \circ \lab ^{-1}(A) 
.\end{align}
Equation \eqref{:42a} follows straightforwardly from \eqref{:42e} and \eqref{:41c}. 

We deduce from \eqref{:42e} 
\begin{align}\label{:42f} 
&\PP \Big(\inf_{t\in[0,T]} |\XRit |\le r \Big) 
\le
\PP \Big(|\XRiz |-r \le \sup_{t\in[0,T]} |\XRit - \XRiz |
\Big) 
\\ \notag \le &
\PP \Big(2\{|\XRiz |-r \}\le 
\sup_{t\in[0,T]}
\{|B^i _t |+| B^i _{T-t}(r_{T}) - B^i _{T}(r_{T})|\}\Big) 
\\ \notag \le &
\PP \Big(|\XRiz |-r \le \sup_{t\in[0,T]} |B^i _t | \Big) 
+
\PP \Big(|\XRiz |-r \le 
\sup_{t\in[0,T]}
| B^i _{T-t}(r_{T}) - B^i _{T}(r_{T})|\Big) 
\\ \notag = & 2 
\PP \Big( | \labix |-r \le 
\sup_{t\in[0,T]} |B^i _t |\Big) 
\\ \notag \le & 4d 
\int_{\sS } \mathrm{Erf} (\frac{|x|-r}{\sqrt{T}}) 
\mu \circ (\labi )^{-1} (dx)
.\end{align}
Then we deduce from \eqref{:28c} and \eqref{:42f} that 
\begin{align}\label{:42g}
 \sup_{R\in\N} \PP \Big( \IrT (\XRs ) \ge l \Big) 
\le \,& 
 \sum_{i> l }^{\infty} \sup_{R\in\N}
\PP \Big(\inf_{t\in[0,T]} |\XRit |\le r \Big) 
\\ \notag \le \, & 
4d \sum_{i> l }^{\infty} 
\int_{\sS } \mathrm{Erf} (\frac{|x|-r}{\sqrt{T}}) 
\mu \circ (\labi )^{-1} (dx) 
.\end{align}
From \As{A4} we deduce 
\begin{align}\label{:42h}& 
\sum_{i=1}^{\infty}
\int_{\sS } \mathrm{Erf} (\frac{|x|-r}{\sqrt{T}}) 
\mu \circ (\labi )^{-1} (dx) = 
\int_{\sS } \mathrm{Erf} (\frac{|x|-r}{\sqrt{T}}) \rho^1(x)dx 
< \infty 
.\end{align}
From \eqref{:42g}--\eqref{:42h} we obtain \eqref{:42c}. 
\PFEND

From the conditions above we have the following lemma. 
\begin{lemma}\label{l:43} 
Make the same assumption as in \lref{l:42}. 
Then for each $ i , a , R \in \N $ such that $ i \le m $ 
\begin{align} \label{:43b}&
\PP \big( L_T^{\Rsss,i}= 0 \, ; \, \XRmbar \le a \big) = 1 
\quad \text{ for } a < R 
.\end{align}
\end{lemma}
\PF 
Recall that by \eqref{:29n} we have 
\begin{align}&\notag 
L_t^{\Rsss,i}= \int_0^t 
\mathbf{1}_{\partial S_R}(\XRiu )dL_u^{\Rsss,i}
.\end{align}
Then $ L^{R,i}=\{L_t^{\Rsss,i}\}$ is non-negative and increases only when $ \{ \XRit \} $ touches the boundary $ \partial \SR =\{ |x|=R \} $. 
Hence $ L_T^{\Rsss,i}= 0 $ for all $ a < R $ on $ \{ \XRmbar \le a \}$, 
which implies \eqref{:43b}. 
\PFEND

Let $ \bb_{\rsp }$ be as in \As{A6} and put 
\begin{align}\label{:45C}&
\mathsf{B}_{\rsp }^{\Rsss,i}(t)= \int_0^t
\bb_{\rsp }(X_u^{\Rsss,i}, \mathsf{X}_u^{\Rsss,\diamond i }) du
.\end{align}
We set for $ m \in \N $
\begin{align}& \notag 
\text{$ \mathbf{X}^{R\sss,m}=(X^{\Rsss,i})_{i=1}^m$,
$ \mathbf{B}_{\rsp }^{\Rsss,m}=
(\mathsf{B}_{\rsp }^{\Rsss,i})_{i=1}^m$, and 
$ \mathbf{L}^{R\sss,m}=(L^{\Rsss,i})_{i=1}^m$.}
\end{align}
Let $ \mathbf{X}^{R\sss}=(X^{\Rsss,i})_{i=1}^{\infty}$ and 
consider random variables 
\begin{align}\label{:45d}&
\mathbb{V}_{\rsp }^{\Rsss,m}=(\mathbf{X}^{\Rsss,m}, \mathbf{B}_{\rsp }^{\Rsss,m},
\mathbf{L}^{\Rsss,m})
,\\\label{:45e}&
\mathbb{W}_{\rsp }^{\Rsss}
=
\big((\mathbf{X}^{\Rsss,n}, \mathbf{B}_{\rsp }^{\Rsss,n},
\mathbf{L}^{\Rsss,n})_{n=1}^{\infty}, \mathbf{X}^{\Rsss }
\big)
.\end{align}
%
By construction, 
$ \mathbb{V}_{\rsp }^{\Rsss,m}$ and 
$ \mathbb{W}_{\rsp }^{\Rsss} $ are 
functionals of $ \mathbf{X}^{\Rsss }$. Hence we can regard 
$ \mathbb{V}_{\rsp }^{\Rsss,m}$ and 
$ \mathbb{W}_{\rsp }^{\Rsss} $ are defined on a common probability space. Let 
\begin{align}&\notag 
\8 = \inf \{0\le t \le T ; \max_{1\le i \le m}\big| X_t^{\Rsss ,i } \big| 
\ge a \} 
.\end{align}
Let $ \Xi ^m = C([0,T];S^m)\times C([0,T];\Rdm )^2 $ and 
$\Xi _0^m = C([0,T];S^m) \times \mathcal{BV} 
\times \mathcal{C}_+ $, where 
\begin{align}
&\notag 
\mathcal{BV} = 
\{ \eta =(\eta ^i)_{i=1}^{m}\in C([0,T];\Rdm );\, 
\eta \text{ is bounded variation}\} 
,\\
&\notag 
\mathcal{C}_+ = 
\{ \zeta =(\zeta^i)_{i=1}^{m} \in C([0,T];\Rdm );\, 
\zeta \text{ is non-decreasing}
\}
.\end{align}

We say a sequence of random variables is tight if for any subsequence 
we can choose a subsequence 
that is convergent in law. 
We also remark that tightness in $ C([0,T];\SN )$ for all $ T \in \N $ is equivalent 
to tightness in $\CS $ 
because we equip $\CS $ with a compact uniform norm. 
\begin{lemma}\label{l:45}
Make the same assumption as in \lref{l:42}. 
Then for $ \mu $-a.s.\! $ \mathsf{s}$, 
the following hold for all $ T\in\N $. 
\\
\thetag{1} 
$\{ \mathbb{V}_{\rsp }^{\Rsss,m} (\cdot \wedge \8 ) \}_{\rspRN }$ 
is tight in $C([0,T]; \Xi^m)$ for each $ m, a \in \N $. 
\\
\thetag{2} 
$\{ \mathbb{V}_{\rsp }^{\Rsss,m} \}_{\rspRN }$ 
is tight in $C([0,T]; \Xi^m)$ for each $m \in \N $. 
\\\thetag{3} 
$\big\{ \mathbb{W}_{\rsp }^{\Rsss}\big\}_{\rspRN }$ 
is tight in 
$\prod_{n=1}^{\infty}C([0,T]; \Xi^n) \ts C([0,T]; \SN )$.
\end{lemma}
\PF
We remark that tightness of 
$ \mathbb{V}_{\rsp }^{\Rsss,m} (\cdot \wedge \8 ) $ 
follows from that of each component 
$ \mathbf{X}^{\Rsss,m}(\cdot \wedge \8 )$, 
$ \mathbf{B}_{\rsp }^{\Rsss,m}(\cdot \wedge \8 )$, and 
$\mathbf{L}^{\Rsss,m}(\cdot \wedge \8 )$. 
Tightness of 
$ \{\mathbf{X}^{\Rsss,m}(\cdot \wedge \8 ) \}_{R\in\N}$ 
follows from \lref{l:42}. 
Tightness of 
$ \{\mathbf{L}^{\Rsss,m}(\cdot \wedge \8 ) \}_{R\in\N}$
follows from \lref{l:43}. 

Recall that $\bb _{r,s,\p} \in C_b (S\times\mathsf{S} ) $ by \As{A6}. 
Then tightness of 
$ \{ \mathbf{B}_{\rsp }^{\Rsss,m}(\cdot \wedge \8 ) \}_{\rspRN }$ follows from \eqref{:45C} with a straightforward calculation. 
We thus obtain \thetag{1}.

In general, a family of probability measures $ m_a $ in a Polish space is compact under the topology of weak convergence if and only if for any $ \epsilon > 0 $ there exists a compact set $ K $ such that $ \inf_a m_a (K) \ge 1- \epsilon $. Using this we conclude \thetag{2} from \thetag{1} combined with \eqref{:42a}. 

With the same reason as the proof of \thetag{1}, we obtain 
\thetag{3} from \thetag{1} and \thetag{2}. 
\PFEND

\lref{l:43} and \lref{l:45} imply that for any subsequence of 
$\big\{\mathbb{V}_{\rsp }^{\Rsss}(\cdot \wedge \8)\big\}_{\rspRN }$, 
$\big\{\mathbb{V}_{\rsp }^{\Rsss}\big\}_{\rspRN }$, and 
$\{ \mathbb{W}_{\rsp }^{\Rsss} \}_{\rspRN }$ 
there exist convergent-in-law subsequences, 
denoted by the same symbols, such that the following 
convergence in law holds: 
\begin{align}\label{:46a}&
\0 
\mathbb{V}_{\rsp }^{\Rsss,m} 
 (\cdot \wedge \8 ) 
= \big(\Xsm _a , \mathbf{B}_a^{\sss ,m},0, \mathbf{X}_a^{\sss } 
\big) \quad \text{ for each $m \in \N $}
,\\\label{:46a'}&
\0 
\mathbb{V}_{\rsp }^{\Rsss,m} 
 = \big(\Xsm , \mathbf{B}^{\sss ,m},0, \mathbf{X}^{\sss } \big) 
\quad \text{ for each $m \in \N $}
,\\\label{:46b}&
\0 
\mathbb{W}_{\rsp }^{\Rsss}
= \big((\Xsn , \mathbf{B}^{\sss ,n},0)_{n=1}^{\infty}, 
\mathbf{X}^{\sss } \big) 
.\end{align}
Here the subscript $ a $ in the right hand side of \eqref{:46a} denotes the dependence on $ a $. We note that the convergence 
$ \limi{R}\mathbf{L}^{\Rsss,m}(\cdot \wedge \8 ) =0$ 
follows from \lref{l:43}. 
From \lref{l:45} \thetag{3}, we have consistency: 
\begin{align}&\notag 
\Xsm = (X^{\sss ,1},\ldots,X^{\sss ,m})
.\end{align}
Here $ X^{\sss ,i}$ in the right hand side is the $ i $-th 
component of $ \mathbf{X}^{\sss }=(X^{\sss ,i})_{i=1}^{\infty}$. 
The same holds for $ \mathbf{B}^{\sss ,n }$ and 
we write $ \mathbf{B}^{\sss ,n}=(\mathsf{B}^{\si })_{i=1}^n $ 
This is the reason why we extend the state space in \thetag{3} 
of \lref{l:45} from that in \thetag{1} and \thetag{2}. 

We next check consistency in $ a $ in the limits in \eqref{:46a} and \eqref{:46a'}. 
Without loss of generality, we can assume 
\begin{align}\label{:46d}&
 P (\{ \Xsmbar = a \} ) = 0
.\end{align}
Indeed, if not, we can choose an increasing sequence 
$ \{ n(a) \}_{a\in\N} $ of positive numbers diverges to infinity 
such that $ P (\{ \XRmbar = n(a) \} ) = 0$ 
instead of $ \{ a \}_{a\in\N} $. 
Let 
\begin{align}&\notag 
\sigma_a^{\sss,m} = \inf\{0\le t \le T ; 
\max_{1\le i \le m}\big| X_t^{\si } \big| \ge a \}
.\end{align}
Then from \eqref{:46d} we deduce that 
the discontinuity points of the stopping time 
$ \sigma_a^{\sss,m}$ are of probability zero. 
Hence from convergence in \eqref{:46a} and \eqref{:46a'} we have 
\begin{align}\label{:46e}&
\big(\Xsm _a , \mathbf{B}_a^{\sss ,m},0,\mathbf{X}_a^{\sss } \big) 
(\, \cdot\, )=
\big(\Xsm , \mathbf{B}^{\sss ,m},0,\mathbf{X}^{\sss } \big) (\, \cdot \wedge \sigma_a^{\sss,m} ) 
.\end{align}

We set $ \XRidt = \sum_{j\not=i} \delta_{\XRjt }$ for 
$ \mathbf{X}^{\Rsss }=(\XRit )_{i=1}^{\infty}$. 
Using reversibility of diffusions, we obtain the following 
dynamic estimates from the static condition 
\As{A6}. 
\begin{lemma} \label{l:44}
Make the same assumption as in \lref{l:42}. 
Furthermore, we assume \As{A6}. Then 
for $ \mu $-a.s.\! $ \mathsf{s}$ and for each $i\in\N$ 
\begin{align}\label{:44b}
& 
\limi{r}\limi{s}\limi{\p }\sup_{R\ge r+s+1} 
\EE \Big[\int_0^T \5 
\big|\{\bb _{r,s,\p} -\bbrs \} \XXXRit \big|^{\phat }dt
\Big]=0
,\\\label{:44c}
 & 
 \limi{r}\limi{s}\limi{\p }
 \EE \Big[\int_0^T 1_{\Sr }(\Xsi ) 
 \big|\{\bb _{r,s,\p} -\bbrs \} \XXXsit \big|^{\phat }dt
\Big]=0
.\end{align}
\end{lemma}
\PF 
Let $\mathsf{X}^{\Rsss} $ be the unlabeled diffusion such that 
$ \mathsf{X}_t^{\Rsss} = \sum_{i=1}^{\infty} \delta_{\XRit }$. 
Because the diffusion $\mathsf{X}^{\Rsss}$ 
is associated with the Dirichlet form 
$(\ER ^{\1 }, \mathcal{D} _R^{\1 })$ 
introduced in \sref{s:24}, 
$\mathsf{X}^{\Rsss}$ is $ \muRs $-reversible. 
Then because of reversibility we have for all $ t $ 
\begin{align}\label{:44p}
&
\EE \big[ 
\5 \big|\{\bb _{\rsp }-\bbrs \}
\XXXRit \big|^{\phat }\big]
\\\notag \le &
\EE \big[\sum_{i=1}^{\infty}
\5 \big|\{\bb _{\rsp }-\bbrs \} 
\XXXRit \big|^{\phat }\big]
\\\notag = &
\EE \big[ \sum_{i=1}^{\infty}
1_{\Sr } (\XRiz ) \big|\{\bb _{\rsp }-\bbrs \}
\XXXRiz \big|^{\phat }\big]
\\\notag = & 
\int_{\mathsf{S} } \sum_{x_i\in\Sr }
1_{\Sr } (x_i ) \big| \{\bb _{\rsp }-\bbrs \}
(x_i, \sum_{j\not=i}^{\infty}\delta_{x_j})
\big|^{\phat }
\muRs (d\mathsf{x})
,\end{align}
where we set $ \mathsf{x}=\sum_i \delta_{x_i}\in \mathsf{S} $. 
Then we obtain \eqref{:44b} from \eqref{:31t} and \eqref{:44p}. 

Recall that 
$ \bb \in L_{\mathrm{loc}}^1(\sS \times \mathsf{S} ,\, \muone )$. 
Then 
$ \bb _{r,s,\p} -\bbrs \in 
L_{\mathrm{loc}}^1(\sS \times \mathsf{S} ,\, \muone )$. 
Hence 
\begin{align}&\notag 
\text{$ \bb _{r,s,\p} -\bbrs \in 
L_{\mathrm{loc}}^1(\sS \times \mathsf{S} ,\, \mu^{\sss ,[1]} )$ 
for $ \mu $-a.s.\! $ \sss $.}
\end{align}
From this and the martingale convergence theorem, we obtain from \eqref{:31t} that 
\begin{align} &\notag 
\limi{r}\limi{s}\limi{\p }
\parallel \bb _{r,s,\p} - \bbrs 
\parallel_{L_{\mathrm{loc}}^1(\sS \times \mathsf{S} ,\, \mu^{\sss ,[1]} )}
 = 0 
\quad \text{ for $ \mu $-a.s.\! $ \sss $}
.\end{align}
Then we can prove \eqref{:44c} in the same way as \eqref{:44b}. 
\PFEND

\noindent{\it Proof of \tref{l:31}. }
From \eqref{:42b} and \eqref{:42a} in \lref{l:42} 
$ \{P_{\labx }^{R} \}_{R\in\N }$ is tight in $\CS $. 
Let $ \Pls $ be an arbitrary limit point of $ \{P_{\labx }^{R} \}_{R\in\N }$. 
Define the map $  \upath $ from $\CS $ to the space of $ \sigma $-finite measure valued paths  defined by  
\begin{align*}& \upath (\mathbf{X})_t :=\ulab (\mathbf{X}_t)= \sum_i \delta_{X_t^i} 
\quad \text{ for $ \mathbf{X}=(X^i)$}
.\end{align*}
Then $ \upath (\mathbf{X}) = \{\upath (\mathbf{X})_t \}$ is the unlabeled process 
associated with $ \mathbf{X}$ by definition. 
By \pref{l:29} the distributions $ \{P_{\labx }^{R} \circ \upath ^{-1}\}_{R\in\N }$ 
of the unlabeled processes converge in finite-dimensional distributions. 
Hence the limit point $ \Pls \circ \upath ^{-1}$ is unique. 
From this we easily see that the limit point $ \Pls $ of $ \{P_{\labx }^{R} \}_{R\in\N }$ 
is unique. 
We therefore obtain \eqref{:31a}, the first claim of \tref{l:31}. 

The condition \eqref{:31c} is obvious from \eqref{:29p}. 
From \lref{l:23} and \pref{l:29} we obtain \eqref{:31d}. 
The condition \As{$\mu $-AC} is clear because 
 the limit unlabeled dynamics are associated with the symmetric Dirichlet form 
$ (\E^{\La},\mathcal{D} ^{\La})$ on $ \Lm $.

We next prove \eqref{:31b}. 
For $\psi\in C_0^\infty(S^m)$, 
let $F:\Xi _0^m \to C([0,T];\R)$ such that 
\begin{align}\label{:46f} 
F(\xi,\eta,\zeta)(t)&=
\psi(\xi(t))-\psi(\xi(0))-
\int_0^t\sum_{j=1}^m \nabla_j\psi(\xi(u)) \cdot d\eta^j(u) 
\\ & \notag 
- \int_0^t\sum_{j=1}^m \nabla_j\psi(\xi(u))\cdot\zeta^j(du) 
- \int_0^t \sum_{j=1}^m \frac{1}{2}\triangle_j\psi(\xi(u))du
.\end{align}
From It$\hat{\mathrm{o}}$-Tanaka formula, 
\eqref{:29m}--\eqref{:29o}, and $ \dmu = 2 \bb $, 
we deduce that for each $ m \in \N $ 
\begin{align}\label{:47e}& 
\sup_{R\ge r+s+1} E \Big[
\sup_{0\le t \le T }
\Big|
\FBL (t) -
\sum_{j=1}^m \int_0^{t} \nc dB _u^j 
\Big|^{\phat } \Big]
\\ \notag & \le \cref{;41}(\sss ,m,r,s,\p ) 
\big\{\sum_{j=1}^m \sup_{x\in\sS ^m} |\nabla_j \psi (x)|\big\}
,\end{align}
where we set 
\begin{align} \label{:47c} &
\cref{;41}(\sss ,m,r,s,\p )= 
\sup_{R\ge r+s+1} \sum_{i=1}^m
\EE \Big[\int_0^T \5 
\big|\{\bb _{r,s,\p} -\bbrs \} \XXXRit \big|^{\phat }dt\Big]
.\end{align}
We deduce from \eqref{:44b} and \eqref{:47c} 
that $ \Ct \label{;41} $ satisfy 
for $ \mu $-a.s.\! $ \mathsf{s}$ and for each $ m \in \N $ 
\begin{align}\label{:47d}&
\limi{r} \limi{s} \limi{\p } \cref{;41} (\sss ,m,r,s,\p ) = 0 
.\end{align}

Take $ \psi = \psi _Q \in C_0^{\infty}(\sS ^m )$ such that 
$ \psi _Q (x_1,\ldots,x_m) = x_i $ for $ \{ |x_i|\le Q \} $. 
Let $ a ,Q,R \in \N $ be such that $ a < Q , R $. 
Recall that $ \mathbf{L}_t^{\Rsss,m}=0 $ by \lref{l:43}. 
Then we deduce from \eqref{:46f} and It$\hat{\mathrm{o}}$-Tanaka formula 
that 
\begin{align}\label{:47o}&
\FBL (t \wedge \8 ) -
\sum_{j=1}^m \int_0^{t \wedge \8 }
 \nabla_j \psi _Q (\mathbf{X}_u^{R\sss ,m}) dB _u^j
\\ \notag &
 = 
X ^{R\si }(t \wedge \8 ) - 
X^{R\si }(0) - \BB _{\rsp }^{R\si } (t \wedge \8 ) - 
B^{i}_{t \wedge \8 },
\end{align}
where we write $Y_t= Y(t)$ for a stochastic process $ Y=\{ Y_t \} $.
We also remark that $ \{ B^{i} \}_{i=1}^{\infty} $ is 
$ (\Rd )^{\N }$-valued Brownian motion 
taken to be independent of $ R $. 

We write $ \mathbf{X}_a^{\sss ,m}=(X_a^{\sss ,i})_{i=1}^{\infty}$ and $ X_a^{\sss ,i} = \{ X_{a,t}^{\sss ,i} \} $. 
We set
\begin{align*}&
\lim_{r,s,\p ,R} = \0 
.\end{align*}
We have from \eqref{:46a}, \eqref{:47o}, \eqref{:47e}, and 
\eqref{:47d} that 
\begin{align}& \label{:47q}
E\big[ 
\sup_{0\le t \le T } \big|
X_{a,t}^{\si } - 
X_{a,0}^{\si } - \BB _{a,t}^{\si } - 
B_{t \wedge \sigmaA }^{i}
 \big|^{\phat } 
\big] 
\\ \notag &
=
\lim_{r,s,\p ,R}
E\big[ \sup_{0\le t \le T } \big| 
X ^{R\si }(t \wedge \8 ) - 
X^{R\si }(0) 
- \BB _{\rsp }^{R\si } (t \wedge \8 ) - 
B^{i}_{t \wedge \8 }
 \big|^{\phat } \big]
\\ \notag & =
\lim_{r,s,\p ,R}
E\big[ \sup_{0\le t \le T } \big| 
\FBL (t \wedge \8 ) 
 -
\sum_{j=1}^m \int_0^{t \wedge \8 }
 \nabla_j \psi _Q (\mathbf{X}_u^{R\sss ,m}) dB _u^j 
 \big|^{\phat } \big]
\\ \notag &
=0 
.\end{align}
Here the second and the third lines follow from \eqref{:46a} and \eqref{:47o}, 
respectively. We obtain the last line from \eqref{:47e} and \eqref{:47d}. 
From \eqref{:47q} we deduce 
\begin{align}\label{:47r}&
X_{a,t}^{\si } - X_{a,0}^{\si } - \BB _{a,t}^{\si } - 
B_{t \wedge \sigmaA }^{i}
 = 0 \quad \text{ for all $ t $}
.\end{align}
Then from \eqref{:46e} and \eqref{:47r} we have for all $ a \in \N $ 
\begin{align}\label{:47s}&
X_{t \wedge \sigmaA }^{\si } - 
X_0^{\si } - \BB _{t \wedge \sigmaA }^{\si } - 
B_{t \wedge \sigmaA }^{i} = 0 \quad \text{ for all $ t $}
.\end{align}
From $ P(\limi{a}\sigma_a^{\sss,m} = \infty )=1$, 
\eqref{:47s} implies 
\begin{align}\label{:47t}&
X_{t}^{\si } - X_0^{\si } - \BB _{t}^{\si } - B_{t}^i = 0 \quad \text{ for all $ t $}
.\end{align}
So it only remains to calculate the representation of $\BB ^{\si }$. 

We now recall 
$ \mathbf{B}_{\rsp }^{\Rsss,m}=
(\mathsf{B}_{\rsp }^{\Rsss,i})_{i=1}^m $
 and 
$ \mathsf{B}_{\rsp }^{\Rsss,i}(t)= \int_0^t
\bb_{\rsp }(X_u^{\Rsss,i}, \mathsf{X}_u^{\Rsss,\diamond i }) du $ 
by definition. 
We then deduce from \eqref{:45d}, \eqref{:46a'}, and 
\eqref{:44c} combined with 
$\bb _{r,s,\p} \in C_b (S\times\mathsf{S} ) $ that 
\begin{align}\label{:47u}
 \BB ^{\si }_t &= 
\lim_{r\to\infty} \lim_{s\to\infty} \limi{\mathsf{p}}\limi{R}
\mathsf{B}_{\rsp }^{\Rsss,i}(t)
&& \text{ by \eqref{:45d} and \eqref{:46a'}}
\\ \notag &=
\lim_{r\to\infty} \lim_{s\to\infty} \limi{\mathsf{p}}\limi{R}
\int_0^t 
\bb_{\rsp }(X_u^{\Rsss,i}, \mathsf{X}_u^{\Rsss,\diamond i }) du 
&& \text{ by definition}
\\ \notag &=
\lim_{r\to\infty} \lim_{s\to\infty} \limi{\mathsf{p}} 
\int_0^t
\bb_{\rsp }(X_u^{\si }, \mathsf{X}_u^{\sss,\diamond i }) du 
&& \text{ by $\bb _{r,s,\p} \in C_b (S\times\mathsf{S} ) $}
\\\notag &= 
\limi{r} \int_0^{t} 
1_{\Sr }(X_u^{\si })
\bb (X_u^{\si } ,\mathsf{X}_u^{\sss ,\diai }) du 
&& \text{ by \eqref{:44c}}
\\ \notag &= 
\int_0^{t} 
\bb (X_u^{\si } ,\mathsf{X}_u^{\sss ,\diai }) du 
&& \text{ in law}
.\end{align}
Putting \eqref{:47t}--\eqref{:47u} together yields 
\begin{align}& \label{:47w}&
X_{t}^{\si } - X_0^{\si } 
- \int_0^{t} \bb (X_u^{\si } ,\mathsf{X}_u^{\sss ,\diai }) du 
 - B_t^i = 0 
.\end{align}
From \eqref{:47w} we obtain \eqref{:31b}.

We finally prove \As{NBJ}. Let 
\begin{align}
\label{:47x}&
P_{\mu }^{\infty} = \int _{\mathsf{S} } \Pls \mu (d\mathsf{x})
.\end{align}
From \eqref{:31b} it is not difficult to see that for each $ l \in \N $
\begin{align}
\label{:47y}&
P_{\mu }^{\infty} (\partial \{ \IrT (\mathbf{X}) \le l \} ) = 0 
.\end{align}
Here $ \partial \{ \IrT (\mathbf{X}) \le l \}$ denotes the boundary of the set $ \{ \IrT (\mathbf{X}) \le l \}$. 
Let $ P_{\mu }^{R}$ be the distribution of $ \XRl$ with initial distribution $ \mu $. 
From \pref{l:29} and \eqref{:31a} we see that $ \{ P_{\mu }^{R} \} $ converge weakly to 
$ P_{\mu }^{\infty}$. 
Combining this with \eqref{:47y} we see 
\begin{align}\label{:47z}
P_{\mu }^{\infty} ( \IrT (\mathbf{X}) \le l ) &
= \limi{R} P_{\mu }^{R} (\IrT (\mathbf{X} )\le l ) 
\ge \inf_{R\in\N } P_{\mu }^{R} (\IrT (\mathbf{X} )\le l ) 
.\end{align}
Using this, \pref{l:29} and \eqref{:42c} we then obtain 
\begin{align}\notag %
P_{\mu }^{\infty} ( \IrT (\mathbf{X})<\infty ) = &\limi{l}
P_{\mu }^{\infty} ( \IrT (\mathbf{X}) \le l ) 
&& 
\\ \notag 
\ge &
\limi{l}\inf_{R\in\N } 
P_{\mu }^{R} (\IrT (\mathbf{X} )\le l ) 
&&\text{by }\eqref{:47z}
\\ \notag 
= &
\limi{l}\inf_{R\in\N } 
\int _{\mathsf{S} }\PP (\IrT (\XRs )\le l ) \mu (d\sss )
&&\text{by \pref{l:29}}
\\ \notag 
= &1 
&&\text{by \eqref{:42c} }
.\end{align}
This implies \As{NBJ}. 
We have thus completed the proof of \tref{l:31}. 
\qed

\section{Proof of \tref{l:32}} \label{s:5}

In this section we prove \tref{l:32}. 
Let $ \muRsone (dxd\mathsf{y})
=\rho^{\Rsss ,1}(x)\mu_x^{\Rsss}(d\mathsf{y})dx $ as in \eqref{:29b}. 
Let 
\begin{align} \notag 
&\mathcal{D} _{R}^{\Rsss ,+}=\{ f\in L^2(\mathsf{S} ,\muRs ) ; 
\mbox{there exists 
$\mathbf{f}_{\Rsss }\in L^2(\SR \times\mathsf{S} ,\muRsone )^d $ 
such that}
\\ \notag 
&\qquad \qquad - \int _{\SR \times \mathsf{S} } 
f(\delta_x+\mathsf{y})
\{ \nabla_x \varphi (x,\mathsf{y}) + 
\dmu (x,\mathsf{y})\varphi(x,\mathsf{y}) \} 
\muRsone (dx d\mathsf{y}) 
\\ \notag 
 &\qquad \qquad \qquad =
 \int _{\SR \times \mathsf{S} } 
 \mathbf{f}_{\Rsss } (x,\mathsf{y})\varphi (x,\mathsf{y})
 \muRsone (dx d\mathsf{y}) 
\mbox{ for 
$\varphi\in C_{0}^{\infty}(\SR ^{\mathrm{int}})\otimes \dib $}\}
.\end{align}
Denote $\mathbf{f}_{\Rsss } $ by $D_{\Rsss , x} f$. 
We introduce the bilinear form $ \ER ^{\Rsss ,+}$ with domain 
$ \mathcal{D} _{R}^{\Rsss ,+}$ such that 
\begin{align}\label{:50a}&
\ER ^{\Rsss ,+} (f,g)=\int_{\SR \times \mathsf{S} }
\half D_{\Rsss , x} f \cdot D_{\Rsss , x} g\, 
\muRsone (dxd\mathsf{y})
,\\\label{:50b}&
\dRcirc ^{\Rsss ,+} = \{ f \in \di \cap L^2(\mathsf{S} ,\muRs )\, ;\, 
\ER ^{\Rsss ,+} (f,f) < \infty \} 
.\end{align}

\begin{lemma} \label{l:50} 
For $ \mu $-a.s.\! $ \mathsf{s}$, the following hold. 
\\\thetag{1}
$ (\ER ^{\Rsss ,+}, \dRcirc ^{\Rsss ,+})$ is closable on 
$ L^2 (\mathsf{S} ,\muRs ) $ for each $ R \in \N $. 
\\\thetag{2} 
$ (\ER ^{\Rsss ,+},\mathcal{D} _R^{\Rsss ,+}) $ is the closure of 
$ (\ER ^{\Rsss ,+}, \dRcirc ^{\Rsss ,+}) $ on $ L^2 (\mathsf{S} ,\muRs ) $ 
 for each $ R \in \N $. 
\\\thetag{3} 
$(\ER ^{\1 },\mathcal{D} _R^{\1 }) 
= (\ER ^{\Rsss ,+},\mathcal{D} _R^{\Rsss ,+})$ for each $ R \in \N $. 
\end{lemma}
\PF
We easily see from \As{A1} that $ \muRsone $ 
has a labeled density 
$ m^{\Rsss,[1]} (x,y_1,\ldots,y_m)$ with respect to the Lebesgue measure on $ \SR ^{1+m}$ such that 
\begin{align}\label{:50c}&
\cref{;50a}^{-1} e^{-\Phi (x)-\sum_{i=1}^m \Phi (y_i) - 
\sum_{i=1}^m 
\Psi (x,y_i) - \sum_{i <i'}^m \Psi (y_i,y_{i'})}
\\ \notag & \quad \quad \quad \le m^{\Rsss,[1]} (x,y_1,\ldots,y_m)
\\ \notag & \quad \quad \quad \quad \quad \quad 
 \le \cref{;50a}
e^{-\Phi (x)-\sum_{i=1}^m \Phi (y_i) - \sum_{i=1}^m 
\Psi (x,y_i) - \sum_{i <i'}^m \Psi (y_i,y_{i'})}
,\end{align}
where $ \sss (\SR )=1+m $, $ m \in \N \cup\{ 0 \} $, 
$ \pi_R^c (\sss ) = \sum_{j=m+2}^{\infty} \delta_{s_j}$ for 
$ \lab (\sss ) = (s_i)_{i=1}^{\infty}$. Furthermore, 
$ \Ct \label{;50a}$ is a positive constant depending on 
$ (\Phi ,\Psi )$, $ \pi_R^c (\sss )$, and $ m $. 
For $ \mu $-a.s.\! $ \mathsf{s}$, 
$ (\ER ^{\Rsss ,+}, \dRcirc ^{\Rsss ,+})$ 
can be regarded as a form on 
$ L^2(\SR ^{1+m}, m^{\Rsss,[1]}dxdy_1\cdots dy_m$). 
Hence we deduce \thetag{1} and \thetag{2} from \eqref{:50c}. 
Once we regard $ (\ER ^{\Rsss ,+}, \dRcirc ^{\Rsss ,+})$ as a form in finite dimensions as above, \thetag{3} is obvious. 
\PFEND

We next introduce the Dirichlet form 
$(\ER ^+, \mathcal{D} _R^+)$. Let 
\begin{align} \label{:51a}
&\mathcal{D} _R^+=\{ f\in \Lm ; \mbox{there exists 
$\mathbf{f}_R\in L^2(\SR \times\mathsf{S} ,\muone)^d $ such that}
\\ \notag 
&\qquad \qquad - \int _{\SR \times \mathsf{S} } 
f(\delta_x+\mathsf{y})
\{ \nabla_x \varphi (x,\mathsf{y}) + \dmu (x,\mathsf{y})\varphi(x,\mathsf{y}) \} 
\muone (dx d\mathsf{y}) 
\\ \notag 
 &\qquad \qquad \qquad =
 \int _{\SR \times \mathsf{S} } 
 \mathbf{f}_R(x,\mathsf{y})\varphi (x,\mathsf{y})
 \muone (dx d\mathsf{y}) 
\mbox{ for 
$\varphi\in C_{0}^{\infty}(\SR ^{\mathrm{int}})\otimes \dib $}\}
.\end{align}
Let us denote $\mathbf{f}_R$ by $D_{R, x} f$, and 
set for $ f,g \in \mathcal{D} _R^+ $ 
\begin{align}
\label{:51b}&
\ER ^+ (f,g)=\int_{\SR \times \mathsf{S} }
\half D_{R, x} f \cdot D_{R, x} g\, \muone (dxd\mathsf{y})
,\\\label{:51c}&
\dRcirc ^{+}= 
 \{ f \in \di \cap \Lm \, ;\, \ER ^{+} (f,f) < \infty \} 
.\end{align}
\begin{lemma} \label{l:59}
\thetag{1} 
$ (\ER ^+, \dRcirc ^{+} )$ is closable on $ \Lm $. 
\\\thetag{2} 
$ (\ER ^+,\mathcal{D} _R^+)$ is the closure of 
$ (\ER ^+, \dRcirc ^{+} )$ on $ \Lm $. 
\\\thetag{3} 
$\9 = (\ER ^+,\mathcal{D} _R^+)$. 

\end{lemma}
\PF
On account of the disintegration \eqref{:29q}, we have 
\begin{align}\label{:59a}&
\muone (dxd\mathsf{y}) = 
\int_{\mathsf{S} }\muRsone (dxd\mathsf{y}) \mu (d\mathsf{s})
.\end{align}
From \eqref{:50a}--\eqref{:50b} and the Fubini theorem, we see that 
$ \dRcirc ^{+}\subset \dRcirc ^{\Rsss ,+} $ 
for $ \mu $-a.s.\! $ \mathsf{s}$. 
From \eqref{:50a}, \eqref{:51b}, and \eqref{:59a} 
we deduce for $ f,g \in \dRcirc ^{+}$ 
\begin{align}\label{:59b}&
\ER ^+ (f,g)=\int_{ \mathsf{S} } \ER ^{\Rsss ,+} (f,g) 
\mu (d\mathsf{s})
.\end{align}
Hence \thetag{1} follows from \lref{l:50} \thetag{1} and 
the argument in \cite[45-46pp]{o.rm}. 
These together with the similar argument in \sref{s:24} yield 
\thetag{2}. 
We obtain \thetag{3} from \eqref{:59b} combined with \eqref{:29s} and \lref{l:50} \thetag{3}. 
\PFEND

\begin{theorem}\label{l:51}
$(\E^{\La}, \mathcal{D} ^{\La})= (\E^+, \mathcal{D} ^+)$.
\end{theorem}
\PF 
Since $(\E^{\La},\mathcal{D} ^{\La})$ and $ (\E^+, \mathcal{D} ^+)$ 
are the increasing limits of 
$\{\9 \}$ and $ \{ (\ER ^+,\mathcal{D} _R^+) \} $ 
as $ R \to \infty$ respectively, we deduce \tref{l:51} from \lref{l:59} \thetag{3}.
\PFEND

\noindent{\it Proof of \tref{l:32}.} 
From \lref{l:27} and \lref{l:28} we see that the diffusion $\mathbf{X}^{\Os}$ 
such that $\mathbf{X}_0^{\Os}= \labx $ 
associated with the Dirichlet form $ (\E^{\Os},\mathcal{D} ^{\Os})$ 
on $ \Lm $ is a {\ws}\ of \eqref{:27b}--\eqref{:27c} satisfying 
\As{$\mu $-AC} and \As{NBJ}. 

By \tref{l:31} the process 
$\XLa $ under $ \Pls $ 
is a {\ws}\ of the ISDE \eqref{:31b}--\eqref{:31c} 
satisfying \As{$\mu $-AC} and \As{NBJ}. 

For $ \mu $-a.s.\! $ \mathsf{x}$, the {\ws}s of ISDE 
\eqref{:27b}--\eqref{:27c} satisfying \As{$\mu $-AC} and \As{NBJ} 
are unique in law by assumption. 
Hence for $ \mu $-a.s.\!\! $ \mathsf{x}$, 
$\mathbf{X}^{\Os}$ and $\XLa $ starting at $ \lab (\xx ) $ 
have the same distribution. 
Hence the associated semi-group coincides with each other. 
This together with \eqref{:31d} implies 
$(\E^{\Os},\mathcal{D} ^{\Os})=(\E^{\La},\mathcal{D} ^{\La}) $. 
From \tref{l:51} we have already obtained 
$(\E^{\La},\mathcal{D} ^{\La}) = (\E ^+,\mathcal{D} ^+) $. 
Combining these we complete the proof of \tref{l:32}. 
\qed


\begin{remark}\label{r:61}
We give a remark on a generalization of \tref{l:31} and \tref{l:32}. 
For this purpose we introduce a function 
$\mathsf{a} :S\times \mathsf{S} \to \R ^{d^2}$ and assume: 

\smallskip 

\noindent {\bf (B1)} 
$\mathsf{a} , \nabla _x \mathsf{a} \in C_b(S\times \mathsf{S} ) $, 
$ \mathsf{a} = {}^t\mathsf{a} $, and 
$\mathsf{a} $ is uniformly elliptic on $S\times\mathsf{S} $. 

\smallskip 

\noindent 
For $f,g \in \di $ we set 
$ \mathbb{D}_r^{\mathsf{a} ,m} [f,g] (\sss) = 0 $ 
for $ \sss\notin \mathsf{S}_ r^m $ and 
\begin{align}\label{:60c}
&\mathbb{D}_r^{\mathsf{a} ,m} [f,g] (\sss) = 
\frac{1}{2} 
\sum_{i=1}^m \mathsf{a} (s_i , \sss^{\diamond i }) 
 \nabla_{s_i} 
 f_{r,\sss}^m(\mathbf{x}_r^m(\sss)) \cdot \nabla_{s_i} 
 g_{r,\sss}^m(\mathbf{x}_r^m(\sss)) 
\quad \text{ for } \sss\in \mathsf{S}_ r^m 
,\end{align}
where we set 
$ \sss^{\diamond i } = \sum_{j\not=i} \delta _{s_j}$ 
for $ \sss=\sum_i \delta_{s_i}$. 
Moreover, we set 
\begin{align}
\label{:60d}&
\mathbb{D}_r^{\mathsf{a} } [f,g] (\sss) = \sum_{m=1}^{\infty}
\mathbb{D}_r^{\mathsf{a} ,m} [f,g] (\sss)
.\end{align}
If we replace the square fields 
$ \mathbb{D}_r^{m}$ and $ \mathbb{D}_r $ in 
\eqref{:21a} and \eqref{:21b} with 
$ \mathbb{D}_r^{\mathsf{a} ,m}$ and $ \mathbb{D}_r^{\mathsf{a} } $ in 
\eqref{:60c} and \eqref{:60d} and add assumption \As{B1}, then 
all results in \sref{s:3} still hold.
\end{remark}

\section{Construction of cut-off coefficients $ \mathsf{b}_{r,s,\p }$}\label{s:7}

In this section we present a sufficient condition of \As{A6}. 
We shall construct $ \mathsf{b}_{r,s,\p } $ in \As{A6}. 
For this purpose we prepare cut off functions. 

Let $ \chi _{t} \in C_0 (\Rd )$ such that $ 0 \le \chi _{t}\le 1 $ and that 
\begin{align}\label{:71d}&
\chi _{t} (x) =
\begin{cases}
1 & \text{ for }|x|\le t -1, \\
0 & \text{ for } |x| \ge t
.\end{cases}
\end{align}
Let $ \upsilon _{\p } \in C_b (\Rd )$ such that $ 0 \le \upsilon _{\p } \le 1 $ and that 
\begin{align}\label{:71e}&
\upsilon _{\p } (x) =
\begin{cases}
0 & \text{ for }|x|\le 1/\p,\\
1 & \text{ for } 2/\p \le |x| < \infty 
.\end{cases}
\end{align}
Let $ \{ \varpi _{\mathbf{a}[r]} \}_{r\in\N } $ be the cut-off function given by \lref{l:2X}. 
Then from \lref{l:2X} we have $ \varpi _{\mathbf{a}[r]} \in C_0(\mathsf{S} )$ and 
\begin{align}
\label{:71a} & 
0 \le \varpi _{\mathbf{a}[r]} \le 1 
,\quad 
\varpi _{\mathbf{a}[r]} (\sss ) = 
\begin{cases}
1 & \text{ for } \sss \in \mathsf{K}(\mathbf{a}[r]) \\
0 & \text{ for } \sss \in \mathsf{K}(\mathbf{a}_+[r]^c )
,\end{cases}
.\end{align}
For $ \chi _{s}$, $ \upsilon _{\p } $, 
$ \Psi \in C^{1}(\Rd \backslash \{ 0 \} )$, and constants $ \varrho _s $ we set 
\begin{align*}&
\bb _{s,\p } (x,\mathsf{y}) = 
\frac{\beta}{2}\Big( 
\Big\{\sumii{i}
\chi _{s}(x-y_i)\upsilon _{\p } (x-y_i) 
\nabla \Psi (x-y_i)\Big\}
 - \varrho _s \Big) 
,\end{align*}
where we set $ \mathsf{y}=\sum_i\delta_{y_i}$ as before. 
We give a sufficient condition of \As{A6}. 
\begin{lemma} \label{l:71}
Assume that for some $ \p > 1 $
\begin{align}\label{:71g}&
\bb (x,\mathsf{y}) = \limi{s}\limi{\p } \bb _{s,\p } (x,\mathsf{y})
\quad \text{ in $ L_{\mathrm{loc}}^p(\sS \times \mathsf{S} ,\muone )$}
.\end{align} 
Let 
\begin{align}&\label{:71l}
\mathsf{b}_{r,s,\p }(x,\mathsf{y} ) = 
 \chi _{r}(x) \varpi _{\mathbf{a}_+[r]} (\mathsf{y}) \bb _{s,\p } (x,\mathsf{y})
.\end{align}
Then $ \mathsf{b}_{r,s,\p }$ $( r,s,\p \in \N )$ satisfy \As{A6}. 
\end{lemma}

\PF 
From $ \chi _{r} \in C_0(\Rd )$ and $ \varpi _{\mathbf{a}_+[r]} \in C_0(\mathsf{S} )$ we have 
$ \chi _{r} \otimes \varpi _{\mathbf{a}_+[r]} \in C_0(\Rd \ts \mathsf{S} )$. 
From this, \eqref{:71e}, and $ \Psi \in C^{\infty}(\Rd \backslash \{ 0 \} )$ 
we obtain $ \mathsf{b}_{r,s,\p }\in C_0 (\sS \times \mathsf{S} )\subset C_b (\sS \times \mathsf{S} )$. 

We next check \eqref{:31t}. Let 
\begin{align}\label{:71mm}
\mathsf{b}_{r}(x,\mathsf{y} ) := &
 \chi _{r}(x) \varpi _{\mathbf{a}_+[r]} (\mathsf{y}) 
 \bb (x,\mathsf{y}) 
.\end{align}
Then from \eqref{:71d}, \eqref{:71e}, and \eqref{:71g} the functions 
$ \mathsf{b}_{r}(x,\mathsf{y} )$ satisfy 
\begin{align}
\label{:71m}&
\bb (x,\mathsf{y}) - \bb _r (x,\mathsf{y}) = 
\{ 1 - \chi _{r}(x) \varpi _{\mathbf{a}_+[r]} (\mathsf{y}) \} \bb (x,\mathsf{y}) 
.\end{align}

Let $ p'$ be such that $ 1/p+1/p'=1$. 
Then the H\"{o}lder inequality and \eqref{:71a}--\eqref{:71m} combined with 
$ 0 \le \chi _{r}(x) \varpi _{\mathbf{a}_+[r]} (\mathsf{y}) \le 1 $ imply 
\begin{align}&\notag 
\parallel \bb - \bb _r \parallel_{L^1(\sS _k \times \mathsf{S} ,\, \muRsone )} 
\le 
\muRsone 
\big(\sS _k \ts \{\mathsf{S} \backslash \mathsf{K}(\mathbf{a}[r] ) \} \big) ^{1/p' }
\parallel \bbrs \parallel_{L^p(\sS _k \times \mathsf{S} ,\, \muRsone )}
.\end{align}
Similarly as \eqref{:71i}, we have for $ \mu $-a.s.\! $ \mathsf{s}$ 
\begin{align}&\notag 
\limi{r}\sup_{R \in\N } 
\muRsone 
\big(\sS _k \ts \{\mathsf{S} \backslash \mathsf{K}(\mathbf{a}[r] ) \} \big)=0 
.\end{align}
By the martingale convergence theorem we deduce for $ \mu $-a.s. $ \sss $
\begin{align}&\notag 
\limi{R}
\parallel \bbrs \parallel_{L^p(\sS _k \times \mathsf{S} ,\, \muRsone )}
=
\parallel \bbrs \parallel_{L^p(\sS _k \times \mathsf{S} ,\, \mu ^{\sss, [1]} )}
.\end{align}
From these we deduce 
\begin{align}\label{:71p}&
\limi{r}\limi{s}\limi{\p }\sup_{R\ge r+s+1}
\parallel \bb - \bb _r 
\parallel_{L_{\mathrm{loc}}^1(\sS \times \mathsf{S} ,\, \muRsone )}
 = 0 \quad \text{ for $ \mu $-a.s.\! $ \sss $}
\end{align}
From \eqref{:71d}--\eqref{:71m} and the martingale convergence theorem we deduce similarly 
\begin{align}\label{:71q}&
\limi{r}\limi{s}\limi{\p }\sup_{R\ge r+s+1}
\parallel \bb _r - \bb _{r,s,\p} 
\parallel_{L_{\mathrm{loc}}^1(\sS \times \mathsf{S} ,\, \muRsone )}
 = 0 \quad \text{ for $ \mu $-a.s.\! $ \sss $}
.\end{align}
From \eqref{:71p} and \eqref{:71q} we obtain \eqref{:31t}. 
We thus see $ \mathsf{b}_{r,s,\p }$ $( r,s,\p \in \N )$ satisfy \As{A6}. 
\PFEND

\begin{remark}\label{r:71}
The construction of $ \mathsf{b}_{r,s,\p }$ as above is robust and can be applied to all examples in this paper. We will take $ \Psi $ as an interaction potential with singularity at the origin such as the logarithmic potential. 
In general, the interaction $ \Psi $ is long ranged, and the sum in \eqref{:71g} makes sense only at the level of conditional convergence. 
\end{remark}

\section{Examples} \label{s:8}

 In this section we present examples satisfying the assumptions of \tref{l:31} and \tref{l:32}. 
The first two examples satisfy the assumptions of both \tref{l:31} and \tref{l:32}. 
Then the following identity holds for these two examples by \tref{l:32}. 
\begin{align}\label{:70a}&
(\E^{\Os},\mathcal{D} ^{\Os})=(\E^{\La},\mathcal{D} ^{\La})= (\E^+, \mathcal{D} ^+)
.\end{align}
Other examples satisfy the assumptions of \tref{l:31}. Hence the limit stochastic dynamics $ \mathbf{X}$ 
satisfy the associated ISDEs \eqref{:31b}. 
We write $L_{\mathrm{loc}}^p (\muone ) = L_{\mathrm{loc}}^p (\sS \times \mathsf{S} ,\muone )$, 
where $ \muone $ is the one-Campbell measure of $ \mu $ as before.

\subsection{Sine$_\beta$ interacting Brownian motion $(\beta=1,2,4)$} \label{s:81}

Let $ d = 1 $ and $ \sS = \mathbb{R}$. 
Let $ \mu _{\mathrm{sin},\beta }$ be a sine$_{\beta }$ random point field 
\cite{Meh04,forrester}, where $ \beta = 1,2,4$. 
By definition, $\mu _{\mathrm{sin},2 }$ is the random point field on 
$ \mathbb{R}$ with $ n $-point correlation function 
$ \rho_{\mathrm{sin},2 }^{n} $ 
with respect to the Lebesgue measure given by 
\begin{align}&\notag 
 \rho_{\mathrm{sin},2 }^{n} (x_1,\ldots ,x_n ) = 
 \det [\mathsf{K}_{\mathrm{sin},2 } (x_i,x_j)]\ijn 
.\end{align}
Here 
$ \mathsf{K}_{\mathrm{sin},2 } (x,y) = \sin \pi (x-y)/\pi (x-y)$ 
is the sine kernel. 
$ \musinone $ and $ \musinfour $ are also defined by correlation functions 
given by quaternion determinants \cite{Meh04}. 
The random point fields $ \mu _{\mathrm{sin},\beta }$ ($ \beta =1,2,4$) satisfy 
\As{A1}--\As{A3} \cite{o.rm,o.col}. 
$ \mu _{\mathrm{sin},\beta } $ clearly satisfy \As{A4} 
because their one-point correlation functions are constant. 

Let $ 1 < p < 2 $. Then \As{A5} is satisfied with the logarithmic derivative given by
\begin{align}\label{:81b}&
 \bb (x,\mathsf{y}) = \frac{\beta}{2} 
\limi{r} \sum_{|x-y_i|<r} \frac{1}{x-y_i}
\quad \text{ in }
L_{\mathrm{loc}}^p (\mu _{\mathrm{sin},\beta }^{[1]})
.\end{align}
Here $ \mathsf{y}=\sum_i\delta_{y_i}$ and 
\lq\lq in 
$ L_{\mathrm{loc}}^p (\mu _{\mathrm{sin},\beta }^{[1]})$'' means convergence in 
$ L^p(\sS _k \ts \mathsf{S},\, \mu _{\mathrm{sin},\beta }^{[1]} )$ for all $ k \in \N $. 
The formula \eqref{:81b} was proved in \cite{o.isde} for $ \beta =1,2,4$. 
The definition of conditional convergence in \eqref{:71g} and \eqref{:81b} are slightly different. 
For $ \beta = 1,2,4$ a representation of correlation functions of 
$ \mu _{\mathrm{sin},\beta } $ are known. 
Using this we see $ \bb $ in \eqref{:81b} satisfies \eqref{:71g} with $ \rho_s=0 $. 
We thus verify the condition \As{A6} from \lref{l:71}. 
We have checked \As{A1}--\As{A6}. 
Then from \tref{l:31} we deduce that the labeled process $\mathbf{X}=(X^i)_{i\in\Z}$ solves the ISDE 
\begin{align}& \label{:81c} 
 dX_t^i = dB_t^i + \frac{\beta }{2} \lim_{r\to\infty }\sum_{j\not= i,\, |X_t^i - X_t^j |< r } 
\frac{1}{X_t^i - X_t^j } dt \quad (i\in\mathbb{Z}) 
\end{align}
and satisfies \As{$\mu $-AC} and \As{NBJ}.

Tsai \cite{tsai.14} proved that the ISDE \eqref{:81c} has a pathwise unique strong solution 
for each $ \beta \ge 1 $ \cite{tsai.14}. 
Hence we can apply \tref{l:32} and obtain \eqref{:70a} for $ \mu _{\mathrm{sin},\beta }$ with $ \beta =1,2,4 $. 
For $ \beta \not=1,2,4$, neither the quasi-Gibbs property nor 
the existence of the logarithmic derivatives of 
$ \mu _{\mathrm{sin},\beta } $ have been proved. 
Hence we can not apply \tref{l:31} and \tref{l:32} for this regime at present. 

If $ \beta = 2 $, then an algebraic construction of the stochastic dynamics 
associated with the upper Dirichlet form $ (\E^{\Os},\mathcal{D} ^{\Os})$ is known \cite{o-t.sm}. 
The distribution of the dynamics are determined by the space-time correlation functions, which is explicitly given by a concrete determinantal kernel. Because of the identity 
$ (\E^{\Os},\mathcal{D} ^{\Os})=(\E^{\La},\mathcal{D} ^{\La})$ in \tref{l:32}, the same holds for the stochastic dynamics associated with the lower Dirichlet form $ (\E^{\La},\mathcal{D} ^{\La})$. 

\subsection{Gibbs measures with Ruelle-class potential}\label{s:85}
Let $ \sS = \Rd $ with $ d \in \N $. 
Let $ \Phi = 0$ and we consider ISDE \eqref{:12a}. 
Assume that $ \Psi $ is smooth outside the origin and is 
a Ruelle-class potential. 
That is, $ \Psi $ is super-stable and regular in the sense of Ruelle \cite{ruelle.2}. Here we say $ \Psi $ is regular if there exists a positive decreasing function $ \map{\psi }{\mathbb{R}^+}{\mathbb{R}}$ and $ R_0 $ such that 
\begin{align}&\label{:75a}
\Psi (x) \ge - \psi (|x|) \quad \text{ for all } x ,
 \quad 
\Psi (x) \le \psi (|x|) \quad \text{ for all } 
 |x| \ge R_0 ,
\\ &\notag 
\int_0^{\infty} \psi (t)\, t^{d-1}dt < \infty 
.\end{align}
Let $ \mupsi $ be canonical Gibbs measures with interaction 
$ \Psi $ satisfying \As{A2} and \As{A4}. %
We do not a priori assume the translation invariance of $ \mupsi $. 

Let $ \rho _x^1 $ be the $ 1 $-point correlation function 
of the reduced Palm measure of $ \mupsi $ conditioned at $ x $. 
Suppose that for each $ \p \in \N $ 
\begin{align} &\label{:75c}
\int_{\sS _k }\int_{|x-y|\ge 1/\p } |\nabla \Psi (x-y) |\rho _x^1 (y)dxdy < \infty 
.\end{align}

For the non-collision property of tagged particles we assume the following. Suppose that $ d \ge 2 $ or that $ d=1$ with $ \Psi $ is 
sufficiently repulsive at the origin in the following sense (see \cite{inu}). 
There exist a positive constant $\Ct \label{;27E} $ and a 
 positive function $ \map{h}{(0,\infty)}{[0,\infty]}$ satisfying that 
\begin{align}&\label{:75e} 
\int_{0 < t \le \cref{;27E}} \frac{1}{h(t)} dt = \infty 
,\\& \notag 
\rho ^m(x_1,\ldots,x_m) \le h (|x_i-x_j|) 
\quad \text{ for all } x_i\not=x_j 
.\end{align}

From the DLR equation, $ \mu _{\Psi }$ satisfies \As{A1}. 
We have \As{A2} and \As{A4} hold by assumption. 
\As{A3} follows from \eqref{:75e} (see \cite{inu} for proof). 
\As{A5} is satisfied with the logarithmic derivative given by
\begin{align}&\notag 
 \dpsi (x,\mathsf{y}) = - 
\beta \sum ^{\infty}_{j=1} \nabla \Psi (x - y^j )
\quad \text{ in }L_{\mathrm{loc}}^p (\mu_{\Psi } ^{[1]})
.\end{align}
We can readily check the condition \As{A6} by \lref{l:71}. 
By \eqref{:75c} $ \chi _{s}$ satisfies \As{A6}, while the construction of 
$ \upsilon _{\p } $ depends on the singularity of $ \Psi $ at the origin. 
In particular, if $ \Psi $ is continuous at the origin, we can dispense with $ \upsilon _{\p } $. 
Hence we can apply \tref{l:31} to obtain that the labeled process 
$\mathbf{X}=(X^i)_{i\in\N }$ given by the lower Dirichlet form 
satisfies the conditions \As{$\mu $-AC} and \As{NBJ}, and solves ISDE 
\begin{align} &\notag 
dX^i_t = dB^i_t - 
\frac{\beta }{2} \sum ^{\infty}_{j=1,\, j\ne i} 
\nabla \Psi (X_t^i - X_t^j ) dt \quad (i\in \N )
.\end{align}

If $ \Psi \in C_0^3(\R^d )$ in addition to the conditions above and $ \mu $ is a unique grand canonical Gibbs measure for $ \beta >0 $, Lang proved that the associated ISDE has a pathwise unique strong solution. 
The assumptions of \tref{l:32} are then fulfilled if $ d \ge 2$. Hence we obtain \eqref{:70a} for this class.

\subsection{Bessel$_{2,\alpha}$ interacting Brownian motion} \label{s:83}

Let $ d = 1 $ and $ \sS = [0,\infty)$. Let $ 1 \le \alpha < \infty $. 
Let $ \mu _{\mathrm{Be},\alpha }$ be 
the Bessel$_{2,\alpha }$ random point field. 
By definition $ \mu _{\mathrm{Be},\alpha }$ is a determinantal random point field 
whose $ n $-point correlation function 
$ \rho _{\mathrm{Be},\alpha }^n $ 
with respect to the Lebesgue measure on $[0,\infty) $ is given by 
\begin{align}\notag 
&
 \rho _{\mathrm{Be},\alpha }^n (x_1,\ldots ,x_n) = 
 \det [ \mathsf{K}_{\mathrm{Be},\alpha }(x_i,x_j)] \ijn 
.\end{align}
Here $ \mathsf{K}_{\mathrm{Be},\alpha } $ is a continuous kernel given by 
\begin{align}\notag 
 &
\mathsf{K}_{\mathrm{Be},\alpha }(x,y) = 
 \frac{J_{\alpha } (\sqrt{x}) \sqrt{y} J_{\alpha }' (\sqrt{y}) - 
 \sqrt{x} J_{\alpha }' (\sqrt{x}) J_{\alpha }(\sqrt{y}) }{2(x-y)}
\quad \text{ for }x\not=y 
.\end{align}
The Bessel$_{2,\alpha}$ random point fields $ \mu _{\mathrm{Be},\alpha }$
satisfy \As{A1}--\As{A3} \cite{h-o.bes,o.col}. \As{A4} is obvious. 
In \cite{h-o.bes}, it was proved that 
\As{A5} is satisfied with the logarithmic derivative given by
\begin{align}&\notag 
 \bb (x,\mathsf{y}) = \frac{\alpha}{2x} + \sum_{i=1}^{\infty} \frac{1}{x-y_i}
\quad \text{ in }L_{\mathrm{loc}}^p (\mub ^{[1]})
.\end{align}
We can readily verify the condition \As{A6} by \lref{l:71}. 
Hence the assumptions of \tref{l:31} are fulfilled. 
The labeled process $\mathbf{X}=(X^i)_{i\in\N }$ 
given by the lower Dirichlet form thus 
satisfies the conditions \As{$\mu $-AC} and \As{NBJ}, and 
solves the ISDE: 
\begin{align} &\notag 
dX_t^i = dB_t^i + \{ \frac{\alpha }{2X_t^i } + 
\sum _{ j\not = i }^{\infty}
\frac{1}{X_t^i - X_t^j} \} dt \quad (i \in \N )
.\end{align}
%

%

\subsection{Ginibre interacting Brownian motion}\label{s:84}

Let $ d=2 $ and $ \sS = \mathbb{R}^2$. Let $ \beta = 2 $. Let $ \mug $ be the Ginibre random point field. 
By definition $ \mug $ is a random point field on $ \mathbb{R}^2$ whose 
$ n $-point correlation function with respect to the Lebesgue measure is given by 
\begin{align}&\notag 
\rgn ^{n}(x_1,\ldots,x_n) = \det [\kg (x_i,x_j)]\ijn 
,\end{align} 
where $ \map{\kg }{\R ^{2}\ts \R ^{2}}{\mathbb{C}}$ is the kernel defined by 
\begin{align}&\notag 
\kg (x,y) = \frac{1}{\pi } 
e^{- \half\{|x|^{2}+|y|^{2}\}} \cdot 
e^{x \bar{y}}
.\end{align}
Here we identify $ \R ^{2} $ as $ \mathbb{C}$ by the obvious correspondence 
$ \R ^{2} \ni x=(x_1,x_2)\mapsto x_1 + \mathrm{i} x_2 \in \mathbb{C}$, and 
$ \bar{y}=y_1-\mathrm{i} y_2 $ is the complex conjugate in this identification, where $ \mathrm{i} = \sqrt{-1}$. 
The random point field $ \mug $ satisfies \As{A1}--\As{A3} \cite{o.isde,o.rm,o.col}. 
Because the one-point correlation function is constant, $ \mug $ clearly satisfies \As{A4}. 
The logarithmic derivative is given by
\begin{align}&\label{:84c} 
 \bb (x,\mathsf{y}) = 
\lim_{r\to\infty } \sum_{|x-y_j|<r } 
\frac{x-y_j}{|x-y_j|^{2}}
\quad \text{ in }L_{\mathrm{loc}}^p (\mug^{[1]})
,\\
&\label{:84d}
 \bb (x,\mathsf{y}) = 
-x +\lim_{r\to\infty } \sum_{|y_j|<r } 
\frac{x-y_j}{|x-y_j|^{2}}
\quad \text{ in }L_{\mathrm{loc}}^p (\mug^{[1]})
.\end{align}
It is known that \eqref{:84c} and \eqref{:84d} define the same logarithmic derivative \cite{o.isde}. 
We have thus checked \As{A5}. 
We can readily check the condition \As{A6} by \lref{l:71}. 
Hence the assumptions of \tref{l:31} are fulfilled.

The labeled process $\mathbf{X}=(X^i)_{i\in\N }$
given by the lower Dirichlet form 
satisfies the conditions \As{$\mu $-AC} and \As{NBJ} and 
solves both ISDEs below: 
\begin{align}&\notag 
dX_t^i = dB_t^i + \lim_{r\to\infty } 
\sum_{j\not=i,\, |X_t^i-X_t^j|<r } 
\frac{X_t^i-X_t^j}{|X_t^i-X_t^j|^{2}} dt \quad (i \in\N )
,\\\notag 
& 
dX_t^i = dB_t^i - X_t^i dt + \lim_{r\to\infty }
\sum_{j\not=i,\, |X_t^j|<r } 
\frac{X_t^i-X_t^j}{|X_t^i-X_t^j|^{2}} dt \quad (i \in\N )
.\end{align}


\section{Concluding remarks and questions}\label{s:9}

\noindent 
{\bf 1. } We have proved that the two natural Dirichlet forms 
$ (\E^{\La},\mathcal{D} ^{\La})$ and $(\E^{\Os},\mathcal{D} ^{\Os}) $ 
are equal under the assumptions in \tref{l:32}. 
The most important condition for this is the non-explosion property of each tagged particle that follows from \As{A4}. Indeed, this condition controls the effect of boundary $ \partial \SR $ as $ R \to \infty$. We have an example of non-coincidence when tagged particles explode. We then conjecture that non-explosion is a necessary and sufficient condition of the coincidence of 
the upper and the lower Dirichlet forms. 

\noindent {\bf Question 1. } 
Can one prove that the upper and the lower Dirichlet forms coincide with each other if and only if each tagged particle does not explode? 
\smallskip 

\noindent 
{\bf 2. }We can naturally formulate the same problem for non-local Dirichlet forms. In particular, the case such that the associated Markov processes have big jump would be interesting. 

\smallskip 

\noindent 
{\bf 3. }In \cite{kuwae,kuwae-s}, 
the uniqueness of the Silverstein extension of Dirichlet forms was studied. 
In particular, it was proved that the Silverstein extension is unique when the Dirichlet form is quasi-regular and equipped with a suitable exhaustion function with bounded energy measure \cite[Theorem 5.1, Theorem 6.1]{kuwae}. 
Our result \eqref{:32c} however can not be derived from this because we do not a priori know whether the lower Dirichlet form 
$ (\E^{\La},\mathcal{D} ^{\La})$ is a Silverstein extension of the upper Dirichlet form $(\E^{\Os},\mathcal{D} ^{\Os}) $. 
As a corollary of \tref{l:32}, we see that 
$ (\E^{\La},\mathcal{D} ^{\La})$ is the Silverstein extension of the upper Dirichlet form $(\E^{\Os},\mathcal{D} ^{\Os}) $ because these two forms are equal. 

\smallskip 

\noindent {\bf 4. }In \cite{Tk92}, Takeda proved the uniqueness of Markovian extension of Dirichlet forms on distorted Brownian motion in a domain in $ \Rd $ (also called a generalized Schr\"{o}dinger operator). 
We refer to \cite[Chapter 3.3]{fot.2} for the Markovian extension of Dirichlet forms. 
This class of Dirichlet forms is a finite-dimensional counterpart of the Dirichlet forms in the present paper. Hence it is natural to discuss the uniqueness of the Markovian extension of the upper Dirichlet form $(\E^{\Os},\mathcal{D} ^{\Os}) $.

\noindent {\bf Question 2. }What is the sufficient condition for the uniqueness of the Markovian extension of the upper Dirichlet form $(\E^{\Os},\mathcal{D} ^{\Os}) $? Is it same as the condition in Question 1? 
%
%
%
%

\noindent {\em Acknowledgements: } 
{Y.K. is supported by JSPS KAKENHI Grant Numbers15J03091. H.O. is supported in part by Grant-in-Aid for Scientific Research (S), No.\!\! 16H06338; Grant-in-Aid for Challenging Exploratory Research No.\!\! ：16K13764 from the Japan Society for the Promotion of Science. H.T. is supported in part by Grants-in-Aid for Scientific Research (C), No.\! 15K04010; Scientific Research (B), No. 19H01793  from the Japan Society for the Promotion of Science. }

\begin{footnotesize}
\begin{description}
\item{1} \qquad
Faculty of Mathematics, 
Kyushu University, \\
Fukuoka 819-0395, Japan.
E-mail: y-kawamoto@math.kyushu-u.ac.jp 
\item{2} \qquad
Faculty of Mathematics, 
Kyushu University, \\
Fukuoka 819-0395, Japan.
E-mail: osada@math.kyushu-u.ac.jp 
\item{3} \qquad
Department of Mathematics and Informatics,
Faculty of Science, Chiba University, \\
1-33 Yayoi-cho, Inage-ku, Chiba 263-8522, Japan. 
E-mail: tanemura@math.s.chiba-u.ac.jp

\end{description}
\end{footnotesize}

\end{document}